\documentclass[twoside,12pt]{article}

\usepackage{amsmath}
\usepackage{amsthm}
\usepackage{amssymb}
\usepackage{latexsym}
\usepackage{array}
\usepackage{longtable}

\usepackage{tikz}
\usetikzlibrary{matrix,arrows}

\usepackage{times}
\usepackage{graphicx}

\usepackage{fancyhdr}
\usepackage{varioref}

\newtheorem{thm}{\bf Theorem}[]
\newtheorem{prop}[thm]{\bf Proposition}
\newtheorem{cor}[thm]{\bf Corollary}
\newtheorem{lem}[thm]{\bf Lemma}

\newtheorem{rem}[thm]{\bf Remark}

\numberwithin{equation}{section}

\newcommand{\atopa}[2]{\genfrac{}{}{0pt}{}{#1}{#2}}

\def\qed{\phantom{a}{\itshape Q.E.D.}}

\def\proof{{\noindent \itshape Proof.\ \ }}

\def\ve {\varepsilon}
\def\d {{\rm d}}
\def\half {{\frac 12}}

\def\res{{\rm Res}}
\def\GH{{\Gamma\backslash\mathfrak H}}

\def\fbar{\hskip 1.25pt{\overline{\hskip -1.25pt  f}}}
\def\pt{\hskip 1.0pt}
\newcount\hours
\newcount\minutes
\def \SetTime{\hours=\time
  \global\divide\hours by 60
  \minutes=\hours
  \multiply\minutes by 60
  \advance\minutes by-\time
  \global\multiply\minutes by-1}
\SetTime
\def \now{\number\hours:\ifnum\minutes<10 0\fi\number\minutes}

\def\today{\ifcase\month\or
  January\or February\or March\or April\or May\or June\or
 July\or August\or September \or October\or November\or December\fi
  \space\number\day,\space\number\year}

\normalsize

\title{Designed pseudo-Laplacians}
\author{E. Bombieri \and P. Garrett}

\begin{document}

\maketitle


\setcounter{section}{0}
\setcounter{subsection}{0}
\setcounter{footnote}{0}
\setcounter{page}{1}

 \pagestyle{fancy}
 \fancyhead{}
 \fancyhead[LE]{\thepage}
 \fancyhead[RO]{\thepage}
 \fancyhead[CO]{Bombieri and Garrett} 
 \fancyhead[CE]{Designed pseudo-Laplacians} 
 \fancyfoot[C]{}
\thispagestyle{plain}

\begin{abstract}
We elaborate and make rigorous various speculations about the
implications of spectral properties of self-adjoint operators on spaces of
automorphic forms for location of zeros of $L$-functions. Some of these
ideas arose in work of Colin de Verdi\`ere, Lax-Phillips, and Hejhal,
from the late 1970s and early 1980s, not to mention semi-apocryphal
attributions to P\'olya and Hilbert. For example, given a complex
quadratic extension $k$ of $\mathbb Q$, we give a natural
self-adjoint extension of a restriction of the invariant Laplacian on
the modular curve whose discrete spectrum, if any, consists of values $s(s-1)$ for
zeros $s$ of $\zeta_k(s)$. Unfortunately, there seems to be no reason
for this discrete spectrum to be large. In fact, Montgomery's pair
correlation, and the behavior of $\zeta(1+it)$, imply that at most
$94\%$ of zeros of $\zeta(s)$ can appear in this discrete
spectrum. Less naively, some preliminary positive results about the
dynamics of zeros do follow from these considerations.
\end{abstract}


\def\Eis{{\mathfrak E}}

{\parskip=-1pt

\noindent 1. Introduction \dotfill 2

\noindent 2. Friedrichs extensions \dotfill 9

2.1 Friedrichs self-adjoint extensions \dotfill 9

2.1 Friedrichs extensions of restrictions \dotfill 10

\noindent 3. Eisenstein-Sobolev spaces \dotfill 15

3.1 Eisenstein series \dotfill 16

3.2 Heegner points and Eisenstein series \dotfill 16

3.3 Spectral decomposition and spectral synthesis \dotfill 18

3.4 Eisenstein-Sobolev spaces ${\mathfrak E}^r$  \dotfill 19

3.5 Automorphic Dirac delta distributions \dotfill 20

3.6 Eisenstein-Heegner distributions \dotfill 21

3.7 Solving $(-\Delta-\lambda_w)u=\theta$ \dotfill 22

3.8 Constant-term distributions $\eta_a$ \dotfill 22

3.7 Solving $(-\Delta-\lambda_w)u=\eta_a$ 

\noindent 4. Pseudo-Laplacians on non-cuspidal automorphic  spectrum
\dotfill 24

4.1 Necessary condition for discrete spectrum \dotfill 24

4.2 Extensions of restrictions to non-cuspidal pseudo-cuspforms
\dotfill 26

4.3 Exotic eigenfunction expansions \dotfill 31

4.4 Extensions of restrictions by Heegner and constant-term
constraints \dotfill 36

\noindent 5. Eigenfunctions of pseudo-Laplacians \dotfill 36

5.1 Determination of eigenfunctions \dotfill 36

5.2 Computing $\eta_a(v_{w,a})$ for $a>1$ and $\Re(w)>1/2$ \dotfill 38

5.3 Computing $\theta_d(v_{w,a})$ for $a\gg_\theta 1$ and
$\Re(w)>1/2$ \dotfill 38

5.4 Computing $\eta_a(u^{}_{\theta,w})$ for $a\gg_\theta 1$ and
$\Re(w)>1/2$ \dotfill 40

5.5 Computing $\theta(u^{}_{\theta,w})$ for $a>1$ and
$\Re(w)>1/2$ \dotfill 41

5.6 Rewriting the determinant condition \dotfill 41

5.7 Meromorphic continuation and location of zeros \dotfill 42

5.8 An important remark \dotfill 45

5.9 Computing $\theta E_s$ in a special case \dotfill 46

\noindent 6. Meromorphic continuations of spectral expansions \dotfill 47

6.1 Vector-valued integrals \dotfill 48 

6.2 Holomorphic vector-valued functions \dotfill 50

6.3 Spaces $M$ of moderate-growth functions \dotfill 51

6.4 Pre-trace formula and $\Eis^{1+\ve}\subset M$ \dotfill 52

6.5 Meromorphic continuation of spectral integrals \dotfill 53

6.6 Spectral corollaries of meromorphic continuation \dotfill 59

6.7 Remarks on appealing but incorrect arguments \dotfill 62

\noindent 7. Spacing of spectral parameters  \dotfill 62

7.1 Exotic eigenfunction expansions of distributions  \dotfill 62

7.2 Exotic eigenfunction expressions for solutions of differential
equations  \dotfill 63

7.3 Interleaving  \dotfill 66

7.4 Spacing of spectral parameters   \dotfill 68

7.5 Juxtaposition with pair-correlation  \dotfill 72

\noindent 8. Spacing of zeros of $\zeta_k(s)$  \dotfill 72

8.1 Exotic eigenfunction expansions and interleaving \dotfill 73

8.2 Dependence of eigenvalues on cut-off height \dotfill 76

8.3 Sample unconditional results on spacing of zeros \dotfill 78

\noindent References \dotfill 79

}


\section{Introduction}

\vskip -10pt

\def\uhat{{\widehat{u}}}
\def\afc{{\mathrm a\!f\!c}}
\def\Ttil{{\widetilde{T}}}
\def\Stil{{\widetilde{S}}}
\def\B{{\mathfrak B}}
\def\inc{{\mathrm{inc}}}
\def\nc{{\mathrm{nc}}}


First, a simple example of the content of theorem
\ref{solvability-necessary-condition}: there is a Hilbert space
$\Eis^0$ (described just below) of
automorphic forms on $SL_2(\mathbb Z)\backslash\mathfrak H$ such that,
for every complex quadratic field 
extension $k$ of $\mathbb Q$, there is a self-adjoint
extension $\widetilde{S}_k$ of a certain restriction $S_k$ of
$-\Delta$, with the invariant Laplace-Beltrami operator
$\Delta=y^2(\partial^2/\partial x^2+\partial^2/\partial y^2)$, such that the
eigenvalues $\lambda_s=s(1-s)>1/4$ of $\widetilde{S}_k$ on $\Eis^0$, {\it
if any}, occur only for $s$ a zero of the Dedekind zeta function
$\zeta_k(s)$. Something like this is suggested by parts of
\cite{CdV-pseudo}, although perhaps often misunderstood as heuristics
rather than as potentially rigorous arguments.
In any case, this is not what one might hope: there is no assertion of
existence of {\it any} eigenvalues of $\widetilde{S}_k$, and
there is no assertion that every zero of $\zeta_k(s)$ gives an
eigenvalue. Corollary \ref{ninety-four} shows that, assuming Montgomery's pair
correlation conjecture \cite{Montgomery-pair}, at most
$94\%$ of the zeros $s$ of $\zeta(s)$ can occur as parameters for
eigenvalues $\lambda_s$ of any particular $\widetilde{S}_k$,
suggesting that perhaps {\it none} do.


Eigenvalues of self-adjoint operators are {\it real}. Thus, when
parametrized as $\lambda_s=s(1-s)\in\mathbb R$, either
$\Re(s)=\frac12$ or $s\in \mathbb R$. Thus, as apocryphally suggested
by G. Polya and by D. Hilbert, one might imagine proving the
Riemann Hypothesis by finding a self-adjoint operator such that,
for every non-trivial zero $s$ of $\zeta(s)$, $\lambda_s$ is an
eigenvalue. (See \cite{O}.)


In \cite{Lax-Phillips} an argument is sketched for discrete
decomposition (i.e., pure point spectrum) of $L^2_a(\Gamma\backslash\mathfrak H)$, 
the latter space defined to be $L^2$ automorphic forms with constant term
vanishing at height $y\ge a>1$, with respect to the Friedrichs self-adjoint
extension $\Stil_\Theta$ of the restriction $S_\Theta$ of $-\Delta$ to
$C^\infty_c(\Gamma\backslash\mathfrak H)\cap
L^2_a(\Gamma\backslash\mathfrak H)$.
The first surprise is that part of the {\it continuous} spectrum of the
original invariant Laplace-Beltrami operator lying inside $L^2_a(\GH)$, 
consisting of suitable wave-packets of Eisenstein series $E_s$, becomes discrete.
Even more surprising, the new eigenfunctions with eigenvalues
$\lambda>1/4$ are certain 
{\it truncated Eisenstein series} $\wedge^a E_s$, namely, for $s$ such
that $a^s+c_s a^{1-s}=0$, where $a^s+c_sa^{1-s}$ is the constant term
of $E_s$ evaluated at height $a$, seemingly contradicting elliptic regularity.
Evidently eigenfunctions of the Friedrichs extension of a restriction of an 
elliptic operator can fail to be {\itshape smooth}.
Thus, $a^s+c_sa^{1-s}=0$ if and only if
$\lambda_s=s(1-s)$ is an eigenvalue of $\widetilde{S}_\Theta$. By
self-adjointness, $a^s+c_sa^{1-s}=0$ only for $s\in \frac12+i\mathbb
R$ (or in $[0,1]$). This idea appears in \cite{Lax-Phillips}, pages
204-6.

In that context, the document \cite{Haas} was provocative: attempting
numerical determination of eigenvalues for $-\Delta$ on $SL_2(\mathbb
Z)\backslash\mathfrak H$, the list of
spectral parameters $s$ for eigenvalues $\lambda_s=s(1-s)$
was observed by H. Stark and D. Hejhal to include low-lying zeros of
$\zeta(s)$ and of $L(s,\chi_{-3})$ with the Dirichlet $L$-function of
conductor $3$.
\cite{Hejhal}
reported attempted reproduction of the numerical results, with the
zeros of zeta and the $L$-function notably missing from the list of
spectral parameters, and observing that the spurious appearance of
these values in the list of \cite{Haas} was due to a mis-application
of the Henrici collocation method \cite{FHM}.
Hejhal further observed that the solution procedure of
\cite{Haas}
in fact allowed as solution a value $u(z)=G_{\lambda_s}(z,\omega)$ of an
{\it automorphic Green's function} $G_{\lambda_s}(z,z_o)$ as a
solution, with $\omega=e^{2\pi i/6}$.
Automorphic Green's functions and related meromorphic families had
been studied in \cite{Neunh}, \cite{Niebur}, and \cite{Fay}.


That is, the values $\lambda_s=s(1-s)$ obtained by
\cite{Haas} {\it not} belonging to genuine cuspforms were values
$\lambda_s$ fitting into an equation $(-\Delta-\lambda_s)u=\delta^\afc_\omega$ for
$\Gamma$-invariant Dirac $\delta^\afc_\omega$ on $\mathfrak H$,
supported on images 
of $\omega=e^{2\pi i/6}$.  Since $\lambda_s$ appearing in
such an equation are not necessarily eigenvalues of a self-adjoint
operator, they need not be {\it real}.  A claim that $\delta^\afc_\omega$
can be disregarded on the grounds that it has support of measure zero
is incorrect. Nevertheless, there were precedents in
\cite{Lax-Phillips}, \cite{CdV-Eis}, and \cite{CdV-pseudo} for
legitimate reinterpretation of certain inhomogeneous equations as
homogeneous equations. 
For example, in \cite{Lax-Phillips} the exotic
eigenfunctions for $\lambda>1/4$ (certain truncated Eisenstein series) are solutions of
an {\it inhomogeneous} distributional equation $(-\Delta-\lambda_w)u=\eta_a$,
where $\eta_a$ evaluates the constant term of an automorphic form at
height $y=a$. This was a precedent for the suggestion at the end of
\cite{CdV-pseudo} that the inhomogeneous equation
$(-\Delta-\lambda_s)u=\delta^\afc_\omega$ on $\Gamma\backslash\mathfrak
H$ could be converted to a 
homogeneous equation involving a self-adjoint operator, so that the
values $\lambda_s$ would become (genuine) eigenvalues, and thus be {\it
real}.


A complication observed in \cite{CdV-pseudo} for an equation
$(-\Delta-\lambda_s)u=\theta$ on $\Gamma\backslash\mathfrak H$, with
distribution $\theta$ on $\Gamma\backslash\mathfrak H$, is that the Friedrichs
extension $\widetilde{S}_\theta$ of the 
restriction $S$ of $-\Delta$ to the kernel of $\theta$ on
$C^\infty_c(\Gamma\backslash\mathfrak H)$ converts this
inhomogeneous equation to a homogeneous one with an auxiliary
condition, only for $\theta$ in a suitable (global automorphic)
Sobolev space $H^{-1}(\Gamma\backslash\mathfrak H)$. This is the
Hilbert-space dual of 
$H^1(\Gamma\backslash\mathfrak H)$, which is the completion of
$C^\infty_c(\Gamma\backslash\mathfrak H)$ with respect to the Sobolev
$H^1$ norm given by $|f|^2_{H^1} = \big\langle (1-\Delta)f,f\big\rangle_{L^2}$
for $f\in C^\infty_c(\GH)$, discussed in subsection
\ref{Eisenstein-Sobolev}.

To solve equations $(-\Delta-\lambda)u=\theta$ with
distributions $\theta$, it is advantageous to use a {\it spectral
characterization} of (global automorphic) Sobolev spaces.
Recall the spectral decomposition of $L^2(\Gamma\backslash\mathfrak H)$ with
$\Gamma=SL_2(\mathbb Z)$, for example from \cite{Faddeev} or
\cite{Iw}:
first, for $f\in C^\infty_c(\Gamma\backslash\mathfrak H)$, 
$$
f \;=\; \sum_F \langle f,F\rangle\cdot F + \frac{\langle f,1\rangle\cdot 1}{\langle
1,1\rangle}
+ \frac{1}{4\pi i}\int_{\frac12-i\infty}^{\frac12+i\infty} \langle f,E_s\rangle\cdot E_s\;ds
$$
where $F$ runs over an orthonormal basis of $L^2$ cuspforms, with the
$L^2$ hermitian product $\langle,\rangle$ notation abused to denote the
corresponding integral for $f\in C^\infty_c(\Gamma\backslash\mathfrak H)$:
$$
\langle f,E_s\rangle \;=\; \int_{\Gamma\backslash\mathfrak H} f(z)\cdot
\overline{E_s(z)}\;\frac{dx\,dy}{y^2}
$$
Then one proves a Plancherel theorem for test functions $f$:
$$
|f|^2
\;=\;
\sum_F |\langle f,F\rangle|^2 + \frac{|\langle f,1\rangle|^2}{\langle
1,1\rangle}
+ \frac{1}{4\pi i}\int_{\frac12-i\infty}^{\frac12+i\infty} |\langle f,E_s\rangle|^2\;ds
$$
Further, $f\to (s\to \langle f,E_s\rangle)$ has image dense in
the space of $L^2$ functions $g$ on $\frac12+i\mathbb R$ such
that $g(1-s)=c_s\cdot g(s)$.
As with the Plancherel-Fourier transform on $L^2(\mathbb R)$, the
spectral expansion extends to a Plancherel theorem for
$L^2(\Gamma\backslash\mathfrak H)$, with the pairings $\langle f,E_s\rangle$ necessarily
interpreted as isometric extensions. The spectral
synthesis integrals $\int_{(\frac 12)} \langle f,E_s\rangle\cdot E_s\;ds$
converge only in an $L^2$ sense, certainly not necessarily
pointwise. Then, for $r\in
\mathbb R$ we can define (global automorphic) Sobolev spaces $H^r(\Gamma\backslash\mathfrak H)$
to be the completion of $C^\infty_c(\GH)$ with respect to the
$H^r$ norms
$$
|f|^2_{H^r}
\;=\;
\sum_F |\langle f,F\rangle|^2\cdot (1+\lambda_{s_F})^r
+ \frac{|\langle f,1\rangle|^2\cdot (1+\lambda_1)^r}{\langle 1,1\rangle}
$$
$$
\hskip100pt+ \frac{1}{4\pi i}\int_{\frac12-i\infty}^{\frac12+i\infty} |\langle f,E_s\rangle|^2
\cdot (1+\lambda_s)^r
\;ds
$$
where $s_F\in \mathbb C$ gives the eigenvalue $\lambda_{s_F}$ of
cuspform $F$ with respect to $\Delta$, first defined for $f\in
C^\infty_c(\Gamma\backslash\mathfrak H)$.
See subsection \ref{Eisenstein-Sobolev}.


An obstacle appears:
subsection \ref{Sobolev-imbedding} observes that an automorphic 
Dirac $\delta$ is not in $H^{-1}(\Gamma\backslash\mathfrak H)$, but only in
$H^{-1-\ve}(\Gamma\backslash\mathfrak H)$ for every $\ve>0$.
A possible way around this obstacle appears
briefly at the 
end of \cite{CdV-pseudo}, to consider the restriction $\theta$ of
$\delta^\afc_\omega$ to a smaller Hilbert space of automorphic forms.
In \cite{CdV-pseudo}, this smaller space is suggested to be the
orthogonal complement to the discrete spectrum, but it turns out to be
necessary to retain the constants, which appear as square-integrable residues of
Eisenstein series at $s=1$.
For precision, a different description of the restriction
$\theta$ is appropriate, as given in detail in section
\ref{Eisenstein-Sobolev}. First, the suitable analogue of test 
functions here is the space $\Eis^\infty_c$ of pseudo-Eisenstein
series 
$$
\Psi_\varphi(z) \;=\; \sum_{\Gamma_\infty\backslash \Gamma}
\varphi(\Im(\gamma z))
$$
with test-function data $\varphi\in C^\infty_c(0,\infty)$.
The $r$-th Eisenstein-Sobolev space $\Eis^r$ is the completion of
$\Eis^\infty_c$ with respect to the $H^r$-norm. Let $S$ be the
restriction of $-\Delta$ to $\Eis^\infty_c$. Initially, let $\theta$ be
the restriction of $\delta^\afc_{z_o}$ to $\Eis^\infty_c$. 
As we see just below, this restriction $\theta$ has finite
$H^{-1}$-norm, so extends continuously to an element of $\Eis^{-1}$. 
Then take $S_\theta$ to be the further restriction of $S$ to have
domain $\Eis^\infty_c\cap\ker \theta$, which is still dense in
$\Eis^0$, since $\theta\not\in \Eis^0$.
Friedrichs' construction applies to the unbounded operator $S_\theta$
on the Hilbert space $\Eis^0$.
In present terms, for the analogous restriction $\theta$ of a finite
real-linear combination of automorphic $\delta$ such that $\theta
E_s=(\sqrt{d_k}/2)^s\cdot\zeta_k(s)/\zeta(2s)$ for $d_k$ the
(absolute value of) the discriminant of $k$, since $\theta$ is a
restriction of a distribution of compact support on
$\Gamma\backslash\mathfrak H$, $\theta$ is
in $\Eis^{-\infty}=\bigcup_r \Eis^r$.
Its $H^r$-norm would be
$$
|\theta|_{H^r}^2 \;=\;
\frac{|\langle \theta,1\rangle|^2\cdot (1+\lambda_1)^r}{\langle 1,1\rangle}
+ \frac{1}{4\pi i}\int_{(\frac12)}
|(\sqrt{d_k}/2)^s\zeta_k(s)/\zeta(2s)|^2\cdot (1+\lambda_s)^r\;ds 
$$
By the second-moment bound of
\cite{HL} and Landau's bounds on the behavior of $\zeta(s)$ on
the edge of the critical strip, the integral is finite for
$r<-3/4$.
Thus, $\theta\in \Eis^{-\frac{3}{4}-\ve}\subset \Eis^{-1}$, and has a
spectral expansion converging in $\Eis^{-1}$:
$$
\theta
\;=\;
\frac{\langle \theta,1\rangle\cdot 1}{\langle 1,1\rangle}
+ \frac{1}{4\pi i}\int_{(\frac12)} \frac{a_k^{1-s}\zeta_k(1-s)}{\zeta(2s)}\cdot
E_s\;ds
$$
Certainly this does not converge pointwise, but does converge in
$\Eis^{-1}$.
Such expansions allow rigorous
solution of differential equations $(-\Delta-\lambda_w)u=\theta$ {\it by
division}: for $\Re(w)>\frac12$ and $w\not\in \mathbb R$, a solution
in $\Eis^1$ is 
$$
u_{\theta,w}
\;=\;
\frac{\langle \theta,1\rangle\cdot 1}{(\lambda_1-\lambda_w)\cdot
\langle 1,1\rangle}
+ \frac{1}{4\pi i}\int_{(\frac12)} \frac{a_k^{1-s}\zeta_k(s)}{\zeta(2s)\cdot (\lambda_s-\lambda_w)}\cdot
E_s\;ds
$$
Since $\theta\in \Eis^{-1}$, the Friedrichs extension
$\widetilde{S}_\theta$ of the restriction $S_\theta$ of $-\Delta$ to
$\Eis^\infty_c\cap \ker\theta$ behaves as desired: for $u\in V$,
$(\widetilde{S}_\theta-\lambda)u=0$ if and only if
$(-\Delta-\lambda)u=c\cdot \theta$ for 
some constant $c$, and $\theta u=0$. That is, the inhomogeneous
distributional equation $(-\Delta-\lambda)u=c\cdot \theta$ is
equivalent to a homogeneous equation
$(\widetilde{S}_\theta-\lambda)u=0$ together with the auxiliary
condition $\theta u=0$.

At this point, we can make precise sense of one 
speculation from \cite{CdV-pseudo}: with $\theta\in \Eis^{-1}$ the
restriction of $\delta^\afc_\omega$ to $\Eis^1$ (after extending by
continuity), 
theorem \ref{solvability-necessary-condition}
shows that the discrete spectrum
$\lambda_s>1/4$ of $\widetilde{S}_\theta$, {\it if any}, is 
$\lambda_s=s(1-s)$ with $s$ a zero of $\zeta(s)\cdot
L(s,\chi_{-3})$. There is no assurance of existence of
any discrete spectrum.


Spectral theory can be applied in less naive ways:
systematic construction of
natural self-adjoint operators $S\ge 0$ exhibits
meromorphic functions whose zeros $s$ 
are on the critical line $\frac12+i\mathbb R$, and perhaps also
on $[0,1]$. The arguably next-simplest continuation is to consider {\it
two} elements $\eta,\theta\in \Eis^{-1}$, such that 
$$
(\mathbb C\cdot \eta+\mathbb C\cdot \theta)\cap \Eis^0 \;\;=\;\; \{0\} 
$$
with $\eta^c=\eta$ and $\theta^c=\theta$, the restriction
$S_{\eta,\theta}$ of $-\Delta$ to domain
$\widetilde{S}_{\eta,\theta}$. Let $u_{\eta,w}$ and $u_{\theta,w}$ be solutions in
$\Eis^1(\Gamma\backslash\mathfrak H)$ of $(-\Delta-\lambda_w)u=\eta$ and
$(-\Delta-\lambda_w)u=\theta$, respectively. The off-line non-vanishing
argument (see subsection \ref{two-condition}) shows that for
$\Re(w)>\frac12$ and $w\not\in\mathbb R$,  
$$
\langle \eta,u_{\eta,w}\rangle \cdot \langle \theta,u_{\theta,w}\rangle 
-  \langle \eta,u_{\theta,w}\rangle\cdot \langle \theta,u_{\eta,w}\rangle
\;\;\not=\;\; 0
$$
where the pairings are on $\Eis^{-1}\times \Eis^1$. Taking $\eta$ to be the
constant-term evaluation
$$
\eta_a f \;=\; \int_0^1 f(x+ia)\;dx
$$
with $a\gg_\theta 1$ gives (see subsection \ref{two-condition} and the
prior calculations in section 5)
$$
\langle \eta_a,u_{\theta,w}\rangle \;=\; \frac{\theta E_w}{2w-1}
\;=\; \langle \theta,u_{\eta_a,w}\rangle
$$
For complex quadratic $k$ over $\mathbb Q$ with absolute value of
discriminant $d_k$, let $a_k=\sqrt{d_k}/2$, take $\theta\in \Eis^{-1}$ such that $\theta
E_w=a_k^w\zeta_k(w)/\zeta(2w)$. Then
the non-vanishing assertion becomes explicit: in $\Re(w)>\frac12$ and
$w\not\in (\frac12,1]$, for $a\ge a_k$, after some natural
simplifications,

\vbox{
\[
a_k^{1-2w}\cdot (a^{2w-1}+c_w)
\cdot \Big(
\frac{h_k^2/a_k}{\lambda_1-\lambda_w}
+
\frac{1}{4\pi i}\int_{(\frac12)}
\Big|
\frac{\zeta_k(s)}{\zeta(2s)}\Big|^2 \frac{ds}{\lambda_s-\lambda_w}
\Big)
\]
\[\hskip100pt\;-\;
\frac{1}{2w-1}\cdot \frac{\zeta_k(w)^2}{\zeta(2w)^2}
\;\;\not=\;\; 0
\]
}
\noindent
where $h_k$ is the class number of $k$. 
In fact, meromorphic continuation and functional equation
Theorem \ref{F-to-G-functional-eq} and
Corollary \ref{nonvanishing_off_line} show that the latter expression
has zeros only on $\Re(w)=\frac12$ and $[0,1]$. 


Returning to the simplest situation $\theta\in
\Eis^{-1}$, with $\theta E_w=a_k^w\zeta_k(w)/\zeta(2w)$, definitive proof
of presence or absence of discrete spectrum for the operator
$\widetilde{S}_\theta$ seems difficult. Apart from numerical
tests, which suggest that there is {\it no} discrete spectrum, some
clarification can be achieved by 
a subtler application of operator theory and distribution theory, as
follows. By direct computation, the constant term of a solution $u_w$
to the distributional equation $(-\Delta-\lambda_w)u=\theta$ vanishes
at height $y\gg_\theta 1$. Thus, 
such $u_w$ lies inside the corresponding Lax-Phillips space
$L^2_a(\Gamma\backslash\mathfrak H)$ (as above,
and as in section \ref{pseudo-cuspforms}), and can be expanded in terms of 
the (exotic) eigenfunctions 
for $\widetilde{S}_\Theta$ (mostly certain truncated Eisenstein series). 
Expanding (the image of) $\theta$ in terms 
of those eigenfunctions, we find (see subsection \ref{interleaving}) an {\it
interleaving property}: there is at most one spectral parameter $w$ for
an eigenvalue of $\widetilde{S}_\theta$ between any two adjacent zeros
$s$ of $a^s+c_sa^{1-s}$. Arguing as in \cite{Backlund}, the {\it average}
vertical spacing of zeros of $\zeta_k(s)$ (on the critical line or
not) at height $T$ is $\pi/\log T$, the same as that of the spacing of zeros
of $a^s+c_sa^{1-s}$, which bodes well. However, from \cite{T}
(5.17.4) page 112 (in an earlier edition, page 98), for $\log\log T$
large, the argument of $\zeta(s)$ on the edge of the
critical strip is relatively regular, so eventually the spacing of the
zeros of $a^s+c_sa^{1-s}$ is similarly regular. That is, given
$\ve>0$, for $\log\log T$ sufficiently large, the space between
consecutive zeros of $a^s+c_sa^{1-s}$ is between $(1-\ve)\pi/\log T$
and $(1+\ve)\pi/\log T$. Adjusting $a$ slightly, the interleaving
property shows that, given $\ve>0$, for $\log\log T$ sufficiently
large the space between consecutive zeros of $\zeta_{\mathbb Q(\omega)}(s)$ on
the critical line is at least $(1-\ve)\pi/\log T$. This would be in
conflict with the pair correlation conjecture \cite{Montgomery-pair}: for
example, under the Riemann Hypothesis and assuming pair correlation,
such a lower bound would allow at most 94\% of zeros $s$ of $\zeta(s)$
to appear as parameters for eigenvalues $\lambda_s$ (see
corollary \ref{ninety-four}).

Corollary \ref{spacing-corollary} proves an illustrative positive
result about spacing of on-the-line zeros of $\zeta_k(w)$, {\it
without} any assumptions about point spectrum of self-adjoint
operators. Namely, let $t<t'$ be large, and such that $\frac12+it$ and
$\frac12+it'$ are adjacent on-line zeros of $\zeta_k(w)$. Take
$\theta\in \Eis^{-1}$ such that $\theta E_w=a_k^w\zeta_k(w)/\zeta(2w)$.
Suppose that neither $\frac12+it$ nor $\frac12+it'$ is a zero of
\[
J_{\theta,w}
\;=\;
\frac{h_k^2}{(\lambda_1-\lambda_w)\cdot \langle 1,1\rangle}
+
\frac{1}{4\pi i}\int\limits_{(\frac 12)}
\frac{|\theta E_s|^2-|\theta E_w|^2}{\lambda_s-\lambda_w}\;\d s 
\]
Suppose there is a {\it unique} zero
$\frac12+i\tau_o$ of $J_{\theta,\frac12+i\tau}$ between $\frac12+it$
and $\frac12+it'$, and $\frac{\partial}{\partial
\tau}J_{\theta,\frac12+i\tau}>0$. Then $|t'-t|\ge \frac{\pi}{\log
t}\cdot (1+O({1\over \log\log t}))$. That is, in this configuration,
the distance between consecutive zeros of $\zeta_k(w)$ must be {\it at
least} the average.


Analogous discussions with similar proofs apply to a broad class of
self-adjoint operators on spaces of automorphic forms.


\section{Friedrichs extensions}
This section deals with the construction and properties of
self-adjoint Friedrichs extensions of operators on a
complex Hilbert space.
\vskip 0.2in

\subsection{Friedrichs self-adjoint extensions}\label{2.1}  Consider complex Hilbert spaces $V$ 
with inner product $\langle\ ,\ \rangle$ and required
to have a complex-conjugate linear {\itshape conjugation} map
$v\to \overline  v$, with expected properties:
\[
\overline{\overline v} =v,\quad \overline{\alpha\cdot v} 
= \overline{\alpha}\cdot\overline v \qquad (\alpha\in\mathbb C),
\quad \langle v,\overline w\rangle=\langle w,\overline v\rangle.
\]
Spaces of  $L^2$-complex-valued functions on measure spaces, for example,
have natural conjugations given simply by pointwise conjugation of
functions.

Let $S$ be an unbounded, symmetric, densely-defined operator on $V$ with
domain $D$ dense in $V$.  Assume $S$ is {\itshape semi-bounded}, specifically, that
\[
\inf_{x\in D,\,\langle x,x\rangle = 1} \langle Sx,x\rangle \geqslant c >0.
\]
Suppose $D$ is stabilized by conjugation and that $S$ commutes
with the conjugation $v\rightarrow \overline v$.
For $x,\,y\in D$ let us write $\langle x,y\rangle_1=\langle Sx,y\rangle$ 
and let $i:D\longrightarrow V_1$ denote the completion of $D$ with 
respect to the new inner product $\langle x,y\rangle_1$.  The space $V_1$
has a canonical continuous linear map $j:V_1\to V$ extending by
continuity the identity map $D\to D$, because $\langle x,x\rangle_1
\geqslant c \langle x,x\rangle$ for $x\in D$. In fact, $j$ is
injective: $\langle w,iv\rangle_1=\langle jw,Tv\rangle$ for $w\in V_1$ and $v\in
D$, so $jw=0$ implies that $w$ is orthogonal to the image of $D$ in
$V_1$, which is dense.  Whenever possible
we suppress the inclusions $i$ and $j$ from the notation.  We write
$\inc=j\circ i$, in a commutative diagram
\[
\begin{tikzpicture} 
  \matrix (M) [matrix of math nodes, row sep=3.0em, 
column sep=3.0em, text height=1.5ex,text depth=0.0ex]
{ V_1 &  V \\
  D   &  \\ };
\path[>=latex,right hook->] (M-1-1) edge node[auto]{$\scriptstyle j$} (M-1-2);
\path[>=latex,right hook->] (M-2-1) edge node[auto]{$\scriptstyle i$}(M-1-1);
\path[>=latex,right hook->] (M-2-1) edge node[below=0.50ex, right]{\hskip 0.75ex 
$\scriptstyle \inc$}(M-1-2);
\end{tikzpicture}
\]

Now assume $S$ is genuinely unbounded, so that $V_1\ne V$.  Recall that
the {\itshape Friedrichs extension} 
$(S,D){\,\widetilde{\phantom{s}}}$ of the pair 
$(S,D)$ is a new self-adjoint operator $\widetilde S : V_1\longrightarrow V$ 
with a new domain
$\widetilde D \longrightarrow V_1$, extending $S$ in the sense that there is
a diagram
for the composition $\inc=j\circ i$, in which (only) the outer
curvilinear triangle commutes:
\[
\begin{tikzpicture} 
  \matrix (M) [matrix of math nodes, row sep=3.0em, 
column sep=3.0em, text height=1.5ex,text depth=0.0ex]
{\widetilde D  & V_1 & V \\
            D  &     & \\ };
\path[>=latex,->] (M-1-1) edge [bend left] node[above=1.5ex, right]
{\kern-2.5ex $\scriptstyle (S,D)
{\,\widetilde{\phantom{s}}}$\kern 2.5ex} (M-1-3);
\path[>=latex,right hook->] (M-1-1) edge node[below=1.0ex, right]
{\kern-2ex $\scriptstyle \rm{inc}$\kern 2ex} (M-1-2);
\path[>=latex,right hook->] (M-1-2) edge node[below=1.0ex, right]
{\kern-2ex $\scriptstyle \rm{j}$\kern 2ex} (M-1-3);
\path[>=latex,->] (M-2-1) edge node[below=1.0ex, right]{$\scriptstyle (S,D)$} (M-1-3);
\path[>=latex,right hook->] (M-2-1) edge (M-1-1);
\path[>=latex,right hook->] (M-2-1) edge node[above=0.25ex, left]
{$\scriptstyle i$}(M-1-2);
\end{tikzpicture}
\]
\noindent 
For simplicity, usually we shall write only $\widetilde S$ without
mention of the domain $\widetilde D$, but 
the domain $\widetilde D$ is part of the description
of the Friedrichs extension. 

The Friedrichs extension is characterized by its inverse 
$\widetilde S^{-1}$ being an everywhere defined, continuous, self-adjoint operator 
$\widetilde S^{-1} : V \rightarrow V_1$, 
with the $\langle\ ,\ \rangle_1$
topology on $V_1$, with the property
\[
 \langle x,\widetilde S^{-1}y\rangle_1 = \langle jx,y\rangle\qquad
(x\in V_1,\, y\in V)
\]
with $j$ the embedding $j:V_1\rightarrow V$ defined before.
Thus $\widetilde S$ is self-adjoint, $\widetilde S \geqslant S$, and
\[
\inf_{x\in D,\, \langle x,x\rangle = 1} \langle Sx,x\rangle \leqslant
\inf_{x\in \widetilde D,\,\langle x,x\rangle = 1} \langle \widetilde Sx,x\rangle.
\]
When $\widetilde S \ne S$ the spectra of $\widetilde S$
and $S$ can be different.  If the spectrum of $S$ is discrete and $(\lambda_\nu)$
is the associate sequence of eigenvalues, and similarly for $\widetilde S$,
we have $\lambda_\nu \leqslant \widetilde \lambda_\nu$ for all $\nu$.

\vskip 0.2in


\subsection{Friedrichs extensions of restrictions}\label{extns-of-restrns}
Let $V_{-1}$ be the complex-linear dual of $V_1$, with norm
\[
|\mu|_{-1} := \sup_{x\in V_1,\,\langle x, x\rangle \leqslant 1} |\mu(x)|.
\]
Since $V_1$ is a Hilbert space, the norm $|\cdot|_{-1}$ gives an inner
product $\langle\ ,\ \rangle_{-1}$ by polarization, and $V_{-1}$ is a
Hilbert space.  
Using the conjugation map on $V$, let $\Lambda:V\rightarrow V^*$
be the complex-linear isomorphism of $V$ with its complex-linear dual
by means of  $\Lambda(x)(y) =
 \langle y,\overline x\rangle = \langle x,\overline y\rangle$.

The inclusion $j:V_1\rightarrow V$ dualizes to $j^*:V^*\rightarrow
V_1^*=V_{-1}$ 
by means of $(j^*\mu)(x)=\mu(jx)$ for $\mu\in V^*$ and $x\in V_1$. Thus
we have
\[ 
\begin{tikzpicture} 
  \matrix (M) [matrix of math nodes, row sep=3.0em, 
column sep=3.0em, text height=1.5ex,text depth=0.0ex]
{ V_1 & V & V^* & V_{-1} \\ };
\path[>=latex,right hook->] (M-1-1) edge node[above=1.5ex, right] {\kern-1.0ex $\scriptstyle j$\kern 1.0ex} (M-1-2);
\path[>=latex,->] (M-1-2) edge node[above=1.5ex, right] 
{$\kern-1.0ex \scriptstyle \Lambda$\kern 1.0ex} 
node [below=1.5ex, right] {\kern-1.0ex $\scriptstyle \widetilde{\widetilde{\  }}$\kern 1.0ex}(M-1-3);
\path[>=latex,->] (M-1-3) edge node[above=1.5ex, right] {\kern-1.0ex $\scriptstyle j^*$\kern 1.0ex} (M-1-4);
\path[>=latex,->] (M-1-1) edge [bend right] (M-1-4);
\end{tikzpicture}
\]
\noindent Conjugation acts on $V_{-1}$ by
$\overline\lambda(x) = \lambda(\overline x)$.

Define a continuous linear $S^\#:V_1\longrightarrow V_{-1}$, with
$\langle\ ,\ \rangle_1$ and $\langle\ ,\ \rangle_{-1}$ topologies,
respectively, by
\[
S^\# (x)(y) =
\langle x,\overline y\rangle_1 \qquad ( x,\,y \in V_1).
\]
By the Riesz-Fr\'echet theorem, $S^\#$ is a topological isomorphism.

\begin{prop}\label{domains}\ \ The restriction of $S^\#$
to the domain of $\widetilde S$ is $j^*\circ\Lambda\circ \widetilde S$.
The domain of $\widetilde S$ is $\widetilde D = \{x\in V_1\,:\, S^\#x \in
(j^*\circ\Lambda)V\}$.
\end{prop}
\proof \ \ By construction of the Friedrichs extension, its domain is 
$\widetilde D=\widetilde S^{-1} V$.  Thus, for $x=\widetilde S^{-1} x'$
with $x'\in V$ and all $y\in V_1$ we have
\[
\begin{split}
(S^\#x)(y)&=(S^\#\widetilde S^{-1} x')(y) =
\langle\widetilde S^{-1} x',\overline y\rangle_{-1} 
= \langle x,\overline y\rangle\\
& = ((j^*\circ\Lambda)x')(y) 
=((j^*\circ\Lambda\circ \widetilde S)x)(y).
\end{split}
\]
This shows that  $S^\#\widetilde D =(j^*\circ \Lambda\circ\widetilde S) \widetilde D$.

On the other hand, for $S^\#x = (j^*\circ \Lambda)y$ with
$y\in V$ we have, for all $z\in V_1$:
\[
\langle z,\overline x\rangle_1 = (S^\#x)(z) =
((j^*\circ \Lambda)y)(z) = (\Lambda y)(jz) = (jz,\overline y)
=\langle z,\widetilde S^{-1} \overline y\rangle_1.
\]
Therefore, $\overline x=\widetilde S^{-1}\overline y$, proving the
second statement of the proposition. \qed

Let $\Theta$ be
a finite-dimensional subspace of $V_{-1}$ with
\begin{equation}\label{disjointedness}
\Theta \cap (j^*\circ \Lambda)V =\{0\}.
\end{equation}
Since $\Theta$ consists of linear functionals on $V_1$ continuous in the
$\langle\ ,\ \rangle_1$-topology, the simultaneous kernel $\ker \Theta$
is a closed subspace of $V_1$.
\begin{lem}\label{denseness}\ \ 
$D\cap \ker \Theta$ is dense in $V$.
\end{lem}
\proof\ \  This follows from the general fact that for
a continuous inclusion of Hilbert spaces $j:X\longrightarrow Y$, for
$D\subset X$ dense in $Y$, and for a finite-dimensional subspace
$\Theta\subset X^*$ such that $j^*(Y^*)\cap\Theta=\{0\}$,
we have that $D\cap \ker \Theta \subset X$ is dense in $Y$.

For completeness, we recall the simple proof.  Consider first
$\Theta$ of dimension $1$, spanned by $\theta$.  Since
$\theta\notin j^*(Y^*)$,  $\theta$ cannot be continuous
in the $Y$-topology on dense $D$.  This provides
for each $\ve >0$ an element $x_{\ve}\in D$ with $|x_\ve|_Y < \ve$
and $\theta(x_\ve)=1$.  Given $y\in Y$, density of $D$ in $Y$ yields
a sequence $z_n$ in $D$ approaching $y$ in the $Y$-topology.  If
$\theta(z_n)=0$ for infinitely many $n$, there is nothing to prove.  Otherwise,
the sequence $z'_n=z_n-\theta(z_n)\cdot x^{}_{2^{-n}/\theta(z_n)}$ 
is in $\ker \Theta$ and still $z_n' \to y$ in the $Y$-topology, because
\[
|\theta(z_n)\cdot x^{}_{2^{-n}/\theta(z_n)}|_Y < 
|\theta(z_n)|\cdot \frac{2^{-n}}{|\theta(z_n)|} = 2^{-n}\to 0.
\]
Induction on dimension completes the proof. \qed

Let $S_{\Theta}$ be $S$ restricted to the smaller domain $D_{\Theta}:=
D\cap \ker \Theta$ and let $(S_{\Theta},D_{\Theta}){\widetilde{\,\phantom{s}}}$ be 
the Friedrichs extension
associated to $(S_{\Theta},D_{\Theta})$, with domain $\widetilde D_{\Theta}$,
which is indeed dense, by the preceding Lemma \ref{denseness}.
Our next goal is the analogue of Lemma \ref{domains} for
the domain $\widetilde D_{\Theta}$.  In order to do this, we need some preparatory
observations.

The extension 
\[
(S_{\Theta})^\#\,:\, V_1\cap \ker \Theta\longrightarrow
(V_1\cap \ker \Theta)^*
\]
is defined in the same way as $S^\#$, by
\[
(S_{\Theta})^\#(x)(y) = \langle x,y\rangle_1 \qquad
(x,\,y \in V_1\cap \ker \Theta).
\]
Let 
\[
t_\Theta: V_1\cap \ker \Theta\longrightarrow V_1,\quad
t_\Theta^*: (V_1)^*\longrightarrow (V_1\cap \ker \Theta)^*
\]
be the inclusion and its adjoint.
Assume that $\Theta$ is stable under the extension of the conjugation
map to $V_{-1}$.

\begin{lem}\label{S-Theta-bd}\ \ 
$
(S_\Theta)^\# = t_\Theta^*\circ S^\#\circ t_\Theta.
$
\end{lem}
\proof\ \  Lemma \ref{denseness} shows that $D_\Theta$ is dense
in $V$ in the $V$-topology, so formation of the Friedrichs extension as
an unbounded self-adjoint operator (densely defined) on $V$ is
legitimate.  For $x,\,y\in V_1\cap\ker \Theta$ we have
\[
\begin{split}
(t_\Theta^*\circ S^\#\circ t_\Theta)(x)(y) &=
S^\#(x)(y)
=\langle t_\Theta x,t_\Theta \overline {y}\rangle_1
=\langle x,\overline y\rangle_1
=(S_\Theta)^\#(x)(y),
\end{split}
\]
which is the statement of the lemma. \qed

The following re-characterization of the Friedrichs extension of the
restriction is straightforward, but essential.
\begin{thm}\label{domain-S-Theta-tilde}\ \ 
The domain $\widetilde D_\Theta$ of $\widetilde S_\Theta$ is
\[
\begin{split}
\widetilde D_\Theta
=
\{ x\in V_1\cap \ker\Theta\,:\, (S^\#\circ t_\Theta) x\in(j^*\circ\Lambda)V + \Theta\}.
\\
=
\{ x\in V_1\cap \ker\Theta\,:\, S_\Theta^\# x\in (t^*_\Theta\circ
j^*\circ\Lambda)V\}.
\hskip33pt
\end{split}
\]
We have $\widetilde S_\Theta x = y$, with $x\in V_1\cap \ker\Theta$ and
$y\in V$, if and only if \hfill\break
$(S^\#\circ t_\Theta)x=(j^*\circ\Lambda)y +
\theta$ for some $\theta\in\Theta$.
\end{thm}
\proof\ \  The Friedrichs extension $\widetilde S_\Theta$ is characterized by
\[
\langle z,(\widetilde {S_\Theta})^{-1} y\rangle_1 = \langle z,y\rangle
\qquad (z\in D_\Theta,\, y\in V).
\]
Given $S^\#x=(j^*\circ\Lambda)y + \theta$ with $x\in V_1\cap \ker\Theta$,
$y\in V$, and $\theta\in\Theta$, take $z\in D_\Theta$ and compute
\[
\begin{split}
\langle x,\overline z\rangle_1 &=(S^\#x)(z) =((j^*\circ\Lambda)y +
\theta)(z)
= (j^*\overline{y})(z) + \theta(z)
=\langle z,\overline y\rangle + 0 \\
&= \langle y,\widetilde {S_\Theta}^{-1} S\pt \overline z\rangle
=\langle \widetilde {S_\Theta}^{-1} y,S\pt \overline z\rangle 
=\langle \widetilde {S_\Theta}^{-1} y,\overline z\rangle_1,
\end{split}
\]
thus showing that $\widetilde S_\Theta^{-1}x = y$. On the other hand,
by Lemma \ref{S-Theta-bd}, $(S_\Theta)^\#x=y$ if and only if
$(S^\#\circ t_\Theta)x=y+\theta$ for some $\theta\in \ker t_\Theta^*$,
and $\ker t_\Theta^*$ is the closure of $\Theta$ in $V_{-1}$. Since
$\Theta$ is finite-dimensional, this closure is $\Theta$ itself.

%

The second description of the domain of the Friedrichs extension is
immediate from the previous lemma and from the fact that $\Theta$ is
the kernel of $t_\Theta^*$.
\qed

\begin{cor}\label{cor-eigen}
With the hypotheses of the theorem, for $x\in V_1$ with
$(S^\#-\lambda)x = \theta$ with $\theta\in\Theta$,  if $x\in \ker \Theta$
then $x$ is a $\lambda$-eigenfunction for the self-adjoint operator
$\widetilde S_\Theta \geqslant S \geqslant c>0$.  In that case, $\lambda$ is real
and $\lambda \geqslant c>0$.
If $S$ is merely non-negative the same conclusion holds, except for
the weaker inequality $\lambda \geqslant 0$.
\end{cor}

\proof  The first part of the corollary is immediate.  For the second part,
the condition $S \geqslant c>0$ imposed on $S$ at the beginning of  
Subsection \ref{2.1} can be removed, replacing it by non-negativity, by applying the
first conclusion of the corollary to the operator $S+c$ and noting that the 
new eigenvalues are obtained by making a shift by $c$. \qed


In the previous extension-of-restriction scenario, the ambient Hilbert
space did not change. Now we consider a scenario in which the ambient Hilbert space
changes.
Such a situation arose in \cite{Lax-Phillips}, and was used in
\cite{CdV-Eis}. Consider symmetric $S$ with dense domain $D$ on Hilbert
space $V$ with conjugation $v\to \overline v$, and $\langle
Sv,v\rangle\ge \langle v,v\rangle$ for all $v\in D$. 
Let $\Theta$ be
an $S$-stable, conjugation-stable, not necessarily finite-dimensional
subspace of $D$. Let $V^\Theta$ be the orthogonal complement of
$\Theta$ in $V$, with respect to the hermitian inner product on
$V$. Let $D_\Theta=D\cap V^\Theta$. Let $\Theta_{-1}$ be the closure
of $\Theta$ in $V_{-1}$ in the $V_{-1}$ topology. The relevance of the
$S$-stability of $\Theta$ is as expected, namely, that $S$ restricted
to $D_\Theta$ really does map to $V^\Theta$:

\begin{lem}\label{relevance-of-stability}
$S(D_\Theta)\subset V^\Theta$.
\end{lem}

\proof For $v\in D_\Theta$ and $\theta\in \Theta\subset D$, 
\[
\langle Sv,\theta\rangle
=
\langle v,S\theta\rangle
\in
\{\langle v,\theta'\rangle:\theta'\in\Theta\}=\{0\}
\]
giving the indicated inclusion.
\qed

Unlike the previous situation, we must {\itshape assume} that
$D_\Theta$ is dense in $V^\Theta$, and dense in $V_1\cap V^\Theta$ in
the $V_1$ topology.
We prove the requisite density for cases of interest to
us in Lemma \ref{density-by-smooth-cutoffs}.
Let $(V_1)^\Theta$ be the closure of $D_\Theta$ in $V_1$ with respect to
the $V_1$ topology. Let $S_\Theta$ be the restriction
of $S$ to $D_\Theta$, and $S^\#:V_1\to V_{-1}$ by $(S^\#v)(w)=\langle 
v,\overline w\rangle_1$. Let $t_\Theta:V_1^\Theta\to V_1$ be the inclusion,
and $t_\Theta^*:V_{-1}\to (V^\Theta_1)^*$ the adjoint. Let $\widetilde
S_\Theta$ be the Friedrichs extension of $S_\Theta$. Let
$(S_\Theta)^\#:V^\Theta_1\to (V_1^\Theta)^*$ by $((S_\Theta)^\#
v)(w)=\langle v,\overline w\rangle_1$ for $v,w\in V_1^\Theta$. The
present analogue of Lemma \ref{S-Theta-bd} is

\begin{lem}\label{S-Theta-bd-variant}
$(S_\Theta)^\# = t_\Theta^*\circ S^\#\circ t_\Theta$.
\end{lem}

\proof  Identical to that of Lemma \ref{S-Theta-bd}. 
\qed

The analogue of Theorem
\ref{domain-S-Theta-tilde} has a nearly identical form:

\begin{thm}\label{domain-S-Theta-tilde-variant} The domain of
$\widetilde S_\Theta$ is 
\[
\begin{split}
\widetilde D_\Theta
=
\{v\in V^\Theta_1: (S^\#\circ t_\Theta)v\in (j^*\circ \Lambda)V^\Theta+\Theta_{-1}\}
\\
=
\{v\in V^\Theta_1: (S_\Theta)^\# v\in (t_\Theta^*\circ j^*\circ \Lambda)V^\Theta\}
\hskip28pt
\end{split}
\]
We have $\widetilde S_\Theta x=y$ for $x\in
\widetilde D_\Theta$ and $y\in (j^*\circ \Lambda)V^\Theta$ if and only if $(S^\#\circ
t_\Theta)x=y+\theta$ for some $\theta\in \Theta_{-1}$.
\end{thm}

\proof The argument is essentially identical to that of Theorem
\ref{domain-S-Theta-tilde}.
While $\Theta$ gives (continuous) functionals on $V$, via
the hermitian inner product on $V$, the closure $\Theta_{-1}$ in
$V_{-1}$ gives (continuous) linear functionals on $V_1$, via
duality. Because the pairing $V_1\times V_{-1}\to \mathbb C$ extends
the restriction to $D\times D$ of the $V\times V\to \mathbb C$
pairing, $\ker \Theta_{-1}=V_1\cap V^\Theta$. By assumption,
$V^\Theta=V_1\cap V^\Theta$.

As usual, the extension $\widetilde S_\Theta$ is
characterized by  
$
\langle z,(\widetilde {S_\Theta})^{-1} y\rangle_1 = \langle z,y\rangle
$
for $z\in D_\Theta$ and $y\in V^\Theta$.
Given $S^\#x=y + \theta$ with $x\in V_1^\Theta$, $y\in V^\Theta$, and
$\theta\in\Theta_{-1}$, take $z\in D_\Theta$ and compute
\[
\begin{split}
\langle x,\overline z\rangle_1 &=(S^\#x)(z) =((j^*\circ\Lambda)y +
\theta)(z)
= (j^*\overline{y})(z) + \theta(z)
=\langle z,\overline y\rangle + 0 \\
&= \langle y,\widetilde {S_\Theta}^{-1} S\pt \overline z\rangle
=\langle \widetilde {S_\Theta}^{-1} y,S\pt \overline z\rangle 
=\langle \widetilde {S_\Theta}^{-1} y,\overline z\rangle_1,
\end{split}
\]
Thus, $\widetilde S_\Theta^{-1}x = y$. On the other hand,
by Lemma \ref{S-Theta-bd}, $(S_\Theta)^\#x=y$ if and only if
$(S^\#\circ t_\Theta)x=y+\theta$ for some $\theta\in \ker t_\Theta^*$,
and $\ker t_\Theta^*$ is the closure $\Theta_{-1}$ of $\Theta$ in
$V_{-1}$.
The second description of the domain of the Friedrichs extension is
immediate from the previous lemma and from the fact that $\Theta$ is
the kernel of $t_\Theta^*$. \qed


\section{Eisenstein-Sobolev spaces}\label{Eisenstein-Sobolev}
\vskip 0.2in
Let $\Delta$ be the hyperbolic Laplacian
\begin{equation}\label{Laplacian}
\Delta = y^2\Big(\frac{\partial^2}{\partial x^2} 
+ \frac{\partial^2}{\partial y^2}\Big)
\end{equation}
on the upper half-plane $\mathfrak H$,
Let $\Gamma=SL_2(\mathbb Z)$. The standard inner product is
\begin{equation}\label{inner-product}
\langle  f,\overline g\rangle=\int_{\GH}f(z)\overline{g(z)}\,\d\omega
\end{equation}
(the Petersson inner product) with respect to
the hyperbolic area element $\d \omega_z=y^{-2}\d x\,\d y$. 
The hyperbolic area of a fundamental domain of $\GH$ is $\langle
1,1\rangle = \pi/3$.

Let $S$ be $-\Delta$ restricted to have domain
\begin{equation}\label{domain D}
\Eis^\infty_c \;=\; \Big\{ \Psi_{\varphi}(z)=
\sum_{\gamma\in\Gamma_\infty\backslash\Gamma}
\varphi(\Im(\gamma z))\;:\;
\varphi\in C^\infty_c(0,\infty)\Big\}.
\end{equation}
These $\Psi_\varphi$ are {\itshape pseudo-Eisenstein series} (the incomplete
theta series of other authors) with test-function data $\varphi$.
The {\itshape conjugation map} $f\to\fbar$ on $\Eis^\infty_c$ is the expected 
pointwise conjugation.  It commutes with $\Delta$ and $S$, and $\Eis^\infty_c$ is
stable by conjugation. The ambient Hilbert space is the $L^2(\GH)$
completion $\Eis^0$ of $\Eis^\infty_c$.
The operator $S$ is a non-negative
operator, so the previous discussion of Friedrichs extensions applies
to $1+S$, for example, and thereby to $S$ itself, as unbounded
operator on $\Eis^0$, with domain $\Eis^\infty_c$.


\subsection{Eisenstein series}\label{3.1}
Let $\Gamma_\infty$ be the stabilizer of $i\infty$. The Eisenstein
series associated to the cusp at $i\infty$ and $z\in\GH$ is explicitly
given when $\Re(s)>1$ by \begin{equation}\label{Eisenstein}
E_s(z) := 
\sum_{\gamma\in\Gamma_\infty\setminus \Gamma}{\Im(\gamma z)}^s 
\;=\;
\frac 12 \sum_{(c,d)=1}\kern-5pt{}'\kern 5pt \frac {y^s}{|cz+d|^{2s}} 
\end{equation}
and by analytic continuation for general $s$,
where the sum is over coprime integer pairs $c,d$, and the pair $0,0$
is also excluded. 

The Eisenstein series are functions of the two complex variables $z$
and $s$, automorphic in $z$, while $s$ is the {\itshape spectral
parameter}.
The family of functions $z\to E_s(z)$ of $z$ form a
meromorphic automorphic-function-valued function of $s$, with a
simple pole at $s=1$ and infinitely many  
poles for $\Re(s)<\frac 12$.  For $s$ not a pole, $z\to E_s(z)$ is a
real-analytic function of the variable $z\in \GH$.

By the $SL_2(\mathbb R)$-invariance of the Laplacian,
$(-\Delta - \lambda_s)E_s = 0$,
where 
$\lambda_s = s(1-s)$.
The Eisenstein series satisfy the functional equation
$E_s = c^{}_s E_{1-s}$
with
\begin{equation}\label{c_s}
c^{}_s =
\sqrt \pi\,\frac{\Gamma(s-\frac 12)\zeta(2s-1)}{\Gamma(s)\zeta(2s)}
=\frac{\xi(2-2s)}{\xi(2s)},
\end{equation}
where $\xi(s)$ is the completed Riemann zeta function
and we have used 
the functional equation of the zeta function.  This also yields
$c^{}_s\pt c^{}_{1-s}=1$.

Lastly, the residue of the simple pole of the Eisenstein series $E_s$
at $s=1$ is the constant $1/\langle 1,1\rangle = 3/\pi$.


\subsection{Heegner points and Eisenstein series}
It is a well known yet remarkable fact that, for $\Gamma$ an arithmetic group,
values of certain Eisenstein series for $\GH$ have arithmetical significance.
Recall that a {\itshape fundamental discriminant} is a product of relatively prime
factors of the form
\[-4,\ 8,\ -8,\ (-1)^{(p-1)/2}p, \]
where $p$ is an odd prime.  Associated to a fundamental discriminant $d$\,
there is a real, primitive character 
\[ \chi_{d}(n) =\left(\frac d n\right)
\]
where $(d/n)$ is the Kronecker symbol, which enjoys the multiplicativity
\[ \chi_{d}(n)\chi_{d'}(n)=\chi_{dd'}(n) \quad \text{for $d$ and $d'$ coprime}. \]
The absolute value $|d|$ is the {\itshape modulus} of the character. The fundamental 
discriminants are all numbers, positive or negative, of the form $N$ with $N$ square-free 
and $N\equiv 1 \bmod 4$ or of the form $4N$ with
$N\equiv 2\ \text{or}\  3 \bmod 4$.

From now on $d$ will denote a negative fundamental discriminant.
The integral, positive-definite, Lagrange-reduced quadratic
forms of discriminant $d$ are  
\[
Q(x,y):=Ax^2+Bxy+Cy^2
\]
with
\[
d=B^2-4AC,\ \ |B|\le A\le C, \quad 
(\text{and when $A=|B|=1$ then $B=-1$)}.
\footnote{In this case the two quadratic forms $x^2\pm xy +\frac{|d|+1}4 y^2$ are 
equivalent by $(x,y)\to(x,x\mp y)$ and here
we choose the one with $B=-1$ as a representative, the so-called
{\itshape ambiguous form}.}
\]

\noindent The root 
\[
z^{}_Q=\frac {-B+i\sqrt{|d|}}{2A} \in \GH
\]
of the equation $Az^2+Bz+C=0$ is the {\itshape Heegner point} 
associated to that reduced quadratic form. We have
\[
y_Q^{-1}\cdot|m z^{}_Q+n|^2 = \bigg(\frac{\sqrt {|d|}}{2}\bigg)^{-1}Q(n,-m).
\]
From this, one shows that for discriminants $d<-4$ the value
$E_s(z^{}_Q)$
of the Eisenstein series is
\begin{equation}\label{ideal-classes}
E_s(z^{}_Q) = \bigg(\frac{\sqrt {|d|}}2\bigg)^{s}
 \zeta(2s)^{-1}\zeta(s,\mathfrak z^{}_Q)
\end{equation}
where $\mathfrak z^{}_Q$ is the ideal class of the fractional ideal
$[1,z^{}_Q]$ of the imaginary quadratic field $\mathbb Q(\sqrt{d})$.
\footnote{
If the discriminant is $-3$ or $-4$ the associated quadratic forms $Q(x,y)$ are 
the ambiguous form $x^2-xy+y^2$ and $x^2+y^2$.  Besides the obvious 
automorphism $(x,y)\to(-x,-y)$ arising from $-I\in\Gamma_\infty$ they have
 the automorphisms $(x,y)\to (x-y,x)$ and $(x,y)\to (y,-x)$ of order $6$
and $4$, so equations \eqref{ideal-classes} and \eqref{Eisenstein_decomp}
must be corrected by factors $3$ and $2$ in the right-hand side.}
Therefore, summing over the $h(d)$ ideal classes we have
\begin{equation}\label{Eisenstein_decomp}
\sum_{i=1}^{h(d)} E_s(z_{Q_i}) = 
\bigg(\frac{\sqrt {|d|}}2\bigg)^{s}\zeta(2s)^{-1}\zeta(s,\mathbb Q(\sqrt{d}))=
\bigg(\frac{\sqrt {|d|}}{2}\bigg)^{s} \frac{\zeta(s)}{\zeta(2s)}\,L(s,\chi_{d})
\end{equation}
where $z_{Q_i}$ runs over the $h(d)$ Heegner points for the
fundamental discriminant $d$ and where 
$\chi_{d}$ is the quadratic character associated to $\mathbb
Q(\sqrt{d})$.


\subsection{Spectral decomposition and spectral synthesis}
\vskip 0.1in

The spectral theory of the Laplacian on $L^2(\GH)$,
where $\Gamma$ is a discrete subgroup of $SL_2(\mathbb R)$ acting on the upper
half-plane $\mathfrak H$
and of finite covolume, is well understood: see for
example \cite{Faddeev}, or Iwaniec' monograph \cite{Iw}. This 
decomposes $L^2(\GH)$ as 
the direct orthogonal sum of the $L^2$ cuspidal discrete spectrum,
constants ($L^2$ residues of Eisenstein series), and eigen-packets associated to a 
continuous spectrum generated by the Eisenstein series.
Here we only consider $\Gamma=SL_2(\mathbb Z)$, which has
just one family of Eisenstein series giving the continuous 
spectrum, attached to the single cusp $i\infty$.
For $f\in
C^\infty_c(\GH)$, the spectral expansion is
\[
f(z) \;=\; \sum_F \langle f,F\rangle\cdot F(z)
+ \frac{\langle f,1\rangle\cdot 1}{\langle 1,1\rangle}
+ \frac{1}{4\pi i}\int_{(\frac12)} \langle f,E_s\rangle\cdot E_s(z)\;ds
\]
where $F$ runs over an orthonormal basis for $L^2$ cuspforms. As usual,
$\langle f,E_s\rangle$ cannot be the $L^2$ pairing, because
$E_s$ is not in $L^2$, but by standard abuse of notation refers to the
integral $\int_{\GH} f(z)\,E_{1-s}(z)\,\d\omega_z$, which converges
absolutely for automorphic test functions $f$. For test functions $f$,
the right-hand side converges uniformly pointwise in $z$. We have
Plancherel for test functions:
\[
|f|^2_{L^2} \;=\; \sum_F |\langle f,F\rangle|^2
+ \frac{|\langle f,1\rangle|^2}{\langle 1,1\rangle}
+ \frac{1}{4\pi i}\int_{(\frac12)} |\langle f,E_s\rangle|^2\;ds
\]
Extend the spectral expansion to $f\in L^2(\GH)$ by isometry, and
write $\mathcal E f$ for the extension of $f\to (s\to \langle
f,E_s\rangle)$. 
Pseudo-Eisenstein series $f\in \Eis^\infty_c$ with test
function data are orthogonal to cuspforms, so for such
automorphic forms the spectral expansion and Plancherel become
\[
f \;=\; 
\frac{\langle f,1\rangle\cdot 1}{\langle 1,1\rangle}
+ \frac{1}{4\pi i}\int_{(\frac12)} \mathcal E f(s)\cdot E_s\;ds
\]
and
\[
|f|^2_{L^2} \;=\;
\frac{|\langle f,1\rangle|^2}{\langle 1,1\rangle}
+ \frac{1}{4\pi i}\int_{(\frac12)} |\mathcal E f(s)|^2\;ds
\]
The operator $S$, $-\Delta$ restricted to
$\Eis^\infty_c$, is {\itshape symmetric} because the compact support of
pseudo-Eisenstein series with test-function data allows integration by
parts.  For $f\in \Eis^\infty_c$ we have the {\itshape spectral
relation}, intertwining of $S$ and multiplication by $\lambda_s$, 

\[
\begin{split}
\mathcal E(Sf)(s)&=\int_{\GH}(-\Delta)f(z)\cdot E_{1-s}(z) 
\,\d \omega_z
=\int_{\GH}f(z)\cdot (-\Delta)E_{1-s}(z) \,\d \omega_z\\
&=\int_{\GH}f(z)\cdot \lambda_s\, E_{1-s}(z) \,\d \omega_z
=\lambda_s\cdot\mathcal Ef(s).
\end{split}
\]
Thus, for all $0\leqslant\ell\in\mathbb Z$,
\begin{equation}\label{S^ell}
\mathcal E(S^\ell f)(s) = \lambda_s^\ell\cdot\mathcal Ef(s) 
\qquad (\text{for\ \ } f\in D).
\end{equation}


\subsection{Eisenstein-Sobolev spaces $\Eis^r$}
\label{Eisenstein-Sobolev}
\vskip 0.1in

For $r\in \mathbb R$, the $r$-th global automorphic Sobolev space
$H^r(\GH)$ is the completion of $C^\infty_c(\GH)$ with respect to the
$r$-th Sobolev norm, defined on automorphic test functions $f$ by

\vbox{
\[
|f|^2_{H^r} \;=\;
\sum_F |\langle f,F\rangle|^2\cdot (1+\lambda_{s_F})^r
+ \frac{|\langle f,1\rangle|^2\cdot (1+\lambda_1)^r}{\langle 1,1\rangle}
\]
\[+ \frac{1}{4\pi i}\int_{(\frac12)} |\langle f,E_s\rangle|^2\cdot (1+\lambda_s)^r\;ds
\]
}
The $r$-th Eisenstein-Sobolev space is
\[
\Eis^r = \text{completion of $\Eis^\infty_c$ with respect to $|\cdot|_{H^r}$}
\]

Let $X^r=\mathbb C \oplus X^r_o$ where $X^r_o$ is the weighted $L^2$-space of
measurable functions $g$ on $\frac12+i\mathbb R$ such that 
\[
\int_{-\infty}^\infty |g(\frac12+it)|^2\cdot (\frac14+t^2)^r\;dt \;\;<\;\; +\infty
\]
For all $r\in\mathbb R$, Plancherel restricts (for $r\ge 0$) or
extends (for $r\le 0$) to an isometry
$\Eis^r \to X^r$, by $f\to \frac{\langle f,1\rangle}{\langle
1,1\rangle^{\frac12}} \oplus \mathcal Ef$. 
Let
\[
\Eis^\infty = \bigcap_r\, \Eis^r = \lim_r \Eis^r\qquad
\text{and}\qquad
\Eis^{-\infty} = \bigcup_r \, \Eis^r =
\mathop{\rm colim}_r \Eis^r
\]
in the category of locally convex topological spaces.
There is the expected hermitian pairing, 
\[
\langle f,\theta\rangle_{\Eis^r\times \Eis^{-r}}
=
\frac{\langle f,1\rangle\cdot
\overline{\langle \theta,1\rangle}}{\langle 1,1\rangle} 
+ \frac 1{4\pi i} \int_{(\frac 12)}
\mathcal Ef(s)\cdot \overline{\mathcal E\theta(s)}\,\d s
\]
\noindent For all $r\in \mathbb
R$ the map $S:\Eis^\infty_c\to \Eis^\infty_c$ is continuous 
when the domain is given the $\Eis^r$ topology and the
range is given the $\Eis^{r-2}$ topology. Extending by continuity
defines {\itshape $L^2$-differentiation} $\Eis^r\rightarrow \Eis^{r-2}$.


\subsection{Automorphic Dirac delta distributions}\label{Sobolev-imbedding}
\vskip 0.1in
The {\itshape pre-trace formula} (as in \cite{Iw} and elsewhere) is
\[
\sum_{F:|\lambda(F)|\leqslant T} |F(z_0)|^2
+\frac{|\langle F,1\rangle|^2}{\langle 1,1\rangle}
+ \frac 1{4\pi i} \int_{(\frac 12)} |E_s(z_0)|^2\, \d s \;\;\ll_C\;\; T^2
\]\label{standard-estimate}

\noindent for $z_0$ in a fixed compact subset $C \subset \GH$, where $F$ runs over
an orthonormal basis for the cuspidal
spectrum and $SF =\lambda(F)\, F$.  In particular,
dropping the cuspidal part we have
\[
\frac{|\langle F,1\rangle|^2}{\langle 1,1\rangle}
+
\frac 1{4\pi i} \int_{(\frac 12)} |E_s(z_0)|^2\, \d s \ll_C T^2
\]
Integrating by parts, the function
$s\rightarrow E_s(z_0)$ is in the weighted $L^2$ space $X^{-1-\ve}$,
so $E_s(z_0) 
=\mathcal E \theta(s)$ for some $\theta\in \Eis^{-1-\ve}$, for all $\ve > 0$,
and 
\[
\theta =
\frac{\langle \theta,1\rangle\cdot 1}{\langle 1,1\rangle}
+
\frac 1{4\pi i}
 \int_{(\frac 12)} E_{1-s}(z_0)\cdot E_s \,\d s 
\qquad\text{(as an element of $\Eis^{-1-\ve}$)}
\]
Define the non-cuspidal Dirac $\delta$ distribution, or
Eisenstein-Dirac $\delta$ distribution, by
\begin{equation}\label{delta-nc}
\delta_{z_0}^{\rm nc} := \frac 1{\langle 1,1 \rangle} 
+
\frac 1{4\pi i}
 \int_{(\frac 12)} E_{1-s}(z_0)\cdot E_s \,\d s 
\qquad \text{(as an element of $\Eis^{-1-\ve}$)}
\end{equation}
That is, $\mathcal \delta^{\mathrm nc}_{z_o}(s)=E_{1-s}(z_o)$. 
By design, its action on $f\in \Eis^{1+\ve}$ is $\delta^{\rm
nc}_{z_0}f = f(z_0)$. From 
evaluating the hermitian pairing on $\Eis^{1+\ve}\times
\Eis^{-1-\ve}$:
\[
\begin{split}
\langle f,\delta^{\rm nc}_{z_o}\rangle_{\Eis^{1+\ve}\times \Eis^{-1-\ve}}
&=
\frac{\langle f,1\rangle\cdot
\overline{\langle \delta^{\rm nc}_{z_o},1\rangle}}{\langle 1,1\rangle} 
+ \frac 1{4\pi i} \int_{(\frac 12)}
\mathcal Ef(s)\cdot \overline{\mathcal E\delta^{\rm nc}_{z_o}(s)}\,\d s
\\
&=
\frac{\langle f,1\rangle\cdot 1 }{\langle 1,1\rangle} 
+ \frac 1{4\pi i} \int_{(\frac 12)}
\mathcal E f (s)\cdot \overline{E_s(z_o)}\,\d s
\end{split}
\]
At least for $f$ a test-function pseudo-Eisenstein series, this is
$f(z_o)$.
The estimate is uniform for $z_0$ in compact
subsets of $\GH$, so the map $z_0\rightarrow \delta^{\rm nc}_{z_0}$ is a continuous 
$\Eis^{-1-\ve}$-valued function of $z_0$. Thus, for test-function
pseudo-Eisenstein series $f$, by the Cauchy-Schwarz-Bunyakowsky inequality,
\[
\begin{split}
\sup_{z_o\in C} |f(z_o)|
=
\sup_{z_o\in C} |\langle f,\delta^{\rm nc}_{z_o}\rangle|
\le
\sup_{z_o\in C} |f|_{\Eis^{1+\ve}}\cdot |\delta^{\rm
nc}_{z_o}|_{\Eis^{-1-\ve}}
\\
=
|f|_{\Eis^{1+\ve}}\cdot \sup_{z_o\in C}  |\delta^{\rm
nc}_{z_o}|_{\Eis^{-1-\ve}}
\ll_{C,\ve} |f|_{\Eis^{1+\ve}}
\end{split}
\]
That is, the seminorms obtained by taking suprema on compacta are
dominated by the $\Eis^{1+\ve}$ norm. This proves the {\itshape
Sobolev embedding} $\Eis^{1+\ve} \subset C^0(\GH)$, with the latter
topologized by suprema on compact subsets. Thus, $\delta^{\rm
nc}_{z_o}(f)=f(z_o)$ for all $f\in \Eis^{1+\ve}$.

As a corollary of the above argument, we again see the expected
$\mathcal E \delta^\nc_{z_o}(s)=E_{1-s}(z_o)$.

\vskip 0.2in

\subsection{Eisenstein-Heegner distributions}\label{Heegner-distributions}
\vskip 0.1in
For a fundamental discriminant $d<-4$, let $H_d$ be the set of Heegner points in $\GH$
representing the ideal classes of the ring of integers $\mathbb Q(\sqrt d)$, and let
$\theta_d$ be the functional
\begin{equation}\label{theta_d}
\theta_d = \sum_{z\in H_d} \delta_z^{\rm nc} \in \Eis^{-1-\ve}
\end{equation}
for all $\ve >0$. We call this the {\it Eisenstein-Heegner distribution} attached
to the fundamental discriminant $d$. 
The cardinality  $h(d)=|H_d|$ is the class number of the quadratic field 
$\mathbb Q(\sqrt d)$.  Let $\chi_d$ be the quadratic character
attached to $\mathbb Q(\sqrt d)/\mathbb Q$.  This is a primitive character because
$d$ is a fundamental 
discriminant.  For $d < -4$ we have 
\begin{equation}\label{spectr_transf_H}
\mathcal E\theta_d(s) \;=\;
\bigg(\frac{\sqrt {|d|}}{2}\bigg)^s \frac{\zeta(s)}{\zeta(2s)} L(s,\chi_d).
\end{equation}

The Eisenstein-Heegner distributions $\theta_d$ belong to the space
$\Eis^{-\frac 12-\ve}$ for any $\ve>0$, from the Landau bound 
$1/\zeta(1+it)=O(\log T)$ and the deeper second moment bound for 
$\zeta(s)L(s,\chi)$ on the critical line $\Re(s)=\frac 12$.
The $\Eis^{-\frac 12-\ve}$ spectral expansion is (for $d\ne -3$ or $-4$):
\begin{align*}
\theta_d&=\frac{\theta_d(1)\cdot 1}{\langle 1,1\rangle} +\frac {1}{4\pi i}
\int_{(\frac 12)} \theta_d E_{1-s}\cdot E_s\, \d s\\
&=\frac{h(d)}{\langle 1,1\rangle} +\frac {1}{4\pi i}
\int_{(\frac 12)} 
\bigg(\frac{\sqrt{|d|}}{2}\bigg)^{1-s}\,\frac{\zeta(1-s)L(1-s,\chi_d)}{\zeta(2-2s)}
\cdot E_s\,\d s,
\end{align*}
\label{special-spectral-exp-theta}


\subsection{Solving $(-\Delta-\lambda_w)u=\theta$}\label{solving}
Let $\theta$ be a finite real-linear combination of Eisenstein-Heegner
distributions $\theta_d$.
For $\Re(w)>\half$, the equation $(-\Delta-\lambda_w) u =\theta$
has a unique solution $u_{\theta,w}$ in $\Eis^{\frac 32-\ve}$ for every $\ve>0$,  
with spectral expansion obtained directly from that of $\theta$
via \eqref{S^ell}:
\begin{equation}\label{spec-exp-u}
u_{\theta,w}=\frac{\theta(1)\cdot 1}{\langle 1,1\rangle\cdot (\lambda_1-\lambda_w)}
+\frac {1}{4\pi i}\int_{(\frac 12)} 
\theta E_{1-s}\cdot E_s \pt\frac{\d s}{\lambda_s-\lambda_w}.
\end{equation}
\subsection{Constant-term distributions $\eta_a$}
Let $\eta_a$ denote the constant term distribution at height $a>1$ on
$f\in \Eis^\infty_c $, namely:
\begin{equation}\label{eta-a}
\eta_a\pt  f =\int_0^1 f(x+ia)\,\d x.
\end{equation}
This functional is a compactly-supported, real-valued, regular, Borel
measure on $\GH$, so is a continuous functional on $C^0(\GH)$. By the
remark at the end of \ref{Sobolev-imbedding}, there is a continuous
injection $\Eis^{1+\ve}\subset C^0(\GH)$ for every $\ve>0$, so $\eta_a$
restricts to a continuous functional on $\Eis^{1+\ve}$, still denoted
$\eta_a$. Thus, $\eta_a\in \Eis^{-1-\ve}$ for all $\ve>0$.

As with automorphic Dirac $\delta$ and Eisenstein-Dirac $\delta$, we
remove a potential ambiguity 
about correct determination of spectral coefficients $\mathcal
E\eta_a$. We could again use a variant pre-trace formula, but,
instead, we give an argument relevant to subsequent developments:

\begin{prop}\label{remove-ambiguity-C^0} A continuous functional $\mu$
on the Fr\'echet space $C^0(\GH)$ given by a compactly-supported
real-valued regular Borel measure on $\GH$, restricted to a functional
on $\Eis^\infty$, is in $\Eis^{-1-\ve}$ for every $\ve>0$, and
$\mathcal E \mu(s)=\mu(E_{1-s})$. 
\end{prop}

\proof Fix $\ve>0$. We have a Sobolev imbedding 
$H^{1+\ve}(\GH)\subset C^o(\GH)$ for every $\ve>0$. The continuous
dual of $C^o(\GH)$ is exactly compactly-supported regular Borel
measures $\mu$. Thus, $\mu$ has a natural image in
$\Eis^{-1-\ve}=(\Eis^{1+\ve})^*$, since $\Eis^{1+\ve}\subset
H^{1+\ve}(\GH)$. 
Thus, there is a spectral expansion in $\Eis^{-1-\ve}$:
\[
\mu  = \frac{\langle \mu ,1\rangle\cdot 1}{\langle 1,1\rangle}
+\frac{1}{4\pi i}\int_{(\frac 12)} \mathcal E \mu (s)\cdot E_s\;\d s
\]
For any $u\in \Eis^{1+\ve}$, on one hand,
using $\overline{\mu}=\mu$, 
\[
\mu  (u) = \langle u,\mu \rangle_{\Eis^{1+\ve}\times \Eis^{-1-\ve}}
=
\frac{\langle u,1\rangle\cdot\overline{\langle \mu ,1\rangle}}{\langle 1,1\rangle}
+\frac{1}{4\pi i}\int_{(\frac 12)} \mathcal E u(s)\cdot \overline{\mathcal E \mu (s)}\;\d s
\]
On the other hand, the spectral integral for $u$ converges in
$\Eis^{1+\ve}$, and is the limit of compactly supported integrals of
$C^0(\GH)$-valued functions, so
\[
\begin{split}
\mu  (u)
=
\mu \Big(
\frac{\langle u,1\rangle\cdot 1}{\langle 1,1\rangle}
+\frac{1}{4\pi i}\int_{(\frac 12)} \mathcal E u(s)\cdot E_s\;\d s
\Big)
\\
=
\mu \Big(
\lim_{T\to\infty}
\frac{\langle u,1\rangle\cdot 1}{\langle 1,1\rangle}
+\frac{1}{4\pi i}\int_{|\Im(s)|\le T} \mathcal E u(s)\cdot E_s\;\d s
\Big)
\\
=
\lim_{T\to\infty}
\frac{\langle u,1\rangle\cdot \mu (1)}{\langle 1,1\rangle}
+\frac{1}{4\pi i}\int_{|\Im(s)|\le T} \mathcal E u(s)\cdot \mu
E_s\;\d s
\\
=
\frac{\langle u,1\rangle\cdot \mu (1)}{\langle 1,1\rangle}
+\frac{1}{4\pi i}\int_{(\frac12)} \mathcal E u(s)\cdot \mu  E_s\;\d s
\end{split}
\]
because the continuous functional $\mu $ passes inside the
compactly-supported $C^0(\GH)$-valued integral of the continuous
$C^0(\GH)$-valued integrand $s\to \mathcal E u(s)\cdot E_s$. (Such standard
properties of Gelfand-Pettis vector-valued integrals are recalled in section
\ref{Gelfand-Pettis}.) Since $\overline{\mu }=\mu $, we have
$\mu (1)=\overline{\langle \mu ,1\rangle}$, and
$\mu  E_s=\overline{\mu  E_{1-s}}$ for $\Re(s)=\frac12$. The
two expressions for $u$ agree for all $u\in \Eis^{1+\ve}$, giving the
proposition. \qed

In fact, $\eta_a\in \Eis^{-\frac 12-\ve}$ for all $\ve>0$, because the $|\cdot|_{-\frac 12-\ve}$
norm is
\[
|\eta_a|_{-\frac12-\ve}^2 =\frac {1}{\langle1,1\rangle}+\frac 1{4\pi}\int_{-\infty}^\infty
\frac{|a^s+c_sa^{1-s}|^2}{(1+4t^2)^{\frac 12+\ve}}\,\d t
\ll_a \int_{-\infty}^\infty \frac{\d t}{(1+4t^2)^{\frac 12+\ve}} <\infty.
\]
Thus, the $\Eis^{-\frac 12-\ve}$ spectral expansion convergent in
$\Eis^{-\frac 12-\ve}$ is
\begin{equation}\label{spec-eta-a}
\eta_a=\frac {1}{\langle 1,1\rangle}+\frac{1}{4\pi i}
\int_{(\frac12)} (a^{1-s}+c_{1-s}a^s)\pt E_s\,\d s
\end{equation}

\subsection{Solving  $(-\Delta-\lambda_w)u=\eta_a$}
For $\Re(w)>\frac 12$ the equation  $(-\Delta-\lambda_w)u=\eta_a$
has an unique solution $v_{w,a}\in \Eis^{\frac 32-\ve}$ for all $\ve>0$,
with spectral expansion
obtained directly from that of $\eta_a$ via \eqref{S^ell}
\begin{equation}\label{a-spectral}
v_{w,a}=\frac{1}{\langle 1,1\rangle}
+\frac{1}{4\pi i}\int_{(\frac12)} 
(a^{1-s}+c^{}_{1-s}a^s)\cdot E_s\,\frac{\d s}{\lambda_s-\lambda_w}.
\end{equation}

\section{Pseudo-Laplacians on non-cuspidal automorphic\\ spectrum}

\subsection{Necessary condition for discrete spectrum}\label{necessary-condition}
For any $\theta\in \Eis^{-1}$ with $\overline{\theta}=\theta$ and
$\theta\not\in \Eis^0 $, let
$S_\theta$ be $-\Delta$ restricted to the domain $\Eis^\infty_c\cap \ker
\theta$. Symmetry of $S_\theta$ is inherited from $-\Delta$. The
{\itshape pseudo-Laplacian} $\widetilde S_\theta$ is the Friedrichs
extension of $S_\theta$.

As above, for $\theta\in \Eis^{-1}$, the distributional equation $(-\Delta
-\lambda)u=\theta$ has a unique solution $u\in \Eis^1$ for all $\lambda$
not in $\{0\}\cup [\frac 14,+\infty)$, via spectral expansions and
\eqref{S^ell}.

\begin{thm}\label{solvability-necessary-condition} 
For real $\lambda_w>\frac 14$, if the equation
$(-\Delta -\lambda_w)u=\theta$ has a solution in $\Eis^1$ then $\mathcal E
\theta(w)=0$. More precisely, existence of a solution implies that,
for all $\ve>0$,   
\[
\int_{\Im (w)-\ve}^{\Im(w)+\ve} \Big|\mathcal E \theta\Big(\frac 12+it\Big)\Big|\;\d t
\ll_{w,\ve} \ve^{\frac 32}
\]
\end{thm}

\proof The solution $u$ has spectral expansion
\[
u = \frac {\langle u,1\rangle\cdot 1}{\langle 1,1\rangle}
+ \frac{1}{4\pi i}\int_{(\frac 12)} \mathcal Eu(s)\cdot E_s \;\d s
\quad (\text{in $\Eis^1$})
\]
and
\[
(-\Delta-\lambda_w)u = (\lambda_1-\lambda_w)\frac{\langle u,1\rangle\cdot 1}{\langle 1,1\rangle}
+ \frac{1}{4\pi i}\int_{(\frac 12)} \mathcal Eu(s)\cdot (\lambda_s-\lambda_w)\cdot E_s \;\d s
\quad (\text{in $\Eis^{-1}$})
\]
Since $\theta$ itself has a spectral expansion in $\Eis^{-1}$, by
\eqref{S^ell} necessarily
\[
\mathcal E \theta(s) = \mathcal E u(s)\cdot (\lambda_s-\lambda_w)
\quad (\text{at least as locally-$L^2$ functions on $\textstyle\frac 12+i\mathbb R$})
\]
Further, from the Cauchy-Schwarz-Bunyakowsky inequality,
\begin{align}
&\int_{v_0-\ve}^{v_0+\ve} \Big|\mathcal E \theta\Big(\frac 12+it\Big)\Big|^2\;\d t
\le
\int_{v_0-\ve}^{v_0+\ve} \Big|\mathcal E u\Big(\frac 12+it\Big)\Big|\cdot |t^2-v_0^2|\;\d t
\notag \\
&\le
\left(\int_{v_0-\ve}^{v_0+\ve} \Big|\mathcal E u\Big(\frac 12+it\Big)\Big|\;\d t\right)^{\frac 12}
\cdot
\left(\int_{v_0-\ve}^{v_0+\ve} \left|t^2-v_0^2\right|^2\;\d t\right)^{\frac 12}
\ll_{w,\ve}
\|\mathcal E u\|\cdot \ve^{\frac 32}
\notag
\end{align}

\noindent as asserted.
\qed

As a corollary, we have a {\it necessary}, but in general not {\it
sufficient}, constraint on possible discrete spectrum of $\widetilde
S_\theta$: 
\begin{cor}\label{discrete-spectrum-if-any}
The discrete spectrum $\lambda_w>\frac 14$ of $\widetilde S_\theta$,
{\bf if any}, is of the form 
$w(1-w)$ for $w\in \frac 12+i\mathbb R$ such that $\mathcal E
\theta(w)=0$, in the sense of the previous theorem.
\end{cor}

\proof From Theorem \ref{domain-S-Theta-tilde}, any solution $u\in
\widetilde D_\theta$ to $(\widetilde S_\theta -\lambda_w)u=0$ is a
solution to a distributional equation $(-\Delta -\lambda_w)u=c\cdot
\theta$ for some constant $c$. For $u$ not identically $0$ and
$\lambda_w\not=0$, the constant $c$ cannot be $0$, since $(-\Delta
-\lambda_w)u=0$ has no non-zero solution in $\Eis^1$ for
$\lambda_w\not=0$. Thus, without loss of generality, take $c=1$, and
apply the theorem.
\qed

\begin{rem}\label{disclaimer-sufficient-condition}
For $\theta\in \Eis^{-1+\ve}$ for some $\ve>0$, theorem
\ref{mero-contn-in-V} will show that on $\Re(w)=\frac 12$ with
$w\not=\frac 12$, $\mathcal E \theta(w)=0$ is also {\itshape
sufficient} for existence of a solution $u\in \Eis^1$ to the
distributional equation
$(-\Delta-\lambda_w)u=\theta$. However, such $u$ is not in the domain
of the self-adjoint operator $\widetilde S_\theta$ unless also $\theta u=0$, which does not
follow from $\mathcal E\theta(w)=0$ in general. 
\end{rem}

\begin{rem}\label{disclaimer-about-if-any}
The corollary gives a definite relation between
the spectrum of a natural self-adjoint operator and the zeros of
$\mathcal E \theta(s)$, which, as in section
\ref{Heegner-distributions}, in many interesting cases 
an $L$-function or a finite linear combination of such.
However, there appears to be no general
assurance of {\itshape existence} of any discrete spectrum 
whatsoever.
\end{rem}


\begin{rem}\label{trivial-non-vanishing}
The larger point of our discussion of self-adjoint
operators is to prove that various quantities do not
vanish in $\Re(w)>\frac 12$ (off the real
line). Unsurprisingly, some non-vanishings are more trivial than might
be anticipated. For example, for any $\theta\in \Eis^{-1}$, with
$u^{}_{\theta,w}\in \Eis^1$ a solution of $(-\Delta-\lambda_w)u=\theta$, by
spectral theory
\[
\theta u^{}_{\theta, w} = \frac{\langle u,1\rangle\cdot \overline{\langle
\theta,1\rangle}}{\lambda_1-\lambda_w}
+\frac{1}{4\pi i}\int_{(\frac 12)} \frac{\mathcal E \theta(s)\cdot
E_s}{\lambda_s-\lambda_w}\;\d s
\]
For $\Re(w)>\frac 12$ and $w$ off $(\frac 12,1]$, it is elementary that
the imaginary part of that expression is non-zero. That is, $\theta
u^{}_{\theta,w}\not=0$ off the critical line and the real line. That is,
although facts about self-adjoint operators do yield these
particular conclusions, some of these conclusions are elementary.
\end{rem}

\subsection{Extensions of restrictions to non-cuspidal
pseudo-cuspforms}\label{pseudo-cuspforms}

Here we consider families of restrictions of $-\Delta$ similar to
\cite{Lax-Phillips}, pages 204--206, with attention to details. 
For fixed $a>1$, let $\Theta\subset L^2(\GH)$ be the space of
pseudo-Eisenstein series $\Psi_\varphi$ with test function $\varphi$
supported on $[a,\infty)$.  
Since $\Delta \Psi_\varphi = \Psi_{\Delta \varphi}$ the space $\Theta$
is stable under $\Delta$. Let 
$\Eis^0_\Theta$ be the orthogonal complement to $\Theta$ in $\Eis^0$.
Let $S_\Theta$ be the restriction of $-\Delta$ to $\Eis^\infty_c \cap
\Eis^0_\Theta$, and $\widetilde S_\Theta$ its Friedrichs
extension. To avoid potential ambiguities, we should be sure that
$S_\Theta$ is densely defined on $\Eis^0_\Theta$:

\begin{lem}\label{density-by-cutoffs} $\Eis^\infty_c \cap \Eis^0_\Theta$ is
dense in $\Eis^0_\Theta$, and $\Eis^\infty_c\cap \Eis^0_\Theta$ is
dense in $\Eis^1\cap \Eis^0_\Theta$ with the $\Eis^1$ topology.
\end{lem}

\proof To show that $\Eis^\infty_c\cap \Eis^0_\Theta$ is dense in
$\Eis^0_\Theta$, given a
sequence of pseudo-Eisenstein 
series $\Psi_{\varphi_i}\in \Eis^\infty_c$ converging to $f\in \Eis^0_\Theta$, we
produce a sequence of pseudo-Eisenstein series in $\Eis^\infty_c\cap \Eis^0_\Theta$
converging to $f$.  We will do so by cutting off the constant
terms of the $\Psi_{\varphi_i}$ at height $a$. Since the limit $f$ of the
$\Psi_{\varphi_i}$ has constant term vanishing above height $y=a$,
that part of the constant terms of the $\Psi_{\varphi_i}$
must become small. The explicit details are routine. \qed


Essentially as in \cite{Lax-Phillips} but restricting to the orthogonal
complement $\Eis^0$ to cuspforms, we have

\begin{thm}\label{discretization} The resolvent $(\widetilde S_\Theta
-\lambda_w)^{-1}$ is compact for $\lambda_w\not\in\mathbb R$, and
has a meromorphic continuation to $w\in\mathbb C$, giving a compact
operator for $w$ off a discrete subset of $(\frac 12+i\mathbb R)\cup
[0,1]$. In particular, $\widetilde S_\Theta$ has purely discrete spectrum.  
\end{thm}

\proof Since $(1+\widetilde S_\Theta)^{-1}:\Eis^0_\Theta \to
\Eis^1_\Theta$ is continuous with the (finer) $\Eis^1$-topology on
$\Eis^1$, we see below that it suffices to demonstrate a
Rellich-lemma-type compactness, namely, that the inclusion
$\Eis^1_\Theta\to \Eis^0_\Theta$ is a compact linear map. The
corresponding compactness for compact Riemannian manifolds, possibly
with boundary, is standard.

The total boundedness criterion for relative compactness
requires that, given $\ve>0$, the image of the unit ball $B\subset
\Eis^1_\Theta$ by the inclusion into $\Eis^0_\Theta$ can be covered by
finitely-many balls of radius $\ve$. The 
Rellich lemma on compact
Riemannian manifolds reduces the issue to an estimate on the $a$-tail
of the quotient $\GH$, that is, the image in $\GH$ of $\{z\in
\mathfrak H:\Im(z)\ge a\}$. Then the necessary estimate on the
$a$-tail will follows from the $\Eis^1_\Theta$- boundedness.  

We prove that, given $\ve>0$, there
is $c$ sufficiently large so that $\varphi_\infty\cdot B$ lies in a single
ball of radius $\ve$ inside $L^2(\GH)$, that is,
\[
\lim_{c\to\infty} \int_{y>c} |f(z)|^2\, \frac{\d x\;\d y}{y^2} 
\;\longrightarrow\; 0 
\]
uniformly for $|f|_{\Eis^1}\le 1$.

A precise choice of smooth truncations and understanding of their $\Eis^1$ norms is
needed. Fix a smooth real-valued function $\psi$ on $\mathbb R$ with
$\psi(y)=0$ for $y\le 0$, $0\le\psi(y)\le 1$ for $0<y<1$, and
$\psi(y)=1$ for $y\ge 1$. For $t>1$, let $\psi_t(y)=\psi(\frac yt -1)$, and
form a pseudo-Eisenstein series $\Psi_{\psi_t}$, a locally
finite sum. Then $\Psi_{\psi_t}\cdot f(x+iy)$ is a smoothly cut-off
tail of $f$ starting gradually at height $t$: $(1-\Psi_{\psi_t})\cdot f$ is
identically $0$ in the region where $y\ge 2t$, and in all images of
this region under $SL_2(\mathbb Z)$.

\begin{lem} The smooth truncation $\Psi_{\psi_t}\cdot f$ has
$\Eis^1$-norm dominated by that of $f$ itself, with implied constant
uniform in $t\ge 2$ and in $f\in\Eis^1$.  \end{lem}

\proof Routine computation. \qed

Returning to the proof of Theorem \ref{discretization}, given $c>1$,
we can assume that $f$ has been smoothly truncated so that in the
fundamental domain it is supported inside the region where $y\ge c$, and increasing
its $\Eis^1$ norm by at most a uniform factor. Let the
Fourier coefficients of $f(x+iy)$ in $x$ be $\widehat{f}(n)$, functions of
$y$. Take $y\ge c>a$, so the 
$\widehat{f}(0)$ vanishes. Using Plancherel for the Fourier expansion in $x$,
integrating over the part of $Y_\infty$ above $y=c$, letting $\mathcal
F$ denote Fourier transform in $x$, there is a direct computation
\begin{align*}
&\int\int_{y>c} |f|^2\;\frac{\d x\;\d y}{y^2}
\;\le\; \frac{1}{c^2} \int\int_{y>c} |f|^2\;\d x\,\d y
= \frac{1}{c^2} \sum_{n\not=0} \int_{y>c} |\widehat{f}(n)|^2\,\d y
\\
&
\;\le\; \frac{1}{c^2} \sum_{n\not=0} (2\pi n)^2\int_{y>c} |\widehat{f}(n)|^2\,\d y
=
 \frac{1}{c^2} \sum_{n\not=0} 
\int_{y>c} \Big|\,\mathcal F\,\frac{\partial f}{\partial x}(n)\Big|^2\,\d y
\\
&=
 \frac{1}{c^2} \int\int_{y>c} \Big|\frac{\partial f}{\partial
 x}\Big|^2\,\d x\,\d y
=
\frac{1}{c^2} \int\int_{y>c} -\frac{\partial^2 f}{\partial x^2} \cdot
\fbar(x)\; \d x\,\d y
\\
&\le
 \frac{1}{c^2} \int\int_{y>c}
\left(-\frac{\partial^2 f}{\partial x^2} \cdot \fbar(x)
 -
\frac{\partial^2 f}{\partial y^2} \cdot \fbar(x)\right)
\; \d x\,\d y
=
\frac{1}{c^2} \int\int_{y>c} -\Delta f \cdot \fbar \;\frac{\d x\,\d
 y}{y^2}
 \\
&\le
 \frac{1}{c^2} \int\int_{\GH} -\Delta f \cdot \fbar \;\frac{\d x\,\d
 y}{y^2}
= 
\frac{1}{c^2} \, |f|^2_{\Eis^1} \;\le\; \frac{1}{c^2}
\end{align*}
This uniform bound completes the proof that the image of the unit ball
of $\Eis^1_\Theta$ in $\Eis^0_\Theta $ is totally bounded. Thus, the
inclusion is a compact map.

As earlier, Friedrichs' construction shows that
$(\widetilde S_\Theta-\lambda)^{-1}:\Eis^0_\Theta \to \Eis^1_\Theta$
is continuous even 
with the stronger topology of $\Eis^1_\Theta$. Thus, the composition
\hbox{$\Eis^0_\Theta \to \Eis^1_\Theta \subset \Eis^0_\Theta$}\hfill by
\[
f\longrightarrow (\widetilde S_\Theta-\lambda)^{-1}f \to (\widetilde S_\Theta-\lambda)^{-1}f
\]
is the composition of a continuous operator with a compact operator,
so is compact. Thus, $(\widetilde S_\Theta-\lambda)^{-1}:\Eis^0_\Theta
\to \Eis^0_\Theta$ is a compact operator.
Thus, for $\lambda$ off a discrete set of points in $\mathbb C$,
$\widetilde S_\Theta$ has compact resolvent $(\widetilde
S_\Theta-\lambda)^{-1}$, and the parametrized family of compact
operators  
$(\widetilde S_\Theta-\lambda)^{-1} : \Eis^0_\Theta \longrightarrow \Eis^0_\Theta $
is meromorphic in $\lambda\in\mathbb C$.

We recall the standard argument (see \cite{Kato}, page 187 and
preceding, for example) for the fact that that, for a (not necessarily
bounded) normal operator $T$, if $T^{-1}$ exists and is compact, then
$(T-\lambda)^{-1}$ exists and is a compact operator for $\lambda$ off
a discrete set in $\mathbb C$, and is {\it meromorphic} in $\lambda$.
First, from the spectral theory of normal compact
operators, the non-zero spectrum of compact $T^{-1}$ is all point
spectrum. We claim that the spectrum
of $T$ and non-zero spectrum of $T^{-1}$ are in the obvious bijection
$\lambda\leftrightarrow \lambda^{-1}$. From the algebraic identities 
$T^{-1}-\lambda^{-1}=T^{-1}(\lambda-T)\lambda^{-1}$
and $T-\lambda = T (\lambda^{-1}-T^{-1}) \lambda$, 
failure of either $T-\lambda$ or $T^{-1}-\lambda^{-1}$ to be {\it injective}
forces the failure of the other, so the point spectra are
identical. For (non-zero) $\lambda^{-1}$ not an eigenvalue of {\it
compact} $T^{-1}$, $T^{-1}-\lambda^{-1}$ is injective {\it and} has a
continuous, everywhere-defined inverse. 
That $S-\lambda$ is {\it surjective} for compact normal $S$ and $\lambda\not=0$
not an eigenvalue is an easy part of Fredholm theory.
For such $\lambda$, inverting the
relation $T-\lambda = T (\lambda^{-1}-T^{-1}) \lambda$ gives
\[
(T-\lambda)^{-1} = \lambda^{-1} (\lambda^{-1}-T^{-1})^{-1} T^{-1}
\]
from which $(T-\lambda)^{-1}$ is continuous and everywhere-defined. That
is, $\lambda$ is {\it not} in the spectrum of $T$. Finally, $\lambda=0$ is
not in the spectrum of $T$, because $T^{-1}$ exists and is
continuous. This establishes the bijection.  
\qed


Essentially as in \cite{Lax-Phillips}, identification of the
eigenfunctions for $\widetilde S_\Theta$ depends on the cut-off
height $a>1$ and reduction theory.
The {\itshape truncation} $\wedge^a E_s$ of an Eisenstein series at
height $\Im(z)=a$ is as follows. With $y=\Im(z)$, let
$\tau_s(z)$ be $y^s+c_sy^{1-s}$ for $y\ge a$ and $0$ for
$0<y<a$, and form a pseudo-Eisenstein series  
\[
\Psi_{s,a}(z) = \sum_{\gamma\in\Gamma_\infty\backslash\Gamma}
\tau_s(\gamma z)
\]
Even though $\tau_s$ is {\itshape not} a test function, by reduction
theory this is a locally finite sum, so converges uniformly absolutely 
pointwise. For $a>1$, inside the standard fundamental domain for
$SL_2(\mathbb Z)$, by reduction theory, $\Psi_{s,a}(z)$ is $0$ unless
$y\ge a$, and is $y^s+c_sy^{1-s}$ for $y\ge a$. The truncated
Eisenstein series is  
\[
\wedge^a E_s = E_s - \Psi_{s,a}
\]
By design, for $a>1$, inside the standard fundamental domain for
$SL_2(\mathbb Z)$, this truncation makes the constant term vanish
above height $y=a$.

\begin{thm}\label{exotic-eigenfunctions} For $a>1$, spectral
parameters $w$ for eigenvalues $\lambda_w>\frac 14$ of $\widetilde
S_\Theta$ are exactly the zeros of the constant term $a^w+c_wa^{1-w}$
of the Eisenstein series $E_w$ with $\Re(w)=\frac 12$. The
corresponding eigenvalues $\lambda_w$ are {\itshape simple}, the
corresponding eigenfunctions are solutions $u\in \Eis^1_\Theta$ of the
equation $(-\Delta-\lambda)u=\eta_a$, and up to constants these are
truncated Eisenstein series $\wedge^a E_w$. Specifically,
\[
(-\Delta-\lambda_w) \wedge^a E_w = 2(1-2w)a^{w+1}\cdot \eta_a
\]
\end{thm}

\begin{rem} In particular, all eigenfunctions with eigenvalues
$\lambda>\frac 14$ fail to be smooth, since their constant terms will
be continuous, but have discontinuous first derivative (in $y=\Im(z)$)
at height $a$. In this regard, such eigenfunctions are {\itshape exotic}.
\end{rem}

\proof Because the homogeneous equation $(-\Delta-\lambda)u=0$ has no
non-zero solution in $\Eis^1_\Theta$, it suffices to identify the
possible $\theta$ in the $\Eis^{-1}$ closure $\Theta_{-1}$ of $\Theta$
that could fit into an equation
$(-\Delta-\lambda)u=\theta$ with $u$ in the $\Eis^1$ closure of $\Eis^\infty_c\cap
\Eis^1_\Theta$. On one hand, because $a>1$, 
$\Theta_{-1}$ consists of distributions which, on the
standard fundamental domain, have support inside the Siegel set $\mathfrak
S_a=\{x+iy\in\mathfrak H:y\ge a\}$. Further, on
$C_a=\Gamma_\infty\backslash \mathfrak S_a$, the circle $S^1=\mathbb
Z\backslash \mathbb R$ acts by translations, descending to the
quotient from $\mathfrak H$.  By reduction theory, the restrictions to
$C_a$ of every pseudo-Eisenstein series $\Psi_\varphi$ with
$\varphi\in C^\infty_c[a,\infty)$ are invariant under $S^1$, so
anything in the $\Eis^{-1}$ closure is likewise translation-invariant.

On the other hand, $\Eis^\infty_c \cap \Eis^1_\Theta$ consists of functions with
constant term vanishing in $y\ge a$, and taking $\Eis^1$ completion
preserves this property. Since $\Theta_{-1}$ is $S^1$-invariant and
the Laplacian commutes with $S^1$, it suffices to look at
$S^1$-integral averages of $v$ restricted to some cylinder $C_b$ with
$a>b>1$. Such an integral is a restriction of the constant term $c_Pv$
to $C_b$, and vanishes in $y>a$.

Thus, in the standard fundamental domain, the support of 
distributions $\theta\in \Theta_{-1}$ fitting into an equation
$(-\Delta-\lambda)u=\theta$, with $u$ in the $\Eis^1$ closure of
$\Eis^\infty_c\cap 
\Eis^1_\Theta$, is inside the image of the set where $y=a$, and $\theta$ is
$S^1$-invariant. By the classification of distributions
supported on submanifolds, such $\theta$ is a derivative
normal to the circle $\Gamma_\infty\backslash \{x+iy:y=a\}$ followed
by application of an $S^1$-invariant distribution on the circle. The latter
must be integration over the circle, by classification of invariant
distributions. By standard Sobolev theory, 
there can be no actual derivatives, for the resulting distribution 
to lie in the $-1$ Sobolev space. Thus, up to a constant multiple, $\theta$
is the evaluation of constant term at height $a$ functional $\eta_a$,
namely, $\eta_a f=c_Pf(ia)$.

Since $|a^s+c_sa^{1-s}|\leq 2\sqrt{a}$ for $\Re(s)=\frac 12$, in
fact $\eta_a\in \Eis^{-\frac 12-\ve}$ for every $\ve>0$.

As in Theorem \ref{solvability-necessary-condition} and its proof, from the spectral
relation \eqref{S^ell}, for $v\in \Eis^1$ with
$(-\Delta-\lambda_w)u=\eta_a$, the spectral expansions  
\begin{align}
&u = \frac{\langle u,1\rangle\cdot 1}{\langle 1,1\rangle}
+ \frac{1}{4\pi i}\int_{(\frac 12)} \mathcal E u(s)\cdot E_s\;\d s
\\
&\eta_a = \frac{\langle \eta_a,1\rangle\cdot 1}{\langle 1,1\rangle}
+ \frac{1}{4\pi i}\int_{(\frac 12)} \mathcal E \eta_a(s)\cdot E_s\;\d s
\end{align}
\noindent converging in $\Eis^1$ and $\Eis^{-1}$, respectively, have $\mathcal E
\eta_a(s)=\eta_a(E_{1-s})=a^{1-s}+c_{1-s}a^s$, and the expansions are
related by 
$(\lambda_s-\lambda_w)\cdot \mathcal E u(s) =\mathcal E \eta_a(s)$,  
so $\mathcal E\eta_a(w)=0$ in a strong sense. Since $\Re(w)=\frac12$,
by $c_wc_{1-w}=1$, we have $a^w+c_wa^{1-w}=0$.

To finish the proof, it suffices to show that
$(-\Delta-\lambda_w)\wedge^a E_w$ is a scalar multiple of $\eta_a$
when $a^w+c_wa^{1-w}=0$. Indeed, with $a>1$, in a fundamental domain,
away from $y=a$ we have $(-\Delta-\lambda_w)\wedge^a E_w=0$
locally. Further, in $y>1$, the differential operator annihilates all
Fourier components of $E_w$ but the constant term, and in both $1<y<a$
and $y>a$ does also annihilate the constant term. To compute near
$y=a$, let $H$ be the Heaviside function $H(y)=0$ for $y<0$ and
$H(y)=1$ for $y>0$. Thus, near $y=a$, as functions of $y$ independent
of $x$,
\[
\begin{split}
(-\Delta-\lambda_w)\wedge^a E_w
=
(-\Delta-\lambda_w)
\Big(
H(a-y)\cdot (y^w+c_wy^{1-w})
\Big)
\hskip80pt
\\
=
(-y^2\frac{\partial^2}{\partial y^2}-w(1-w))
\Big(
H(a-y)\cdot (y^w+c_wy^{1-w})
\Big)
\hskip80pt
\\
=-y^2\Big(
H''(a-y)(y^w+c_wy^{1-w})
+
2H'(a-y)(y^w+c_wy^{1-w})'
\hskip80pt
\\
+
H(a-y)(y^w+c_wy^{1-w})''
\Big)
- w(1-w)H(a-y)(y^w+c_wy^{1-w})
\\
=
-y^2\Big(
\delta'_a\cdot (y^w+c_wy^{1-w})
- 2\delta_a\cdot (wy^{w-1}+(1-w)c_wy^{-w})
\Big)
\hskip80pt
\end{split}
\]
Since $a^w+c_wa^{1-w}=0$, the term with $\delta'_a$ vanishes, and the
rest simplifies to
\[
\begin{split}
(-\Delta-\lambda_w)\wedge^a E_w
=
-2a\delta_a \cdot (wa^w+(1-w)c_wa^{1-w})
\\
=
-2\delta_a \cdot (2w-1)a^{w+1}
\hskip100pt
\end{split}
\]
on functions of $y$ independent of $x$. Thus, this is
$2(1-2w)a^{w+1}\cdot \eta_a$. \qed

\subsection{Exotic eigenfunction expansions}

Keep $a>1$ fixed, and, as above, $\Theta$ the collection of
pseudo-Eisenstein series formed from test function data $\varphi$
supported in $[a,+\infty)$.
Any $u\in \Eis^1$ with $\eta_b u=0$ for all $b\ge a>1$ lies in
$\Eis^0_\Theta $, so admits a spectral expansion in terms of the
eigenfunctions for $\widetilde S_\Theta$, converging in the topology
of $\Eis^0$. However, we want to apply
functionals in $\Eis^{-1}$ (termwise) to such an spectral expansion, which requires that
the expansion converge in $\Eis^1$. Thus, we need the following stronger
analogue of Lemma \ref{density-by-cutoffs}:

\begin{thm}\label{density-by-smooth-cutoffs} 
With $a>1$, $\Eis^\infty_c\cap \Eis^0_\Theta $ is dense in $\Eis^1\cap
\Eis^0_\Theta $ with the $\Eis^1$ topology.
\end{thm}

\proof Given a sequence of pseudo-Eisenstein series
$\Psi_{\varphi_i}\in \Eis^\infty_c$ converging to $f\in \Eis^1_\Theta$ in the
topology of $\Eis^1$, we produce a sequence of pseudo-Eisenstein series
in $\Eis^\infty_c\cap \Eis^0_\Theta $ converging to $f$ in the
topology of $\Eis^1$.  We will do so by smooth cut-offs of the
constant terms of the $\Psi_{\varphi_i}$. Since the limit $f$ of the
$\Psi_{\varphi_i}$ has constant term vanishing above height $y=a$ and
is in $\Eis^0_\Theta $, that part of the constant terms of the
$\Psi_{\varphi_i}$ becomes small. More precisely, we proceed as follows.

Let $g$ be a smooth real-valued function on $\mathbb R$ with $g(y)=0$
for $y<-1$, $0\le g(y)\le 1$ for $-1\le y\le 0$, and $g(y)=1$ for
$y\ge 0$. For $\ve>0$, let $g_\ve(y)=g((y-a)/\ve)$.
Fix real $b$ with $a>b>1$. Given $\Psi_{\varphi_i}\to f\in \Eis^0_\Theta $,
the $b$-{\itshape tail} of the constant term of $\Psi_{\varphi_i}$ is
$\tau_i(y)=c_p\Psi_{\varphi_i}(y)$ for 
$y\ge b$, and $\tau_i(y)=0$ for $0<y\le b$. By design,
$\Psi_{\varphi_i}-\Psi_{g_\ve\cdot \tau_i}\in \Eis^\infty_c\cap
\Eis^0_\Theta $ for small 
$\ve$. We will show that, as $i\to+\infty$, for $\ve_i$ sufficiently small
depending on $i$, the $\Eis^1$-norms of $\Psi_{g_{\ve_i}\cdot \tau_i}$ go
to $0$, so $\Psi_{\varphi_i}-\Psi_{g_{\ve_i}\cdot \tau_i}\to f$ in the
$\Eis^1$ topology.

For $b>1$, let $\mathfrak S_b=\{x+iy\in \mathfrak H:y\ge b,\;|x|\le \frac
12\}$. By reduction theory, the cylinder $C_b=\Gamma_\infty\backslash
(\Gamma_\infty\cdot \mathfrak S_b)$ maps homeomorphically to its image
in $\GH$. For $f\in C^\infty_c(\GH)$, let
\[
|f|^2_{H^1(C_b)}
=
\int_{C_b} |f(z)|^2- \Delta f\cdot \fbar\;\frac{\d x\,\d y}{y^2}
\le
\int_{\GH} |f(z)|^2- \Delta f\cdot \fbar\;\frac{\d x\,\d y}{y^2}
\]
For each $b>1$, let $H^1(C_b)$ be the completion of
$C^\infty_c(\GH)$ with respect to the semi-norm $|\cdot|_{H^1(C_b)}$,
allowing for collapsing. The cylinders $C_b$ admit natural actions of
the circle group $S^1=\mathbb Z\backslash \mathbb R$, by translation,
inherited from the translation of the real part of $x+iy\in\mathfrak
H$. As usual, this induces a continuous action of $S^1$ on
$H^1(C_b)$. Thus, the map $F \to c_PF$ gives continuous maps of the 
spaces $H^1(C_b)$ to themselves. Thus, $c_P\Psi_{\varphi_i}$ goes to
$c_P f$ in $H^1(C_b)$, and $c_P\Psi_{\varphi_i}\to c_P f=0$ in
$H^1(C_a)$. 


To have a useful Leibniz rule for differentiation, it is convenient to
rewrite the norms: for $f\in
C^\infty_c(\GH)$, put  
\[
|f|^2_{H_1} = |f|^2_{L^2(\GH)} + |(|\nabla f|_{\mathfrak s})|^2_{L^2(\GH)}
\]
where $\nabla$ is the left $SL_2(\mathbb R)$-invariant, right
$SO_2(\mathbb R)$-equivariant tangent-space-valued gradient on
$SL_2(\mathbb R)$, which therefore descends to $\mathfrak H$ and to
$\GH$, and $|\cdot|_{\mathfrak s}$ is a natural $SO_2(\mathbb
R)$-invariant norm on the tangent space(s). 
More explicitly, let $\mathfrak s$ be the space of symmetric
$2$-by-$2$ matrices of trace $0$, identified with the tangent space at every
point of $\mathfrak H$ via left translation of the exponential map:
for $\beta\in \mathfrak s$, as usual the associated left
$SL_2(\mathbb R)$-invariant differential operator $X_\beta$ is 
\[
(X_\beta f)(g) = \frac{\partial}{\partial t}\Big|_{t=0} f(ge^{t\cdot \beta})
\]
It is easy to describe $\nabla$ in coordinates, even though it is
provably independent of coordinates: let $h=\begin{pmatrix}1
& \hfill 0 \cr 0 & -1\end{pmatrix}$ and $\sigma=\begin{pmatrix}0 & 1
\cr 1 & 0\end{pmatrix}$, and put 
\[
\nabla F(g) = X_h F(g)\cdot h + X_\sigma F(g)\cdot
\sigma
\;\in\; \mathfrak s \otimes_{\mathbb R} \mathbb C
\]
Up to a scalar, the $SO_2(\mathbb R)$-invariant hermitian inner
product $\langle,\rangle_{\mathfrak s}$ on the complexified
$\mathfrak s$, makes $h,\sigma$ an orthonormal basis for $\mathfrak
s$. Let $|\cdot|_{\mathfrak s}$ be the associated norm. The essential
property is the integration by parts identity 
\[
\int_{\GH} \langle \nabla F_1, \nabla F_2\rangle_{\mathfrak s}
=
\int_{\GH} -\Delta F_1\cdot \overline F_2
\]
for $F_1,F_2\in C^\infty_c(\GH)$. The advantage of this formulation is
that, extending $\nabla$ by continuity in the $H^1$ topology, $\nabla
F$ exists (in an $L^2$ sense) for $F\in H^1(C_b)$. Thus, we can say
that
\[
|F|^2_{H^1(C_b)}
=
|F|^2_{L^2(\GH)} + |\;|\nabla F|_{\mathfrak s}|^2_{L^2(\GH)}
\]
Then
\[
\begin{split}
|\Psi_{g_\ve\cdot \tau_i}|_{\Eis^1}
=
|\Psi_{g_\ve\cdot \tau_i}|_{H^1(C_{a-\ve})}
=
|g_\ve\cdot \tau_i|_{H^1(C_{a-\ve})}
\hskip80pt
\\
\le
|(g_\ve-1)\cdot \tau_i |_{H^1(C_{a-\ve})}
+ |\tau_i - c_Pf|_{H^1(C_{a-\ve})}
+ |c_Pf|_{H^1(C_{a-\ve})}
\end{split}
\]
The middle summand goes to $0$:
\[
|\tau_i - c_Pf|_{H^1(C_{a-\ve})}
\le
|c_P\Psi_{\varphi_i}-c_Pf|_{\Eis^1}
\le
|\Psi_{\varphi_i}-f|_{\Eis^1}
\longrightarrow 0
\]
The first and third summands require somewhat more care. Estimate 
\[
\begin{split}
|(g_\ve-1)\cdot \tau_i|^2_{H^1(C_{a-\ve})}
=
\int_{C_{a-\ve}} |(g_\ve-1)\tau_i|^2
+ |\nabla(g_t-1)\tau_i|_{\mathfrak s}^2
\hskip80pt
\\
\le
\!\!\!\int\limits_{C_{a-\ve}} \!\! |g_\ve-1|^2\cdot (|\tau_i|^2
+|\nabla \tau_i|_{\mathfrak s}^2)
+
\!\!\!\!\int\limits_{C_{a-\ve}} \!\!|\nabla g_\ve|_{\mathfrak s}^2\cdot |\tau_i|^2
+
\!\!\!\!\int\limits_{C_{a-\ve}} \!\!
2|g_\ve|\cdot|\nabla g_\ve|_{\mathfrak s}\cdot |\tau_i|\cdot|\nabla \tau_i|_{\mathfrak s}
\\
\end{split}
\]
The first summand in the latter expression goes to $0$ as $\ve\to 0^+$
because $g_\ve-1=0$ when $y\ge a$, and $\tau_i$ and $|\nabla
\tau_i|_{\mathfrak s}$ are continuous.

In terms of the coordinates $z=x+iy$ on $\mathfrak H$, for a smooth
function $F$ a standard computation gives
\[
\nabla F
=
y \frac{\partial F}{\partial x}\cdot \sigma
+
y \frac{\partial F}{\partial y}\cdot h
\]
so
\[
|\nabla g_\ve(x+iy)|_{\mathfrak s}
=
|\frac{1}{\ve}\cdot y\,g'((y-a)/\ve)\cdot h|_{\mathfrak s}
=
\frac{1}{\ve}\cdot |y\,g'((y-a)/\ve)|
\ll_g
\frac{1}{\ve}
\]
Similarly, since $\tau_i$ is a function of $y$ independent of $z$,
$\nabla \tau_i=y\tau'_i(y)\cdot h$. By the fundamental theorem of
calculus and the Cauchy-Schwarz-Bunyakowsky inequality, we recover an easy instance 
of a Sobolev inequality:
\[
\begin{split}
|\tau_i(a-v)| = \Big|0 - \int_0^v \tau'_i(a-v)\;\d v\Big|
\hskip100pt
\\
\le
\Big(\int_0^v |\tau'_i(a-v)|^2\;\d v\Big)^{\frac 12}
\cdot
\Big(\int_0^v 1^2\;\d v\Big)^{\frac 12}
=
o(\sqrt{v})
\end{split}
\]
with Landau's little-$o$ notation, since $\tau'_i$ is locally $L^2$. Thus,
\[
\begin{split}
\int_{C_{a-\ve}} 
|g_\ve| \cdot |\nabla g_\ve|_{\mathfrak s}\cdot |\tau_i|\cdot |\nabla
\tau_i|_{\mathfrak s}
\le
\frac{1}{\ve}\cdot o(1)\cdot \sqrt{\ve}\cdot
\int_0^\ve
|\nabla \tau_i|_{\mathfrak s}
\hskip80pt
\\
\le
\frac{1}{\ve}\cdot o(1)\cdot \sqrt{\ve}\cdot
\Big(\int_0^\ve |\tau'_i|^2\Big)^{\frac 12}
\cdot
\Big(\int_0^\ve 1^2\Big)^{\frac 12}
\ll_{\tau_i}
\frac{1}{\ve}\cdot o(1)\cdot \sqrt{\ve}\cdot\sqrt{\ve}
= o(1)
\end{split}
\]
That is, the summand $\int_{C_{a-\ve}}|g_\ve|\cdot|\nabla g_\ve|_{\mathfrak s}\cdot
|\tau_i|\cdot|\nabla \tau_i|_{\mathfrak s}$ goes to $0$.
Using the same subordinate estimates, 
\[
\int_{C_{a-\ve}}|\nabla g_\ve|^2_{\mathfrak s}\cdot
|\tau_i|^2
\ll
\frac{1}{\ve^2} \int_0^\ve \big(o(1)\cdot\sqrt{v}\big)^2 \;\d v
=
\frac{1}{\ve^2} \cdot o(1)\cdot \frac{\ve^2}{2}
\longrightarrow 0 
\]
Thus, taking the $\ve_i$ sufficiently small, the smooth truncations
$\Psi_{\varphi_i}-\Psi_{g_{\varepsilon_i}\cdot \tau_i}$ of the
$\Psi_{\varphi_i}$ are in $\Eis^\infty_c\cap \Eis^0_\Theta $, and
still converge to $f$ in $\Eis^1$. 
\qed


\begin{cor}\label{compare-norms} An orthogonal basis for $\Eis^0_\Theta $
consisting of $\widetilde{S}_\Theta$-eigenfunctions is an orthogonal
basis for $\Eis^1\cap \Eis^0_\Theta $, as well. In particular, for
eigenfunction $f$ with eigenvalue $\lambda$, we have $\langle 
f,f\rangle_{\Eis^1}=\lambda\cdot \langle f,f\rangle$.
\end{cor}

\proof Since $\widetilde S_\Theta^{-1}(\Eis^0_\Theta )\supset
\Eis^\infty_c\cap \Eis^0_\Theta $,
by the theorem $\widetilde S_\Theta^{-1}(\Eis^0_\Theta )$ is dense in
$\Eis^1\cap \Eis^0_\Theta $. Since finite linear combinations of the
$\widetilde S_\Theta$-eigenfunctions are dense in $\Eis^0_\Theta $ and
$\widetilde S^{-1}_\Theta$ is continuous, their images in $\Eis^1\cap
\Eis^0_\Theta $ are dense there. 
From Proposition \ref{domains}, for two $\widetilde S_\Theta
$-eigenfunctions $u,v$ with (real) eigenvalues $\lambda,\mu$, 
\[
\langle u,v\rangle_{\Eis^1}
=
\langle S_\Theta^\#u,v\rangle_{\Eis^{-1}\times V}
=
\langle \widetilde S_\Theta u,v\rangle_{\Eis^{-1}\times V}
=
\langle \widetilde S_\Theta u,v\rangle
=
\lambda \langle u,v\rangle
\]
Symmetrically, $\langle u,v\rangle_{\Eis^1}=\mu\langle u,v\rangle$.
Thus, orthogonality of eigenfunctions in $\Eis^0_\Theta \subset V$ implies
orthogonality in $\Eis^1\cap \Eis^0_\Theta \subset \Eis^1$, and for $u=v$ we have
$\langle u,u\rangle_{\Eis^1}=\lambda\cdot \langle u,u\rangle$.
\qed

\begin{cor}\label{exotic-expansion} Let $\{u_k:k=1,2,\ldots\}$ be the
eigenfunctions for $\widetilde S_\Theta$, with eigenvalues
$\lambda_k$. For $f\in \Eis^1\cap \Eis^0_\Theta $, 
\[
f
=
\sum_{k\ge 1}
\langle f,u_k\rangle_V\cdot \frac{u_k}{\langle u_k,u_k\rangle_V}
=
\sum_{k\ge 1} \langle f,u_k\rangle_{\Eis^1}\cdot \frac{u_k}{\langle u_k,u_k\rangle_{\Eis^1}}
\]
and these expansions converge to $f$ not only in $\Eis^0_\Theta $, but also
in the finer topology of $\Eis^1\cap \Eis^0_\Theta $. 
\end{cor}

\proof Again, by the theorem, since $f\in \Eis^1\cap \Eis^0_\Theta $ is in the
closure of $\Eis^\infty_c\cap \Eis^0_\Theta $, and $\{u_k\}$ is an
orthogonal basis for 
$\Eis^1\cap \Eis^0_\Theta $, such $f$ has an expansion
\[
f = \sum_{k\ge 1} \langle f,u_k\rangle_{\Eis^1}\cdot \frac{u_k}{\langle u_k,u_k\rangle_{\Eis^1}}
\]
convergent in $\Eis^1\cap \Eis^0_\Theta $. As in the previous proof, from
Proposition \ref{domains}, 
\[
\langle f,u_k\rangle_{\Eis^1}
=
\langle f,S_\Theta^\#u_k\rangle_{V\times \Eis^{-1}}
=
\langle f,\widetilde S_\Theta u_k\rangle_{V\times \Eis^{-1}}
=
\langle f,\widetilde S_\Theta u_k\rangle
=
\lambda_k \langle f,u_k\rangle
\]
Thus, for every $u_k$, 
\[
\langle f,u_k\rangle_{\Eis^1}\cdot \frac{u_k}{\langle
u_k,u_k\rangle_{\Eis^1}}
=
\lambda_k\langle f,u_k\rangle\cdot \frac{u_k}{\lambda_k \langle
u_k,u_k\rangle}
=
\langle f,u_k\rangle\cdot \frac{u_k}{\langle u_k,u_k\rangle}
\]
giving the termwise equality of the expansions.
\qed


\subsection{Extensions of restrictions by Eisenstein-Heegner and constant-term
constraints}\label{two-condition}
Fix a finite real-linear combination $\theta$ of Eisenstein-Heegner
distributions $\theta_d\in \Eis^{-\frac12-\ve}$. Fix $a>1$.
Let $S_{\theta,a}$ be $-\Delta$ restricted to domain
\[
\Eis^\infty_c \cap \ker \theta \cap \ker\eta_a.
\]
Note that $S_{\theta,a} \geqslant \frac 14$ on the continuous spectrum,
and its domain excludes constants. Symmetry of $S_{\theta,a}$ is
inherited from $S$.  The {\itshape pseudo-Laplacian} $\widetilde
S_{\theta,a}$ is the Friedrichs extension of $S_{\theta,a}$ on 
$\Eis^1\cap \ker \theta \cap \ker\eta_a$, with
$\widetilde S_{\theta,a} \geqslant 0$.


\begin{thm}\label{equivalence-thm} Given
a finite real-linear combination $\theta$ of Eisenstein-Heegner distributions
$\theta_d$, for all 
$a$ with $\Im z\not=a$ for all the Heegner points $z$ involved, 
the Friedrichs extension
${\widetilde S}_{\theta,a}$ ignores $\theta$ and $\eta_a$,
in the sense that for $u$ in the domain of ${\widetilde S}_{\theta,a}$ the
eigenvector condition $({\widetilde S}_{\theta,a} -\lambda_w)\, u = 0$
is equivalent to the satisfaction of the equation
\[
(S_{\theta,a}^\# -\lambda_w)\, u 
= A\cdot \theta+B\cdot \eta_a \quad (\text{\rm for some }
A,B\in\mathbb C).
\]
\end{thm}
\proof The point is to show that $(\mathbb C\cdot \theta\cap \mathbb
\eta_a)\cap(j^*\cap \Lambda)\Eis^0 = 
\{0\}$. Then Theorem \ref{domain-S-Theta-tilde} applies.

There is an unique {\itshape highest} Heegner point $z_0$ 
appearing in $\theta$.  In fact, if $4|d$ the highest Heegner point in $\theta_d$ is
$i\sqrt{|d|}/2$, while otherwise it is $(1+i \sqrt{|d|})/2$.  Then it suffices to
take pseudo-Eisenstein series $f_n = \Psi_{\varphi_n}$ with $\varphi_n(\Im z_0)=1$
and $\varphi_n(\Im z')=1$ for the other Heegner points $z'$ appearing in
%
%
$\theta$. Thus, $A\cdot\theta+B\cdot \eta_a$ is not in $\Eis^0$
for $A\not=0$. Thus, it suffices to show that $\eta_a$ is not in
$\Eis^0$, which is even simpler.~\qed

\section{Eigenfunctions of pseudo-Laplacians}

\subsection{Determination of eigenfunctions}
We continue to keep fixed a finite real-linear combination $\theta$ of
Eisenstein-Heegner distributions $\theta_d$. 
From the preceding Theorem \ref{equivalence-thm},
a solution $u\in \Eis^1$ of an equation
\[
(S^\# - \lambda_w) u = A\cdot  \theta + B\cdot \eta_a
\]
is a $\widetilde S_{\theta,a}$-{\itshape eigenfunction}
precisely when $u\in \ker \theta\cap \ker \eta_a$.
As above, let $v_{w,a}$ be the unique solution in $\Eis^{\frac 32-\ve}$
to $(-\Delta-\lambda_w) = \eta_a$ and similarly let $u_{\theta,w}$
be the unique solution in $\Eis^{1+\ve}$ to 
$(-\Delta -\lambda_w) u =\theta$.  The condition for
existence of a non-zero solution $(A,B)$ to the homogeneous system
\begin{equation}\label{basic-system}
\begin{cases}
\ \theta(A u_{\theta,w}+B v_{a,w}) &= \, 0\\
\ \eta_a (A u_{\theta,w}+B v_{a,w}) & =\, 0
\end{cases}
\end{equation}
is the vanishing of the determinant:
\begin{equation}\label{resolvent}
\det
\begin{pmatrix}
 \theta(u^{}_{\theta,w}) & \theta (v_{w,a}) \\
{} & {}\\
 \eta_a(u^{}_{\theta,w}) &  \eta_a (v_{w,a})
\end{pmatrix}
= \, 0
\end{equation}
We compute the components.


\subsection{Computing $\eta_a(v_{w,a})$ for $a > 1$ and $\Re(w)>\frac 12$}
By the spectral expansion \eqref{spec-eta-a}, this is
\begin{align}
&\eta_a(v_{w,a})=\frac {1}{\langle 1,1\rangle} + \frac {1}{4\pi i}\int_{(\frac 12)}
(a^{1-s}+c^{}_{1-s}a^s)(a^s+c^{}_s a^{1-s}) \frac{\d s}{\lambda_s-\lambda_w}\notag\\
&\ = \frac {1}{\langle 1,1\rangle} + \frac {1}{4\pi i} \int_{(\frac 12)}
(a+c^{}_{1-s}a^{2s}+c^{}_s a^{2-2s}+a) \frac{\d s}{\lambda_s-\lambda_w}
\quad \text{(use $c^{}_s c^{}_{1-s}=1$)}\notag\\
&\ = \frac {1}{\langle 1,1\rangle} +\frac {1}{2\pi i} \int_{(\frac 12)}
(a + c^{}_s a^{2-2s}) \frac{\d s}{\lambda_s-\lambda_w}
\quad\text{($s\to 1-s$ in one term)} \label{eta-D-v}
\end{align}
The behavior of $c_s$ as a function of $s$ in the half-plane $\sigma\geqslant \frac 12$
is easily determined from the second formula in \eqref{c_s}.  We find that
$c^{}_s$ has a simple pole at $s=1$ with residue 
$\displaystyle \frac 3\pi = \frac {1}{\langle 1,1\rangle}$ and is of  order
$\sqrt{\sigma} (|t|+1)^{\frac 12+\ve}$ for any fixed $\ve>0$,
uniformly in $\sigma$ and $t$.  For $a>1$, $|a^{2-2s}|=a^{2-2\sigma}$ goes
exponentially to $0$ (uniformly in $t$) 
as $\sigma\to +\infty$, so we can compute the  integral
in \eqref{eta-D-v} by moving  the line of integration to $\sigma=+\infty$
and see that it tends to $0$ as $\sigma\to +\infty$.
In doing so, we encounter residues at $s=w$ and $s=1$.  The residue at
$s=1$ cancels the constant term $1/\langle 1,1\rangle$. Noting that
$\lambda_s-\lambda_w=-(s-w)(s-(1-w))$, and noting the negative
orientation around $s=w$ of the path integral, the final result is 
\begin{thm}\label{eta-v} For $a>1$ and $\Re(w)>\frac12$,
\[
\eta_a(v_{w,a}) 
=\frac{a^{1-w}(a^w+c^{}_w a^{1-w})}{2w-1}
\]
\end{thm}

\subsection{Computing $\theta_d (v_{w,a})$ for $a\gg_\theta 1$ and $\Re(w)>\frac 12$}
By linearity, it suffices to compute $\delta^{\rm nc}_{z}(v_{w,a})$ when $z$ is a
Heegner point.  Note that $v_{w,a}\in \Eis^{\frac 32-\ve}$ for all
$\ve>0$, so the integral for the pairing $\Eis^{\frac 32-\ve}\times
\Eis^{-1-\ve}$ is absolutely convergent if $\ve$ is sufficiently small.

Using $c^{}_{1-s}E_s=E_{1-s}$ and the spectral expansion, we find
\begin{align}
 & \phantom{=}\delta^{\rm nc}_{z}(v_{w,a}) =\frac{1}{\langle 1,1\rangle} +
\frac{1}{4\pi i} \int_{(\frac 12)}\eta_a E_{1-s}(z) \cdot E_s(z)
\pt\frac{\d s}{\lambda_s-\lambda_w}\notag\\
&=\frac{1}{\langle 1,1\rangle} +
\frac{1}{4\pi i} \int_{(\frac 12)}(a^{1-s}+c^{}_{1-s}a^s) \cdot E_s(z)
\pt\frac{\d s}{\lambda_s-\lambda_w}\notag\\
&=\frac{1}{\langle 1,1\rangle} +
\frac{1}{4\pi i} \int_{(\frac 12)}(a^{1-s} E_s(z)
+a^s E_{1-s}(z)) 
\pt\frac{\d s}{\lambda_s-\lambda_w} \tag{use $c^{}_{1-s}E_s=E_{1-s}$}\\
&=\frac{1}{\langle 1,1\rangle} +
\frac{1}{2\pi i} \int_{(\frac 12)}a^{1-s} E_s(z)
\pt\frac{\d s}{\lambda_s-\lambda_w} \label{delta-v-a}
\end{align}

The computation of the integral requires some care, depending on the height 
of $z$ relative to $a$.  To this end, we proceed as before, moving the line of integration
from $\sigma=\frac 12$ to $\sigma=C$  where $C>1$ (the actual value of
$C$ is immaterial), thereby  
acquiring the contribution of residues at $s=w$ and also at $s=1$ from the Eisenstein series.
The residue of $E_s(z)$ at $s=1$ is $\displaystyle \frac 3\pi = \frac{1}{\langle 1,1\rangle}$,
hence its contribution cancels the constant term and one obtains
\begin{equation}\label{first-step-E}
\delta^{\rm nc}_{z}(v_{w,a}) =
\frac{a^{1-w}E_w(z)}{2w-1}+\frac{1}{2\pi i} \int_{(C)} 
a^{1-s} E_s(z)\pt\frac{\d s}{\lambda_s-\lambda_w}
\end{equation}

The series for $E_s(z)$ with $z=x+iy$ 
is absolutely convergent for $\Re(s)=c>1$ and for
$y\to\infty$ it is asymptotic to $y^s$.   It is obvious that for $\sigma>1$
we have
\[
|E_s(z)| \leqslant \frac 12 \sum_{m,n}{}^{'} \frac {y^\sigma}{|nz+n|^{2\sigma}}
\]
where the dash means that the sum is extended to all pairs $(m,n)$ of coprime
integers. It follows that if $y$ is bounded away from zero (in our case
$y\geqslant \sqrt{3}/2$)  then $|E_s(z)| \ll (\max(1,y^{\sigma})$
uniformly for $\sigma \geqslant c >1$, with the constant involved in the inequality
depending only on $c$.

Therefore, if $y/a < 1$ one may move the line of integration all
the way to $+\infty$,  showing that the integral in question vanishes.
We have proved that in this case
\begin{equation}\label{first-evaluation}
\delta^{\rm nc}_{z}(v_{w,a}) =
\frac{a^{1-w}E_w(z)}{2w-1}\qquad (\Im(z) < a)
\end{equation}
If instead $y/a >1$ the analysis is more complicated.  In this case, we
split the sum for $E_s(z)$ into two components:
\begin{equation}\label{S1+S2}
\begin{split}
E_s(z) &=   \frac 12\, 
\sum_{\atopa{|mz+n|^2 >y/a}{{\rm GCD}(m,n)=1}}\frac {y^s}{|mz+n|^{2s}} 
+   \frac 12\, 
\sum_{\atopa{|mz+n|^2 \leqslant y/a}{{\rm GCD}(m,n)=1}} \frac {y^s}{|mz+n|^{2s}}\\
& = {\Sigma}_1 + {\Sigma}_2.
\end{split}
\end{equation}

The evaluation of the integral $\displaystyle \frac 1{2\pi i} \int_{(C)} a^{1-s}\,
{\Sigma}_1 \pt \frac{\d s}{\lambda_s-\lambda_w}$ can be done as before
by letting $C\to +\infty$, obtaining 
\begin{equation}\label{eval-S1}
\frac 1{2\pi i} \int_{(C)} a^{1-s}\,
{\Sigma}_1\pt \frac{\d s}{\lambda_s-\lambda_w} = 0.
\end{equation}
To deal with the integral involved in ${\Sigma}_2$ we note that the
sum involved is a finite sum. For
$|cz+d|^2 \ne y/a$
we can integrate term-by-term
and move the line of integration backwards all the way to $\to -\infty$, 
with the limit of the integral being $0$.  
In doing this we encounter two residues at $s=w$ and $s=1-w$, and conclude that
\begin{equation}\label{eval-S2}
\begin{split}
\frac 1{2\pi i} \int_{(C)} a^{1-s}\,
{\Sigma}_2 \pt\frac{\d s}{\lambda_s-\lambda_w} = 
-\frac 12 \pt\frac{a^{1-w}}{2w-1}
\sum_{\atopa{|mz+n|^2\leqslant y/a}{{\rm GCD}(m,n)=1}} \frac {y^w}{|mz+n|^{2w}}\\
+\frac 12\pt\frac{a^{w}}{2w-1} 
\sum_{\atopa{|mz+n|^2\leqslant y/a}{{\rm GCD}(m,n)=1}} \frac {y^{1-w}}{|mz+n|^{2-2w}}.
\end{split}
\end{equation}

We  have proved:
If $\Im(z) < a$ then 
\begin{equation}\label{large a}
\delta_z^{\rm nc}(v_{w,a}) = \frac{a^{1-w} E_w(z)}{2w-1}.
\end{equation}
while if $y=\Im(z)\geqslant a$ then
\begin{equation}\label{small a}
\begin{split}
\delta_z^{\rm nc}(v_{w,a}) = \frac{a^{1-w} E_w(z)}{2w-1}-\frac 12 \pt\frac{a^{1-w}}{2w-1}\sum_{\atopa{|mz+n|^2 < y/a}{{\rm GCD}(m,n)=1}} \frac {y^w}{|mz+n|^{2w}}\\
+\frac 12\pt\frac{a^{w}}{2w-1} \sum_{\atopa{|mz+n|^2 < y/a}{{\rm GCD}(m,n)=1}}  \frac {y^{1-w}}{|mz+n|^{2-2w}}.
\end{split}
\end{equation}

The condition $|mz+n|^2 < y/a$ is rather restrictive on the pair $(m,n)$.
In fact, if $m \ne 0$ then  $|mz+n|^2 \geqslant y^2$ and the condition
$(my)^2 < y/a$ implies $y<1/(m^2 a)$, which implies  $m=\pm 1$ or $m=0$.
If $m=\pm 1$ then $|mz+n |^2 = (x\mp n)^2 + y^2 < y/a$, hence
$\frac 14 +   y^2 < y$ because $|x\mp n|\geqslant \frac 12$
and $a>1$ by hypothesis, but this is impossible.  We conclude that
$m=0$, $n=\pm 1$ and we get the much simplified formula
\begin{equation}\label{Dirac-delta-z}
\delta_z^{\rm nc}(v_{w,a}) = \frac{a^{1-w} E_w(z)}{2w-1}
-\frac{a^{1-w}y^w - a^w y^{1-w}}{2w-1}
\end{equation}
in the remaining range $1 < a < y$.
By linearity, this extends to the computation of $\theta_d(v_{w,a})$
for Eisenstein-Heegner distributions $\theta_d$, and to real-linear
combinations of such. We summarize these computations as a theorem: 
\begin{thm}\label{theta-v}
Let $a>1$, and let $d<0$ be a fundamental discriminant.  Then
\begin{equation}\label{tDv}
\theta_d(v_{w,a}) = \frac{1}{2w-1} \left\{a^{1-w}\pt
\theta_d E_w- R_w(d,a)\right\}
\end{equation}
where 
\begin{equation}\label{RD}
R_w(d,a)=
\sum_{\atopa{x+i y\in H_d}{y>a}} (a^{1-w}y^w-a^w y^{1-w})
\end{equation}
\end{thm}


\subsection{Computing $\eta_a(u_{\theta, w})$ for $a \gg_\theta 1$ and
$\Re(w)>\frac 12$}

Let $\theta$ be a finite real-linear combination of
Eisenstein-Heegner distributions $\theta_d$. 
\begin{thm}\label{eta-u}
\begin{equation}\label{eta-u=theta-v}
\eta_a(u^{}_{\theta,w})=\theta(v_{w,a}).
\end{equation}
\end{thm}
\proof One computes
\begin{align*}
\eta_a(u^{}_{\theta,w}) &= \frac{\eta_a(1)\theta(1)}{\langle 1,1\rangle} +
\frac {1}{4\pi i} \int_{(\frac 12)}\theta E_{1-s}\cdot \eta_a E_s
\pt \frac{\d s}{\lambda_s-\lambda_w}\\
&=\frac{\theta(1)}{\langle 1,1\rangle}+
\frac {1}{4\pi i} \int_{(\frac 12)}\theta E_{1-s}\cdot (a^s+c^{}_s a^{1-s})
\pt \frac{\d s}{\lambda_s-\lambda_w}\\
&=\frac{\theta(1)}{\langle 1,1\rangle}
+\frac {1}{4\pi i} \int_{(\frac 12)}(\theta E_{1-s}\cdot a^s + \theta E_s \cdot a^{1-s})
\pt \frac{\d s}{\lambda_s-\lambda_w}\tag{because $c^{}_sE_{1-s}=E_s$} \\
&=\frac{\theta(1)}{\langle 1,1\rangle}
+\frac {1}{2\pi i} \int_{(\frac 12)}\theta E_{1-s}\cdot a^s
\pt \frac{\d s}{\lambda_s-\lambda_w}. \tag{by changing $s\to 1-s$ in one term}
\end{align*}
The theorem follows by linearity and equation \eqref{delta-v-a}. \qed


\subsection{Computing $\theta(u^{}_{\theta,w})$ 
for $a > 1$ and $\Re(w)>\frac 12$}

The outcome here does not admit much simplification, in contrast to
the other cases. That is, from the spectral expansions of subsections
\ref{Heegner-distributions} and \ref{solving}, via the $\Eis^1\times
\Eis^{-1}$ pairing of \ref{Eisenstein-Sobolev}, we obtain
\begin{thm}\label{theta-u}
For $\theta$ a finite real-linear combination of Eisenstein-Heegner
distributions $\theta_d$,
\begin{equation}\label{theta-theta}
\theta(u^{}_{\theta,w})
=
\frac{|\theta(1)|^2}{\langle 1,1\rangle\cdot (\lambda_1-\lambda_w)}
+
\frac{1}{4\pi i}\int_{(\frac 12)}
\left|\theta E_s\right|^2
\frac{\d s}{\lambda_s-\lambda_w}.
\end{equation}
\end{thm}


\subsection{Rewriting the determinant
condition}\label{two-condition-explication}

First, by Theorem \ref{eta-u=theta-v} the determinant-vanishing
condition \eqref{resolvent} becomes
\[
\eta(v_{w,a})\pt\theta(u^{}_{\theta,w}) - \eta_a(u^{}_{\theta,w})^2=0
\]
In view of Theorems \ref{eta-v}, \ref{theta-v},  \ref{theta-u}, 
\begin{cor}\label{key-cor}
Let $\theta=\sum_d \nu_d \theta_d$ be a finite real-linear combination
of Eisenstein-Heegner distributions $\theta_d$ with $d<-4$. 
For all $a>1$,  all $w$ with $\Re(w)>\frac 12$
and off $(\frac 12,1]$,
\begin{align}\label{master equation}
a^{1-w}(a^{w}+c^{}_w a^{1-w})&\times \Big(
\frac{|\theta(1)|^2}{\langle 1,1\rangle\cdot \lambda_1-\lambda_w)}
+ \frac{1}{4\pi i}\int_{(\frac 12)}
\left|\theta E_s\right|^2  
 \frac{\d s}{\lambda_s-\lambda_w}\Big)\notag\\
&\neq \frac{1}{2w-1} \left( a^{1-w} \theta E_w-R_w(\theta,a) \right)^2
\end{align}
where
\begin{equation}\label{L}
\theta E_s=\sum_d \nu_d\pt
\bigg(\frac{\sqrt{|d|}}{2}\bigg)^s \pt \frac{\zeta(s)}{\zeta(2s)}\pt L(s,\chi_d)
\end{equation}
and where
\begin{align}
R_w(\theta,a)=\sum_d \nu_d
\sum_{\atopa{x+i y\in H_d}{y>a}} (a^{1-w}y^w-a^w y^{1-w}).\tag*{\eqref{RD}}
\end{align}
\end{cor}
\proof In $\Re(w)>\frac12$, $u_{\theta,w}$ and $v_{a,w}$ are in
$\Eis^1$. The vanishing condition
$\eta(v_{w,a})\pt\theta(u^{}_{\theta,w}) -
\eta_a(u^{}_{\theta,w})^2=0$ is thus necessary and sufficient for some
non-zero linear combination of $u_{\theta,w}$ and $v_{a,w}$ to be an
eigenfunction for the self-adjoint semibounded operator $\widetilde
S_{\theta,a}$, with eigenvalue $\lambda_w=w(1-w)$. Such an eigenvalue
must be real and satisfy $0\leqslant \lambda_w$. \qed


\subsection{Meromorphic continuation and location of zeros}
\label{relatively-elementary-mero-cont}

Here we give a direct, relatively elementary argument for meromorphic
continuation of the two-by-two determinant above. Theorem
\ref{mero-contn-N-valued-integral} in section
\ref{mero-contn-spectral-synthesis-integrals} gives several stronger
results by less elementary means, useful in the discussion of
unconditional results on zero spacing in Section \ref{unconditional-spacing}.
For brevity, write
\begin{align}\label{F(a,w)}
F(a,w):=(a^{w}&+c^{}_w a^{1-w})\pt
\Big(\frac{|\theta(1)|^2}{\langle 1,1\rangle\cdot
(\lambda_1-\lambda_w)}
+\frac{1}{4\pi i}\int_{(\frac 12)}
\left|\theta E_s\right|^2  
 \frac{\d s}{\lambda_s-\lambda_w}\Big)
\notag \\
& -\frac {a^{1-w}\,\theta(E_w)^2}{2w-1}
\end{align}
for the determinant above, where for the time being $\Re(w)>\frac 12$
and $a$ is real with $a>1$. 

The analytic continuation of $F(a,w)$ beyond the line $\Re(w)=\frac 12$ is 
easily accomplished.
Except for a possible simple pole at $w=1$, the function $s\to\theta
E_s$  is holomorphic for $\Re(s) > \frac 12 - C/\log(2+|s|)$ for some
absolute positive constant $C$. 
Hence, when $\frac 12 < \Re(w) < \frac 12 + C/\log(2+|w|)$ we can evaluate the integral in \eqref{F(a,w)} as follows :
\begin{align}
\frac{1}{4\pi i}\int_{(\frac 12)}|\theta E_s|^2
\pt\frac{\d s}{\lambda_s-\lambda_w}\notag\\
&\hskip -1.0in =\frac{1}{4\pi i}\int_{(\frac 12)}
\Big( \theta E_s\cdot \theta E_{1-s}- \theta E_w\cdot \theta E_{1-w}\Big)
\pt\frac{\d s}{\lambda_s -\lambda_w}\notag\\
&\hskip -1.0in \quad +\theta E_w\cdot  \theta E_{1-w}
\pt\frac{1}{4\pi i}\int_{(\frac 12)}
\frac{\d s}{\lambda_s-\lambda_w}\,.\label{F}
\end{align}
The left-hand side of this equation is a holomorphic function for
$\Re(w)>\frac 12$. For $\Re(w)>\frac 12$ the last integral is evaluated by the calculus of residues,
moving the line of integration to $\Re(s)\to +\infty$:
\begin{equation}\label{integral}
\frac{1}{4\pi i}\int_{(\frac 12)}
\frac{\d s}{\lambda_s-\lambda_w} = \frac{1}{2(2w-1)}.
\end{equation}
Therefore, for $\frac 12 <\Re(w)<\frac 12+C/(2+|w|)$ we have
\begin{align}
F(w,a)
= \;&
(a^w+c_w a^{1-w})\pt\frac{|\theta(1)|^2}{\langle 1,1\rangle\cdot (\lambda_1-\lambda_w)}
\notag\\
& + (a^w+c_w a^{1-w})\pt
\frac{1}{4\pi i}\int_{(\frac 12)}\Big(\theta E_s\cdot \theta E_{1-s}
- \theta E_w\cdot\theta  E_{1-w}\Big) \pt
\frac{\d s}{\lambda_s-\lambda_w}\notag\\
& \qquad +(a^w+c_w a^{1-w})\pt\frac{\theta E_w\cdot \theta E_{1-w}}{2(2w-1)}
-\frac{ a^{1-w}\theta (E_w)^2}{2w-1}
\end{align}
Using the functional equation $c_w\theta E_{1-w}=\theta E_w$
simplifies this into
\begin{align}
F(w,a)
= \;&
(a^w+c_w a^{1-w})\pt\frac{|\theta(1)|^2}{\langle 1,1\rangle\cdot (\lambda_1-\lambda_w)}
\notag\\
& + (a^w+c_w a^{1-w})\pt
\frac{1}{4\pi i}\int_{(\frac 12)}\Big(\theta E_s\cdot \theta E_{1-s}
- \theta E_w\cdot\theta  E_{1-w}\Big) \pt
\frac{\d s}{\lambda_s-\lambda_w}\notag\\
& \qquad +(a^w-c_w a^{1-w})\pt\frac{\theta E_w\cdot \theta E_{1-w}}{2(2w-1)}
\end{align}\label{newF}
This formula, so far proved for $\frac12 <\Re(w)<\frac 12+C/\log(2+|w|)$, extends to 
an open neighborhood
of the line $\Re(w)=\frac 12$ in the complex plane $w\in\mathbb C$,
because now the integral is well defined there
as a continuous function of $w$.  Indeed, the numerator of the integrand vanishes
as $w\to s$ with $\Re(s)=\frac 12$ at least to the first order when $s\ne \frac 12$  and
at at least to the second order when $s=\frac 12$, while the growth of
$\theta E_s\cdot\theta E_{1-s}$ is of order not more than $|s|^{1-\delta}$
there for some fixed $\delta>0$ if the neighborhood is sufficiently small, hence the integral 
is absolutely convergent. 

\begin{thm}\label{F-to-G-functional-eq}
With $a>1$, let 
\begin{equation}\label{G}
G(w,a) := \frac {F(w,a)}{a^w+c_w a^{1-w}}
\end{equation}
Then 
$G(w,a)$ is a meromorphic function in the whole complex $w$-plane 
and satisfies the functional equation
\begin{equation}\label{funct-G}
G(w,a) = G(1-w,a).  
\end{equation}

\end{thm}
\proof A simple computation using the expansion \eqref{newF} shows
that the stated functional equation holds in an open neighborhood of
the critical line. Due to the importance of this symmetry, we carry
out this computation in detail.
\begin{align}
G(w,a)
= \;&
\frac{|\theta(1)|^2}{\langle 1,1\rangle\cdot (\lambda_1-\lambda_w)}
\notag\\
& + 
\frac{1}{4\pi i}\int_{(\frac 12)}\Big(\theta E_s\cdot \theta E_{1-s}
- \theta E_w\cdot\theta  E_{1-w}\Big) \pt
\frac{\d s}{\lambda_s-\lambda_w}\notag\\
& \qquad +\frac{a^w-c_w a^{1-w}}{a^w+c_wa^{1-w}}\pt\frac{\theta
E_w\cdot \theta E_{1-w}}{2(2w-1)}
\end{align}
The first two summands are indeed invariant under $w\to 1-w$, as is
$\theta E_w\cdot \theta E_w$ in the third summand. Finally, under $w\to
1-w$, using $c_w\cdot c_{1-w}=1$, the part
\[
\frac{a^w-c_w a^{1-w}}{a^w+c_wa^{1-w}}\pt
\frac{1}{2w-1}
\]
of the third summand becomes
\[
\frac{a^{1-w}-c_{1-w} a^w}{a^{1-w}+c_{1-w}a^w}\pt
\frac{-1}{2w-1}
\;=\;
\frac{c_w a^{1-w}-a^w}{c_w a^{1-w}+a^w}\pt
\frac{-1}{2w-1}
\;=\;
\frac{a_w-c_w a^{1-w}}{a^w+c_w a^{1-w}}\pt
\frac{1}{2w-1}
\]
giving the claimed invariance. The conclusion of the theorem follows
by analytic continuation. \qed

\begin{cor}\label{nonvanishing_off_line}
The only zeros of the function $G(w,a)$ defined in equation
\ref{G} are on $\Re(w)=\frac12$ and $[0,1]$.  
\end{cor}
\proof Corollary \ref{key-cor} shows that $G(w,a)$ cannot vanish in
$\Re(w)>\frac12$ except possibly on $(\frac12,1]$, because otherwise
$\lambda_w=w(1-w)$ would be an eigenvalue for a non-negative
self-adjoint operator. Then the symmetry of \ref{funct-G} shows
non-vanishing in $\Re(w)<\frac12$ except possibly on $[0,\frac12)$.  
\qed


\subsection{An important remark}  The preceding considerations also apply to
a general real-linear combination $\theta=\sum_\nu b_\nu\cdot
\delta^\nc_{z_\nu}$ of Eisenstein-Dirac distributions $\delta^\nc_{z_\nu}$.
With the simplifying assumption $\sum b_\nu E_s(z_\nu)=0$, we have
\begin{align}\label{master eisenstein}
&a^{1-w}(a^{w}+c^{}_w a^{1-w})\times \frac{1}{4\pi i}\int_{(\frac 12)}
\left| \sum b_\nu E_s(z_\nu)\right|^2
\pt \frac{\d s}{\lambda_s-\lambda_w}\notag\\
& =\frac{1}{1-2w} \left\{\sum b_\nu \left[ a^{1-w}E_w(z_\nu)- 2 \sqrt{a\Im(z_\nu)}
\sinh\left(\frac 12(1-w)\log^+\frac{\Im(z_\nu)}{a}\right)\right]\right\}^2
\end{align}
where $\log^+x = \max(\log x,0)$.

If $a > \max\Im(\zeta_\nu)$ then for almost all $a$ the vanishing of
$\sum b_\nu E_w(z_\nu)$ implies the vanishing of the left-hand term, which means
either $a^{1-w}(a^{w}+c^{}_w a^{1-w})=0$, which vanishes only when
$\Re(s)=\frac 12$, or
\[
\frac{1}{4\pi i}\int_{(\frac 12)}
\left| \sum b_\nu E_s(z_\nu)\right|^2
 \frac{\d s}{\lambda_s-\lambda_w}=0.
\]
Thus if zeros on the critical line of a
function $\sum_\nu b_\nu E_s(z_\nu)$ with $\sum b_\nu=0$ had a spectral interpretation 
for some $a=a_0>\max \Im(z_\nu)$ this would be so for all $a>a_0$, implying that the corresponding eigenvalues 
$w(1-w)$ would be independent of $a$.  For $a\to\infty$ this is analogous to
the condition formulated by Colin de Verdi\`ere in the special case $\{\nu\}=\{1\}$, $b_1=1$,
$z_1=\frac 12+i\frac{\sqrt 3}{2}$, and, tentatively, suggested by him as a spectral 
intepretation of the zeros of $\zeta_{Q(\sqrt{-3})}(s)$ on the critical line.


\subsection{Computing $\theta E_s$ in a special case}  The fundamental 
discriminants $d$ are the odd squarefree numbers $d=m$ with
$m\equiv 1 \pmod 4$ or numbers of the type $d=4m$ with 
$m$ squarefree and $m\equiv 2,3 \pmod 4$.  To $m$, one associates the real primitive
character $\chi_m$ given by
\begin{equation}\label{Kron}
\chi_m(n) =
\begin{cases} 
\displaystyle{\left(\frac mn\right)}\quad & m\equiv 1 \  (\bmod 4)\\
 \\
\displaystyle{\left(\frac {4m}n\right)}\quad & m\equiv 2,3 \  (\bmod 4)\\
\end{cases}
\end{equation}
where on the right-hand side we have the Kronecker symbol.
Then $\zeta(s) L(s,\chi_m)$ is the zeta function of the quadratic
field $\mathbb Q(\sqrt{d})$.

A precise asymptotic formula for a certain simple linear combination
of quadratic $L$-functions has been obtained in the paper \cite{GH} of
Goldfeld and Hoffstein. We recall {\it verbatim\/} their Theorem (1),
where their reference (0.6) is our equation \eqref{Kron}. (Compare
also \cite{VT}.)

\noindent{\bf Theorem (1)}\ \  {\it Let $\ve>0$ be fixed. Let $\chi_m$ be defined as
in (0.6).  Then there exist analytic functions $c(\rho)$ and $c^*_{\pm}(\rho)$
with Laurent expansion 
$c(\rho)=c_{\frac 12}/(\rho-\frac 12)+c'_{\frac 12}+O(\rho-\frac 12)$,
$c^*_{\pm}(\rho)=-c_{\frac 12}$ such that
\begin{align*}
\sum_{\atopa{1<-m<x}{m\ \square -free}}&L(\rho,\chi_m) \\
&\hskip -0.4in=\begin{cases}
c(\rho)\pt x+O(x^{\frac 12+\ve})\qquad  & if\ \ {\rm Re}(\rho)\geq 1\\
c(\rho)\pt x+c^*_\pm(\rho)\pt x^{\frac 32-\rho}+ O(x^{\theta+\ve})\qquad & if\ \ {\rm Re}(\rho)\neq\frac 12, \ \frac 12\leq {\rm Re}(\rho)<1\\
c_{\frac 12}\pt x \log x +(c'_{\frac 12}+{c^*}'_{\hskip -3pt \pm \frac 12}
-c_{\frac 12})\pt x
+O(x^{\frac{19}{32}+\ve})\qquad  & if\ \ \rho=\frac 12.\\
\end{cases}
\end{align*}
Here 
\begin{align*}
c(\rho) &= \frac 34 (1-2^{-2\rho})\pt\zeta(2\rho) \prod_{p\ne 2}
(1-p^{-2}-p^{-2\rho-1}+p^{-2\rho-2}),\\
c_{\frac 12} &=\frac 3{16} \prod_{p\ne 2}(1-2p^{-2}+p^{-3})
\end {align*}
and
\[
\theta = \begin{cases}
\displaystyle {\frac 12} \qquad & if\ \ \displaystyle{{\rm Re}(\rho)>\frac{-5+\sqrt{193}}{12}}\\
 \\
\displaystyle {\frac{19+3{\rm Re}(\rho)-6{\rm Re}(\rho)^2}{24+16{\rm Re}(\rho)}}
\qquad & if\ \ \displaystyle{\frac 12 \leqq \rm Re(\rho)\leqq\displaystyle{\frac{-5+\sqrt{193}}{12}}}
\end{cases}
\]
and all O-constants depend at most on $\rho$, $\ve$.
}

The authors do not give the dependence on $\rho$ in the
proportionality factor involved in the symbol  $O(...)$, but there must be
a function $\omega(x)$ slowly increasing to $\infty$ such that the 
estimates remain uniform in $\rho$ as long as $|\Im(\rho)|<\omega(x)$.

What is of interest to us is not the sum $\sum L(\rho,\chi_m)$ (with $m<0$)
but rather the sum divided by $\zeta(2\rho)$.  This yields the following result.

Let 
\begin{equation}\label{P(s)}
A(s)= \frac 34 (1-2^{-2s})\prod_{p\ne 2}(1-p^{-2}-p^{-2s-1}+p^{-2s-2})
\end{equation}

\begin{thm}\label{average-L} There is a function $\omega(x)$, slowly increasing
to $\infty$ as $x\to\infty$, such that
\begin{equation}\label{main-asymp}
\sum_{\atopa{1<-m<x}{m\ \square -free}} \zeta(2s)^{-1}L(s,\chi_m)
= A(s) x + B(s) x^{\frac 32-s} + O(x^{\frac{19}{32}+\ve})
\end{equation}
holds, uniformly for $\frac 12 \leqslant \Re(s) \leqslant 1$ and $ |\Im(s)| < \omega(x)$.
\end{thm}
The function $B(s)$ is more complicated to describe but it is holomorphic for 
$\frac 12\leqslant \Re(s)$. (Again, compare \cite{VT}.)


\section{Meromorphic continuations of spectral integrals}

Here we prove meromorphic continuation results stronger than the
scalar meromorphic continuation in Theorem \ref{F-to-G-functional-eq}.
The immediate goal is to prove that, for suitable $\theta\in \Eis^{-1+\ve}$, solutions
$u^{}_{\theta,w}$ of equations $(-\Delta-\lambda_w)u=\theta$ expressed by
spectral expansions 
\[
u^{}_{\theta,w}
=
\frac{\langle \theta,1\rangle\cdot 1}{(\lambda_1-\lambda_w)\cdot
\langle 1,1\rangle}
+\frac{1}{4\pi i}\int_{(\frac 12)} \frac{\mathcal E \theta(s) \cdot
E_s}{\lambda_s-\lambda_w}\;\d s
\]
convergent in $\Eis^{1+\ve}$, at first valid only in $\Re(w)>\frac 12$,
meromorphically extend to the whole complex plane, as function-valued
functions. However, the meromorphic continuations do not lie
in
$\Eis^{1+\ve}$, but only in a larger space $M$, large enough to include
Eisenstein series. On the critical 
line, we find that $u^{}_{\theta,w}$ stays in $\Eis^{1+\ve}$ only for
$\mathcal E\theta(w)=0$. By Corollary \ref{discrete-spectrum-if-any},
the only {\itshape possible} discrete spectrum $\lambda_w>\frac 14$ of
$\widetilde S_\theta$ occurs among zeros of $\mathcal E\theta$, this
will show that eigenfunctions in the discrete spectrum, if any, are
meromorphic continuations of these spectral synthesis integrals.
However, this still does not prove that zeros $w$ of $\mathcal
E\theta$ on $\Re(w)=\frac 12$ gives eigenvalues $w(1-w)$ of
$\widetilde S_\theta$, since $u^{}_{\theta,w}$ is not in the domain of
$\widetilde S_\theta$ unless, additionally, $\theta
u^{}_{\theta,w}=0$.
In fact, we discuss meromorphic continuations of images $\Phi
u^{}_{\theta,w}$ under continuous linear maps $\Phi:M\to N$ for
quasi-complete, locally convex topological vector spaces $N$.


\subsection{Vector-valued integrals}\label{Gelfand-Pettis}
We recall some standard results about vector-valued integrals,
mostly without proofs. Original sources are \cite{Gelfand} and
\cite{Pettis}, for which \cite{Rudin} offers a reasonable exposition.
See also \cite{Garrett} chapter 14 for exposition and proofs more
tightly aimed at applications such as those here.
Let $V$ be a topological vectorspace over $\mathbb C$, $f$ a measurable
$V$-valued function on a measure space $X$. 
A {\itshape Gelfand-Pettis integral} of $f$ is a vector $I_f\in V$ so that
$\lambda(I_f)=\int_X \lambda\circ f$ for all $\lambda\in V^*$.
If it exists and is unique, this vector $I_f$ is denoted $\int_X f$.
In contrast to construction of integrals as limits, this
characterization is a property no reasonable notion 
of integral would lack. Since this property is an irreducible minimum, this
characterizes a {\itshape weak integral}.

Uniqueness of the integral is immediate when $V^*$ separates points on
$V$, as for locally convex $V$, by the Hahn-Banach theorem. Similarly, linearity
of $f\to I_f$ follows when $V^*$ separates points. Thus, the issue is
existence.

The functions we integrate are relatively nice: compactly-supported and
continuous, on measure spaces with finite, positive, regular Borel
measures. In this situation, all the $\mathbb C$-valued integrals 
$\int_X\lambda\circ f$
exist for elementary reasons, being integrals of
compactly-supported $\mathbb C$-valued continuous functions on a compact set
with respect to a finite regular Borel measure. 

A topological vector space is {\itshape quasi-complete} or {\itshape locally
complete} if every {\itshape bounded} (in the general topological
vector space sense) Cauchy net is convergent. It is known (for
example, see \cite{BTVS}) that

\vfill\break

\begin{lem} In a quasi-complete, local convex topological vector space,
the convex hull of a compact set has compact closure.
\end{lem}

The latter property ensures existence of certain Gelfand-Pettis
integrals:

\begin{thm}  Let $X$ be a locally compact Hausdorff topological space with a
{\itshape finite}, positive, regular Borel measure. Let $V$ be a locally
convex topological vectorspace in which the {\itshape closure of the convex
hull of a compact set is compact}. Then {\itshape continuous,
compactly-supported} $V$-valued functions $f$ on $X$ have
Gelfand-Pettis integrals. Further,
\[
\int_X\,f \;\in\; \text{meas}(X)\cdot 
\Big(\text{closure of convex hull of } f(X)\Big)
\]
is the basic estimate substituting for estimating a Banach-space norm
of an integral by the integral of the norm of the integrand.
\end{thm}

The legitimacy of passing continuous operators inside such integrals
is an easy corollary:

\begin{cor} Let $T:V\to W$ be a continuous linear map of locally convex
topological vectorspaces, where convex hulls of compact sets in $V$
have compact closures. Let $f$ be a continuous,
compactly-supported $V$-valued function on a finite regular measure
space $X$. Then the $W$-valued function $T\circ f$ has a
Gelfand-Pettis integral, and
$T\Big(\int_X f\Big)=\int_X T\circ f$.
\end{cor}

\proof {\itshape (of corollary)} To verify that the left-hand side of the
asserted equality fulfills the requirements of a Gelfand-Pettis
integral of $T\circ f$, we must show that
\[
\lambda\Big(\hbox{left-hand side}\Big)\;=\; \int_X \lambda\circ(T\circ f)
\]
for all $\lambda\in W^*$. Starting with the left-hand side,
\[
\lambda\Big(T\Big(\int_X f\Big)\Big)
= (\lambda\circ T) \Big(\int_X f\Big)
= \int_X (\lambda\circ T) \circ f
= \int_X \lambda\circ ( T \circ f)
\]
proving that $T\big(\int_X f\big)$ is a weak integral of $T\circ f$. \qed


\subsection{Holomorphic vector-valued functions}
\label{holomorphic-vector-valued}
Now we recall basic facts about holomorphic vector-valued functions,
mostly without proof, for which we refer to the original source
\cite{Gro}, or expositions such as \cite{Rudin} or \cite{Garrett},
chapter 14. Existence and properties of vector-valued integrals
are ingredients in the proofs of the assertions below.

A function $f$ on an open set $\Omega\subset \mathbb C$ taking values in a
quasi-complete, locally convex topological vector space $V$ is
(strongly) {\itshape complex-differentiable} when  
\[
\lim_{z\to z_o}\,\frac{1}{z-z_o}\cdot \big(f(z)-f(z_o)\big) 
\]
exists (in $V$) for all $z_o\in \Omega$, where $z\to z_o$ specificially means
for {\itshape complex} $z$ approaching $z_o$. The function $f$ is
(strongly) {\itshape analytic} when it is locally expressible as a
convergent power series with coefficients in $V$. The function $f$ is
{\itshape weakly holomorphic} when the $\mathbb C$-valued functions $\lambda\circ
f$ are holomorphic for all $\lambda$ in $V^*$.

\begin{thm} For $V$ a locally convex quasi-complete topological vector space,
{\itshape weakly} holomorphic $V$-valued functions $f$ are {\itshape strongly}
holomorphic, in the following senses. First the usual Cauchy-theory
integral formulas apply: 
\[
f(z) = \frac{1}{2\pi i}\, \int_\gamma \, \frac{f(\zeta)}{\zeta - z}
\, \d\zeta
\]
with $\gamma$ a closed path around $z$ having winding number
$+1$. Second, the function $f(z)$ is infinitely differentiable, in fact 
{\itshape strongly analytic}, that is, expressible as a convergent power series
$f(z)=\sum_{n\ge 0} \, c_n \, (z-z_o)^n$ with coefficients $c_n\in V$
given by Gelfand-Pettis integrals echoing 
Cauchy's formulas: 
\[
c_n = \frac{f^{n}(z_o)}{n!} 
= 
\frac{1}{2\pi i}\,\int_\gamma\, \frac{f(\zeta)}{(\zeta - z_o)^{n+1}}\,
\d \zeta 
\]
\end{thm}

The appropriate vector-valued notion of {\itshape meromorphy} is completely
parallel to the scalar-valued version: a $V$-valued function $f$
defined on a punctured neighborhood $N$ of $z_o$ is {\itshape meromorphic}
at $z_o$ when there is a positive integer $n$ such that
$(z-z_o)^m\cdot f(z)$ has an extension to a $V$-valued {\itshape
holomorphic} function on $N\cup\{z_o\}$.

Fix a non-empty open subset $\Omega$ of $\mathbb C$. Let $V$ be quasi-complete,
locally convex, with topology given by seminorms $\{\nu\}$. The space
${\mathrm {Hol}}(\Omega,V)$ of holomorphic $V$-valued functions on $\Omega$ has the
natural topology given by seminorms
$\mu^{}_{\nu,K}(f)=\sup_{z\in K} \nu(f(z))$ for compacts $K\subset
\Omega$, seminorms $\nu$ on $V$.
This is the analogue of the sups-on-compacts seminorms on
scalar-valued holomorphic functions, and there is the analogous
corollary of the vector-valued Cauchy formulas:

\begin{cor} ${\mathrm {Hol}}(\Omega,V)$ is locally convex, quasi-complete.
\end{cor}

A $V$-valued function $f(z,w)$ on a non-empty open subset $\Omega\subset
\mathbb C^2$ is {\itshape complex analytic} when it is locally expressible as a
convergent power series in $z$ and $w$, with coefficients in $V$. The
two-variable version of the above discussion of power series with
coefficients in $V$ succeeds without incident.

Again by the vector-valued form of Cauchy's integral formulas:

\begin{cor} Let $V$ be quasi-complete, locally convex. Let $f(z,w)$ be
complex-analytic $V$-valued in two variables, on a domain
$\Omega_1\times\Omega_2\subset \mathbb C^2$. Then the function
$w\longrightarrow (z\to f(z,w))$ is a holomorphic
${\mathrm {Hol}}(\Omega_1,V)$-valued function on $\Omega_2$. 
\end{cor}


There is a vector-valued version of Abel's theorem on convergent power
series in one complex variable, proven in similar fashion:

\begin{prop} Let $c_n$ be a bounded sequence of vectors in a locally convex
quasi-complete topological vector space $V$. Then on $|z|<1$ the series
$f(z)=\sum_nc_nz^n$ converges and gives a {\itshape holomorphic} $V$-valued
function. That is, the function is infinitely-many-times
complex-differentiable.
\end{prop}


\subsection{Spaces $M$ of moderate-growth functions}

The space of moderate-growth continuous functions on $\GH$ is
\[
M = \{f\in C^0(\GH): \sup_{\Im(z)\ge \sqrt{3}/2}
y^{-r}\cdot |f(x+iy)|<\infty,\;\text{for some $r\in \mathbb R$}\}
\]
The correct topology is as a strict inductive limit of Banach subspaces 
$M^r_0$ each obtained as a completion of $C^0_c(\GH)$, with
respect to norms
\[
|f|_{M^r_0} = \sup_{\Im(z)\ge \sqrt{3}/2}y^{-r}\cdot |f(x+iy)|
\]
Thus, $\lim_{y\to \infty}y^{-r}|f(x+iy)|=0$ for $f\in M^r_0$. That is,
the {\itshape set} $M$ is $M=\bigcup_r M^r_0$, but the
topology is perhaps not quite as expected. Nevertheless, $M$ is a
strict colimit (in the locally convex category) of Banach spaces, so
is an LF-space, so is quasi-complete and locally convex.


\subsection{Pre-trace formula and $\Eis^{1+\ve}\subset M$}

From the Sobolev imbedding at the end of \ref{Sobolev-imbedding}, we
have $\Eis^{1+\ve}\subset C^0(\GH)$ for all $\ve>0$. Certainly $E_s\in
C^0(\GH)$ away from poles, and $s\to E_s$ is a meromorphic
$C^0(\GH)$-valued function. Thus, the Fr\'echet space $C^0(\GH)$
(with seminorms given by suprema on compact subsets) contains both
$\Eis^{1+\ve}$ and Eisenstein series. However, the space $M$
is much smaller, contains Eisenstein series, and $s\to E_s$ is a
meromorphic $M$-valued function. Thus, we want

\begin{lem} $\Eis^{1+\ve}\subset M$, for all $\ve>0$.
\end{lem}

\proof A slightly more refined form of the pre-trace formula from
\ref{standard-estimate} is
\[
\sum_{|s_F|\le T} |F(z)|^2 + \frac{1}{4\pi}\int_{-T}^T
|E_{\frac 12+it}(z)|^2\;dt \;\;\ll\;\; T^2 + T\cdot \Im(z)
\]
as $\Im(z)\to +\infty$, where $F$ runs over an orthonormal basis for
cuspforms, and $F$ has eigenvalue $s_F(1-s_F)$ for the invariant
Laplacian. Thus, certainly, 
\[
\frac{1}{4\pi}\int_{-T}^T |E_{\frac 12+it}(z)|^2\;dt \;\;\ll\;\; T^2 +
T\cdot \Im(z)
\]
as $\Im z\to \infty$. As earlier, this asserts that the functional
$\delta_z^{\mathrm nc}$ given by $\delta^{\mathrm nc}_z f=f(z)$ is in $\Eis^{-1-\ve}$, a
weaker assertion than what we know to be true. However, the
integration by parts does also yield an estimate for the $\Eis^{-1-\ve}$
norm depending on height $\Im(z)$, namely,
$|\delta^{\mathrm nc}_z|_{\Eis^{-1-\ve}} \ll_\ve y$ as $y\to\infty$, for all $\ve$.
For $u\in \Eis^{1+\ve}$, by the Cauchy-Schwarz-Bunyakowsky
inequality extended to the pairing $\Eis^{1+\ve}\times
\Eis^{-1-\ve}\to \mathbb C$, 
\[
|u^{}_{\theta,w}(z)| = |\delta^{\mathrm nc}_z u^{}_{\theta,w}|
\le
|u|_{\Eis^{1+\ve}} \cdot |\delta^{\mathrm nc}_z|_{\Eis^{-1-\ve}}
\ll_\ve
|u|_{\Eis^{1+\ve}} \cdot \Im(z)
\]
Thus, $\Eis^{1+\ve}\subset M^{1+\ve}_0\subset M$ for every $\ve>0$. \qed


Therefore, the continuous dual $M^*$ of $M$ has a natural map to the
duals $\Eis^{-1-\ve}$ of the spaces $\Eis^{1+\ve}$ for all $\ve>0$,
removing a potential ambiguity:


\begin{cor}\label{remove-ambiguity-M}
$\mathcal E \theta(s)=\theta E_{1-s}$ for $\theta=\overline{\theta}\in
M^*$. 
\end{cor}

\proof The proof of Proposition \ref{remove-ambiguity-C^0} applies
here, with the Fr\'echet space $C^0(\GH)$ replaced by the
(quasi-complete, locally convex) LF-space $M$.
\qed

Using the latter, we have

\begin{cor}\label{meromorphy-of-theta(E_s)}
The function $s\to \mathcal E \theta(s)=\theta E_s$ is a meromorphic
function of $s$.
\end{cor}

\proof From the previous corollary, indeed $\mathcal E\theta(s)=\theta
E_{1-s}$. The function $s\to E_s$ is a meromorphic $M$-valued function, so
$s\to \theta E_s$ is a meromorphic $\mathbb C$-valued function. \qed


\subsection{Meromorphic continuation of spectral synthesis integrals}
\label{mero-contn-spectral-synthesis-integrals}

Consider vector-valued integrals
\[
u^{}_{\theta,w} = \frac{\theta(1)\cdot 1}{\lambda_1-\lambda_w}
+
\frac{1}{4\pi i}\int_{(\frac 12)} \frac{\mathcal E \theta(s)\cdot
E_s}{\lambda_s-\lambda_w}\;d\ s
\]
initially defined in $\Re(w)>\frac 12$ (and $w\not=1$), where
$\lambda_s=s(1-s)$, and where $\theta=\overline{\theta}\in \Eis^{-1+\ve}$ for some
$\ve>0$. In that region, $u^{}_{\theta,w}$ solves the differential
equation $(-\Delta-\lambda_w)u=\theta$, and is in $\Eis^{1+\ve}$.
More generally, let $\Phi:M\to N$ be a continuous linear map to a
quasi-complete locally convex topological vector space $N$, and
consider the $N$-valued integrals
\[
u^{}_{\theta,w,\Phi} = \frac{\theta(1)\cdot 1}{\lambda_1-\lambda_w}
+
\frac{1}{4\pi i}\int_{(\frac 12)} \frac{\mathcal E \theta(s)\cdot
\Phi E_s}{\lambda_s-\lambda_w}\;d\ s
\]
Of course, for $\Phi$ the identity map $M\to M$, this just gives
$u^{}_{\theta,w}$ itself. We anticipate that
$\Phi(u^{}_{\theta,w})=u^{}_{\theta,w,\Phi}$, but this needs proof.


From the previous section, $\theta\in \Eis^{-1+\ve}\subset
\Eis^{-1-\ve}$ extends to an element of the dual $M^*$ to the space $M$
of moderate growth continuous functions. Thus, $\theta$ can be applied
to Eisenstein series, and, further, $\mathcal E \theta(s)=\theta
E_s$ by Corollary \ref{remove-ambiguity-M}. Thus, the previous
expression can be rewritten somewhat more concretely as 
\[
u^{}_{\theta,w,\Phi} = \frac{\theta(1)\cdot 1}{\lambda_1-\lambda_w}
+
\frac{1}{4\pi i}\int_{(\frac 12)} \frac{\theta E_{1-s}\cdot
\Phi E_s}{\lambda_s-\lambda_w}\;d\ s
\]
Unsurprisingly, we have

\begin{lem}\label{matching-in-right-half-plane}
$\Phi(u^{}_{\theta,w})=u^{}_{\theta,w,\Phi}$ in the region $\Re(w)>\frac 12$
and $w\not=1$.
\end{lem}

\proof Again, in the region $\Re(w)>\frac 12$ and $w\not=1$, the
hypotheses guarantee, via the extended Plancherel theorem, that the
integral for $u^{}_{\theta,w}$ is an $\Eis^{1+\ve}$-valued holomorphic
function of $w$. In that region, using properties of compactly
supported, continuous-integrand Gelfand-Pettis integrals from
\ref{Gelfand-Pettis}, 
\[
\begin{split}
\Phi\Big(\int_{(\frac 12)} \frac{\theta E_{1-s}\cdot
E_s}{\lambda_s-\lambda_w}\;\d s\Big)
=
\Phi\Big(\lim_{T\to +\infty} \int_{|\Im(s)|\le T} \frac{\theta E_{1-s}\cdot
E_s}{\lambda_s-\lambda_w}\;\d s\Big)
\\
=
\lim_{T\to +\infty} \Phi \int_{|\Im(s)|\le T} \frac{\theta E_{1-s}\cdot
E_s}{\lambda_s-\lambda_w}\;\d s
=
\lim_{T\to +\infty} \int_{|\Im(s)|\le T} \frac{\theta E_{1-s}\cdot
\Phi E_s}{\lambda_s-\lambda_w}\;\d s
\\
=
\int_{(\frac 12)} \frac{\theta E_{1-s}\cdot \Phi E_s}{\lambda_s-\lambda_w}\;\d s
\hskip100pt
\end{split}
\]
since the limit is approached in $\Eis^{1+\ve}\subset M$. \qed

Further specific applications to $u^{}_{\theta,w}$ will be presented
later as corollaries to the following theorem.

\begin{thm}\label{mero-contn-N-valued-integral}
With continuous linear $\Phi:M\to N$, with $N$ quasi-complete and
locally convex, the $\Phi M$-valued function $w\to u^{}_{\theta,w,\Phi}$ has a
meromorphic continuation as an $N$-valued function of $w$. In further
detail, the function 
\[
J_{\theta,w,\Phi} =
\frac{\theta(1)\cdot \Phi(1)}{(\lambda_1-\lambda_w)\cdot \langle 1,1\rangle}
+
\frac{1}{4\pi i}\int\limits_{(\frac 12)}
\frac{\theta E_{1-s}\cdot \Phi E_s-\theta E_{1-w} \cdot \Phi E_w}{\lambda_s-\lambda_w}\;\d s
\]
has a meromorphic continuation to an $N$-valued function with
the functional equation
$J_{\theta,1-w,\Phi} =J_{\theta,w,\Phi}$, and
\[
u^{}_{\theta,w,\Phi}
=
J_{\theta,w,\Phi} + \frac{\theta E_{1-w}\cdot \Phi E_w}{2(1-2w)}
\]
\end{thm}

\begin{rem} The continuation assertion for $u^{}_{\theta,w}$ itself
is stronger than, for example, the assertion of meromorphic continuation of the
numerical, pointwise integrals
\[
u^{}_{\theta,w}(z_o) = \frac{\theta(1)\cdot 1}{(\lambda_1-\lambda_w)
\cdot \langle 1,1\rangle}
+
\frac{1}{4\pi i}\int_{(\frac 12)}
\frac{\theta E_{1-s}\cdot E_s(z_o)}{\lambda_s-\lambda_w}\;\d s
\]
for fixed $z_o\in \mathfrak H$
\end{rem}


\proof From the lemma, in the region $\Re(w)>\frac 12$ and $w\not=1$,
the expression for $u^{}_{\theta,w,\Phi}$ converges as an $N$-valued
integral. The meromorphic continuation will be obtained through
rearrangement of the integral. First, in $\Re(w)>\frac 12$ and
$w\not=1$, we can certainly add and subtract to obtain
\[
\begin{split}
u^{}_{\theta,w,\Phi}=\frac{\theta(1)\cdot \Phi(1)}{(\lambda_1-\lambda_w)\cdot \langle 1,1\rangle}
+
\frac{1}{4\pi i}\int\limits_{(\frac 12)}
\frac{\theta E_{1-s}\cdot \Phi E_s}{\lambda_s-\lambda_w}\;\d s
\hskip100pt
\\
=
\frac{\theta(1)\cdot \Phi(1)}{(\lambda_1-\lambda_w)\cdot \langle
1,1\rangle}
+
\frac{1}{4\pi i}\int\limits_{(\frac 12)}
\frac{\theta E_{1-s}\cdot \Phi E_s-\theta E_{1-w}\cdot \Phi
E_w}{\lambda_s-\lambda_w}\;\d s
\hskip20pt
\\
+
\theta E_{1-w}\cdot \Phi E_w \frac{1}{4\pi i}
\int\limits_{(\frac 12)} \frac{\d s}{\lambda_s-\lambda_w}
\\
=  J_{\theta,w,\Phi} 
+
\theta E_{1-w}\cdot \Phi E_w \frac{1}{4\pi i}
\int\limits_{(\frac 12)} \frac{\d s}{\lambda_s-\lambda_w}
\hskip150pt
\end{split}
\]
By residues,
\[
\begin{split}
\theta E_{1-w}\cdot \Phi E_w \frac{1}{4\pi i}\int\limits_{(\frac 12)}
\frac{\d s}{\lambda_s-\lambda_w}
=
\theta E_{1-w}\cdot \Phi E_w\Big(-\frac 12\cdot \res_{s=w} \frac{1}{\lambda_s-\lambda_w}\Big)
\\
=
\frac{\theta E_{1-w}\cdot \Phi E_w}{2(1-2w)}
\hskip150pt
\end{split}
\]
Since $w\to E_{1-w}$ is a meromorphic $M$-valued function and
$\theta\in M^*$, the function $w\to \theta E_{1-w}$ is a meromorphic
$\mathbb C$-valued function. Similarly, $w\to \Phi E_w$ is meromorphic $N$-valued.
Thus, $\theta E_{1-w}\cdot \Phi E_w$ is a meromorphic
$N$-valued function, with a meromorphic continuation from the
meromorphic continuation of the Eisenstein series. Although the 
numerator is invariant under $w\to 1-w$ by the functional equation of
the Eisenstein series, the denominator is skew-symmetric.


Now we meromorphically continue the integral $J_{\theta,w,\Phi}$.
Constrain $w$ to lie in a fixed compact set $C$, and take
$T$ large enough so that $T\ge 2|w|$ for all $w\in C$. 
At first for $\Re(w)>\frac 12$, make the obvious attempt to cancel
vanishing of the denominator when $s$ is close to $w$, by rearranging
\[
\begin{split}
J_{\theta,w,\Phi} - \frac{\theta(1)\cdot \Phi(1)}{(\lambda_1-\lambda_w)\cdot
\langle 1,1\rangle}
=
\frac{1}{4\pi i}\int\limits_{(\frac 12)}
\frac{\theta E_{1-s}\cdot \Phi E_s-\theta E_{1-w}\cdot
\Phi E_w}{\lambda_s-\lambda_w}\;\d s
\hskip70pt
\\
=
\frac{1}{4\pi i}\int\limits_{|t|\ge T}
\frac{\theta E_{1-s}\cdot \Phi E_s}{\lambda_s-\lambda_w}\;\d s
+
\theta E_{1-w}\cdot \Phi E_w \cdot
\frac{1}{4\pi i}\int\limits_{|t|\ge T}
\frac{1}{\lambda_s-\lambda_w}\;\d s
\\
+
\frac{1}{4\pi i}\int\limits_{|t|\le T}
\frac{\theta E_{1-s}\cdot \Phi E_s-\theta E_{1-w}\cdot
\Phi E_w}{\lambda_s-\lambda_w}\;\d s
\end{split}
\]
The meromorphy of the leading integral in the case of $u^{}_{\theta,w}$
itself is understood via the Plancherel theorem on the continuous
automorphic spectrum. Ignoring constants, the Plancherel theorem for
$\Eis^0$ asserts that, for $t\to A(t)$ in $L^2(\mathbb R)$, the spectral
synthesis integral 
\[
B(z) = \frac{1}{4\pi}\int_{-\infty}^\infty A(t)\cdot E_s(z)\;\d s
\]
for $z\in \mathfrak H$ produces a function $B$ in $\Eis^0$, and the map $A\to B$ gives an
isometry. Since $L^2$ functions in $\Eis^0$
certainly need not have moderate pointwise growth at infinity, to have
a continuous inclusion to $M$ it is necessary to use $\Eis^{1+\ve}$. For
$\Eis^r$ for general index $r$, Plancherel becomes the following. Let
$X_r$ be the measurable functions $t\to A(t)$ (modulo null functions)
on $\mathbb R$ such that
$\int_\mathbb R |A(t)|^2\cdot (\frac 14+t^2)^r\;dt<\infty$. Then the spectral
synthesis integral 
produces a function $B$ in $X_r$, and $A\to B$ gives an isometry
$\Eis^r\to X_r$ (ignoring constants). Since $\theta\in \Eis^{-1+\ve}$, for
$w$ in a fixed compact, 
\[
\int_{|t|\ge T} \Big|
\frac{\theta E_{\frac 12-it}}{
\lambda_{\frac 12+it}-\lambda_w}\Big|^2
\cdot \left(\frac 14+t^2\right)^{1+\ve}\;\d t\;<\; \infty
\]
Indeed, the $V^{1+\ve_o}$-valued function
that is
$w\to(t \to \frac{\theta E_{\frac {1}{2}-it}}{\lambda_{\frac{1}{2}+it}-\lambda_w})$
for $|t|\ge T$, and is $0$ for $|t|< T$, 
is directly seen to be complex differentiable $M$-valued in $w$,
hence, holomorphic. Composition with the Plancherel isometry shows that
\[
w\;\to\;
\frac{1}{4\pi i}\int\limits_{|t|\ge T}
\frac{\theta E_{1-s}\cdot E_s}{\lambda_s-\lambda_w}\;\d s
\]
is a meromorphic $\Eis^{1+\ve}$-valued function in $w$ in the fixed
compact. Since $|w|\ll T$, the meromorphic continuation is given by the same
integral, the invariance of the integrand under $w\to 1-w$ remains.
This treats the first term for $u^{}_{\theta,w}$.

To address the first term for $u^{}_{\theta,w,\Phi}$ for continuous
$\Phi:M\to N$, with inclusion $\Eis^{1+\ve}\subset M$, use properties of
compactly supported, continuous-integrand Gelfand-Pettis integrals
from \ref{Gelfand-Pettis}: 
\[
\begin{split}
\Phi\Big(\int_{|\Im(s)|\ge T} \frac{\theta E_{1-s}\cdot
E_s}{\lambda_s-\lambda_w}\;\d s\Big)
=
\Phi\Big(\lim_{T'\to +\infty} \int_{T'\ge|\Im(s)|\ge T} \frac{\theta E_{1-s}\cdot
E_s}{\lambda_s-\lambda_w}\;\d s\Big)
\\
=
\lim_{T'\to +\infty} \Phi \int_{T'\ge|\Im(s)|\ge T} \frac{\theta E_{1-s}\cdot
E_s}{\lambda_s-\lambda_w}\;\d s
=
\lim_{T'\to +\infty} \int_{T'\ge|\Im(s)|\ge T} \frac{\theta E_{1-s}\cdot
\Phi E_s}{\lambda_s-\lambda_w}\;\d s
\\
=
\int_{|\Im(s)|\ge T} \frac{\theta E_{1-s}\cdot \Phi E_s}{\lambda_s-\lambda_w}\;\d s
\hskip100pt
\end{split}
\]
since the limit is approached in $\Eis^{1+\ve}\subset M$. Thus, the
meromorphic continuation in $\Eis^{1+\ve}$ of the first term for
$u^{}_{\theta,w}$ gives that for $u^{}_{\theta,w,\Phi}$.

In the second summand
\[
\theta E_{1-w}\cdot \Phi E_w \cdot
\frac{1}{4\pi i}\int\limits_{|t|\ge T}
\frac{1}{\lambda_s-\lambda_w}\;\d s
\]
the leading $\theta E_w \cdot \Phi E_w$ has a meromorphic continuation and
is invariant under $w\to 1-w$. Since $|w|\ll T$, the meromorphic
continuation of the integrand is given by the the same integral, and
the invariance under $w\to 1-w$ remains.

The remaining summand
\[
\frac{1}{4\pi i}\int\limits_{|t|\le T}
\frac{\theta E_{1-s}\cdot \Phi E_s-\theta E_{1-w}\cdot
\Phi E_w}{\lambda_s-\lambda_w}\;\d s
\]
is a compactly-supported vector-valued integral. To show that the
integral is a meromorphic $N$-valued function of $w$, we invoke the
Gelfand-Pettis criterion for existence of a weak integral. Let
${\mathrm {Hol}}(\Omega,N)$ be the topological vector space of holomorphic
$N$-valued functions on a fixed open $\Omega$ avoiding 
poles of $E_w$, with compact closure $C$. It suffices to show that the
integrand extends to a continuous ${\mathrm {Hol}}(\Omega,N)$-valued function of $s$, where
${\mathrm {Hol}}(\Omega,N)$ has the natural quasi-complete locally convex topology
as in section \ref{holomorphic-vector-valued}.
Unsurprisingly, to show that the integrand extends to a holomorphic
(hence, continuous)
${\mathrm {Hol}}(\Omega,N)$-valued function of
$s$, it suffices to show that the integrand extends to a holomorphic
$N$-valued function of the two complex variables $s,w$.

By Cauchy-Goursat theory for vector-valued holomorphic functions, 
near a point $s_o$, the $N$-valued function $s\to \Phi E_s$ has
a convergent power series expansion
\[
\Phi E_s = A_0 + A_1(s-s_o) + A_2(s-s_o)^2+\ldots
\]
with $A_i\in N$, and $\theta E_s$ has a scalar power series expansion
\[
\theta E_s = \theta(A_0) + \theta(A_1)\cdot (s-s_o) +
\theta(A_2)\cdot (s-s_o)^2+\ldots
\]
Thus, for suitable coefficients $B_n\in N$, 
\[
\theta E_{1-s}\cdot \Phi E_s =
B_0 + B_1(s-s_o)+B_2(s-s_o)^2+\ldots
\]
Then
\[
\begin{split}
\theta E_{1-s}\cdot E_s - \theta E_{1-w}\cdot \Phi E_w
\hskip220pt
\\
=B_1((s-s_o)-(w-s_o)) + B_2((s-s_o)^2-(w-s_o)^2) + \ldots
\hskip32pt
\\
=
((s-s_o)-(w-s_o))\cdot \hbox{(convergent power series in $s-s_o$, $w-s_o$)}
\\
=
(s-w)\cdot \hbox{(convergent power series in $s-s_o$, $w-s_o$)}
\hskip68pt
\end{split}
\]
Holomorphy is a local property. Thus, the integrand, initially defined
only for $s\not=w$, extends to a holomorphic $N$-valued function
$F(s,w)$ including the diagonal $s=w={\frac 12}+it$ with $|t|\le T$, as
well. Thus, the ${\mathrm Hol}(\Omega,N)$-valued function $f(s)$ given by
$f(s)(w)=F(s,w)$ is holomorphic in $w$.
Thus, there is a Gelfand-Pettis integral $\int_{|t|\le T}
f({\frac 12}+it)\;dt$ in $\mathrm {Hol}(\Omega, N)$, as desired. This proves the
meromorphic continuation. Further, the $w\to 1-w$ symmetry is retained
by the extension of the integrand to the diagonal. \qed

\begin{cor}
For $\Phi:M\to N$ a continuous map from the space $M$ of
moderate-growth continuous functions on $\GH$ to a quasi-complete,
locally convex topological vector space $N$, the meromorphic
continuations satisfy
\[
\Phi(u^{}_{\theta,w})=u^{}_{\theta,w,\Phi}
\]
\end{cor}

\proof In $\Re(w)>\frac 12$, lemma \ref{matching-in-right-half-plane}
proves the asserted equality. Then the vector-valued version of the
identity principle of complex analysis gives equality of the
meromorphic continuations. \qed

\begin{cor}\label{Debunk-illusion-of-symmetry}
With $\theta\in \Eis^{-1+\ve}$ for some $\ve>0$, for $\Re(w)<\frac 12$,
\[
u^{}_{\theta,w} = u^{}_{\theta,1-w}+\frac{\theta E_w\cdot E_{1-w}}{1-2w}
\]
In particular, for $\Re(w)<\frac 12$, $u^{}_{\theta,w}\in \Eis^{1+\ve}$
if and only if $\theta E_{1-w}=0$.
\end{cor}

\proof First, $\theta E_{1-w}\cdot E_w=\theta E_w\cdot E_{1-w}$, using
the functional equations $E_{1-w}=E_w/c_w$ and $c_wc_{1-w}=1$.
Then, from the theorem, with $\Re(w)<\frac 12$, 
\[
\begin{split}
u^{}_{\theta,w}
= J_{\theta,w} + \frac{\theta E_{1-w}\cdot E_w}{2(1-2w)}
= J_{\theta,1-w} - \frac{\theta E_{w}\cdot E_{1-w}}{2(1-2(1-w))}
\hskip100pt
\\
= \Big(J_{\theta,1-w} + \frac{\theta E_{w}\cdot E_{1-w}}{2(1-2(1-w))}\Big)
- 2\frac{\theta E_{w}\cdot E_{1-w}}{2(1-2(1-w))}
= u^{}_{\theta,1-w} 
- 2\frac{\theta E_w\cdot E_{1-w}}{2(1-2(1-w))}
\end{split}
\]
Since $\Re(1-w)>\frac 12$, $u^{}_{\theta,1-w}\in \Eis^{1+\ve}$, so
$u^{}_{\theta,w}\in \Eis^{1+\ve}$ if and only if the extra summand
vanishes. \qed

Again, taking $\Phi=\theta$ in the above, let
\[
J_{\theta,w} = J_{\theta,w,\theta} =
\frac{\theta(1)\cdot 1}{(\lambda_1-\lambda_w)\cdot \langle 1,1\rangle}
+
\frac{1}{4\pi i}\int\limits_{(\frac 12)}
\frac{\theta E_{1-s}\cdot \theta E_s-\theta E_{1-w} \cdot \theta
E_w}{\lambda_s-\lambda_w}\;\d s 
\]

\begin{cor}\label{on-line-imaginary-part}
For $\Re(w)=\frac 12$,
$
\Re\pt \theta u^{}_{\theta,w} = J_{\theta,w} 
$
and
$\Im\pt \theta u^{}_{\theta,w} =
-\theta E_{1-w}\cdot \theta E_w/4\Im(w)$.
\end{cor}

\proof In Theorem \ref{mero-contn-N-valued-integral}, with
$\Phi=\theta$, on $\Re(w)=\frac 12$, complex conjugation applied to
(the meromorphic continuation of the) the integral for
$J_{\theta,w,\theta}$ sends $w\to 1-w$, under which
$J_{\theta,w,\theta}$ is invariant, since $\overline
\theta=\theta$. On the other hand, by the functional equation of
$E_w$, the term $\theta E_{1-w}\cdot \theta E_w/2(1-2w)$ is purely
imaginary. 
\qed


\subsection{Spectral corollaries of meromorphic continuation}

\begin{cor} \label{mero-contn-in-V} With $\theta\in \Eis^{-1+\ve}$ for
some $\ve>0$, at points $w$ with $\Re(w)\le \frac
12$, the meromorphic continuation of $u^{}_{\theta,w}$ as $M$-valued
function of $w$ is in $\Eis^{1+\ve}$ if and only if $\theta E_w=0$.
\end{cor}

\proof The theorem \ref{mero-contn-N-valued-integral} with $\Phi$ the
identity map of $M$ to itself gives the meromorphic continuation of
$u^{}_{\theta,w}$ as an $M$-valued function. For $\theta E_w=0$, the
extra term involving Eisenstein series in the expression
\[
u^{}_{\theta,w}
=
J_{\theta,w} + \frac{\theta E_{1-w}\cdot E_w}{2(1-2w)}
\]
disappears. Further, in the numerator in the integral
\[
J_{\theta,w} =
\frac{\theta(1)\cdot 1}{\lambda_1-\lambda_w}
+
\frac{1}{4\pi i}\int\limits_{(\frac 12)}
\frac{\theta E_{1-s}\cdot E_s-\theta E_{1-w} \cdot E_w}{\lambda_s-\lambda_w}\;\d s
\]
the term $\theta E_w$ vanishes. Thus, the spectral coefficient $s\to \theta
E_{1-s}/(\lambda_s-\lambda_w)$ in the integrand is in $X_{1+\ve}$. By
the extended Plancherel theorem, the integral (extended by continuity)
is in $\Eis^{1+\ve}$. \qed

\begin{cor}\label{representability-of-eigenfunctions} With $\theta\in
\Eis^{-1+\ve}$ for some $\ve>0$, if the Friedrichs
extension $\widetilde S_\theta$ has an eigenfunction $u$ with
eigenvalue $\lambda_w>\frac 14$, then $u$ is a multiple of the
meromorphic continuation of $u^{}_{\theta,w}$ to $\Re(w)=\frac
12$. Further, $\theta E_w=0=\theta E_{1-w}$ and $\theta
u^{}_{\theta,w}=0$.  Conversely, if $\theta E_w=0=\theta E_{1-w}$ and
$\theta u^{}_{\theta,w}=0$, then $u^{}_{\theta,w}$ is an eigenfunction of
$\widetilde S_\theta$.  \end{cor}

\proof From theorem \ref{discrete-spectrum-if-any}, if
$\lambda_w>\frac 14$ is an eigenvalue for $\widetilde S_\theta$, then
$\theta E_w=0$, and $\Re(w)=\frac 12$. Then, also, $\theta E_{1-w}=0$ by
the functional equation of the Eisenstein series. Thus,
$u^{}_{\theta,w}\in \Eis^{1+\ve}$, by the previous. Then for some constant
$c$, $v=u-c\cdot u^{}_{\theta,w}$ satisfies the homogeneous equation
$(-\Delta-\lambda_w)v=0$, which has no non-zero solution in $\Eis^1$. For
$u\in \Eis^1$ to be in the domain of the Friedrichs extension, it is
necessary and sufficient that $\theta u=0$. \qed

As a special case of the theorem \ref{mero-contn-N-valued-integral}
with $\Phi=\theta:M\to \mathbb C$, whose relevance is clearer in light of the
two previous corollaries, we have

\begin{cor} The condition $\theta u^{}_{\theta,w}=0$ on the
meromorphically continued $u^{}_{\theta,w}$ is $u^{}_{\theta,w,\theta}=0$
with the meromorphic continuation of $u^{}_{\theta,w,\theta}$ as in 
theorem \ref{mero-contn-N-valued-integral}.
\end{cor}

\begin{rem}
If we did not know that any solution $u\in \Eis^{1+\ve}$ to the
differential equation were of the form $u^{}_{\theta,w}$, it would be
more difficult to appraise the vanishing requirement $\theta u=0$.
\end{rem}

\begin{cor}\label{final-condition-for-eigenfunction} A necessary and
sufficient condition for $\lambda_w>\frac 14$ to be an eigenvalue for
$\widetilde S_\theta$ is that $\theta u^{}_{\theta,w}=0$.
\end{cor}

\proof By corollary \ref{representability-of-eigenfunctions}, in that range,
the eigenfunctions for $\widetilde S_\theta$ are exactly
$u^{}_{\theta,w}$ with $\theta E_w=0$ and $\theta u^{}_{\theta,w}=0$. As
above, 
\[
\begin{split}
\theta u^{}_{\theta,\theta}
=
\frac{\theta(1)\cdot \theta(1)}{(\lambda_1-\lambda_w)\cdot \langle 1,1\rangle}
+
\frac{1}{4\pi i}\int\limits_{(\frac 12)}
\frac{\theta E_{1-s}\cdot \theta E_s-\theta E_{1-w} \cdot \theta
E_w}{\lambda_s-\lambda_w}\;\d s
\\
+
\frac{\theta E_{1-w}\cdot \theta E_w}{2(1-2w)}
\end{split}
\]
On $\Re(w)=\frac 12$, the first two summands are real, while the third
is purely imaginary. Thus, for this to vanish on $\Re(w)=\frac 12$
requires vanishing of $\theta E_w$. That is, on $\Re(w)=\frac 12$,
vanishing of $\theta u^{}_{\theta,w}$ implies that of $\theta E_w$. \qed

For example, the self-adjoint operator $\widetilde S_\Theta$ on
non-cuspidal pseudo-cuspforms, has purely discrete spectrum, by
Theorem \ref{discretization}. Theorem \ref{exotic-eigenfunctions}
already noted that, for $\Re(w)=\frac 12$, $\lambda_w=w(1-w)$
is an eigenvalue if and only if $a^w+c_wa^{1-w}=0$. Indeed,
$\lambda_w=w(1-w)$ is real if and only if $\Re(w)=\frac 12$ or $w\in
[0,1]$. However, this conclusion does not address the possible
vanishing of $a^w+c_wa^{1-w}$ {\itshape off} $\Re(w)=\frac 12$ and
$[0,1]$. As an instance of desirable corollaries of spectral theory,
we have a result already observed in \cite{Hejhal-Polya}, but there
for non-spectral reasons:

\begin{cor}\label{constant-term-vanishing} $a^w+c_wa^{1-w}=0$ implies
$\Re(w)=\frac 12$ or $w\in [0,1]$. 
\end{cor}

\proof Let $u^{}_{\eta_a,w}$ be the meromorphic continuation of the
spectral expansion
\[
u^{}_{\eta_a,w} = \frac{\eta_a(1)\cdot 1}{(\lambda_1-\lambda_w)\cdot
\langle 1,1\rangle}
+\frac{1}{4\pi i}\int_{(\frac 12)}
\frac{\eta_a E_{1-s}\cdot E_s}{\lambda_s-\lambda_w} \;\d s
\]
By Corollary \ref{representability-of-eigenfunctions}, $u^{}_{\eta_a,w}$ gives an
eigenfunction of $\widetilde S_\Theta$ if and only if $\eta_aE_w=0$
and $\eta_a u^{}_{\eta_a,w}=0$. Of course,
$\eta_aE_w=a^w+c_wa^{1-w}$. The computation of Theorem \ref{eta-v}
gives
\[
\eta_a u^{}_{\eta_a,w} = \frac{a^{1-w}}{1-2w} \cdot (a^w+c_wa^{1-w})
\]
The leading factor $a^{1-w}/(1-2w)$ does not vanish, and does not have
poles except at $w=\frac 12$. Thus, if $a^w+c_wa^{1-w}=0$, then
$u^{}_{\eta_a,w}$ is an eigenfunction for $\widetilde S_\Theta$, and
necessarily $\lambda_w$ is real and non-negative. Thus, for $w$
outside the set where $\Re(w)=\frac 12$ or $w\in [0,1]$, the
expression $a^w+c_wa^{1-w}$ cannot vanish.
\qed


\subsection{Remarks on appealing but incorrect arguments}

The invariance of $\lambda_w$ under $w\longleftrightarrow 1-w$ might
suggest invariance of $u^{}_{\theta,w}$ under $w\longleftrightarrow 1-w$.
However, this is illusory, by corollary
\ref{Debunk-illusion-of-symmetry}. Nevertheless, the implications of  
(the generally incorrect claim) $u^{}_{\theta,1-w}=u^{}_{\theta,w}$ are
striking.
That is, if we were to believe that $u^{}_{\theta,1-w}=u^{}_{\theta,w}$, the
the trivial non-vanishing of remark \ref{trivial-non-vanishing} would
(fairly provocatively, but incorrectly, seem to) show that $\theta
u^{}_{\theta,w}\not=0$ off $\Re(w)=\frac 12$ and $[0,1]$.

Corollary \ref{final-condition-for-eigenfunction} (correctly) shows that for
$\lambda_w>\frac 14$, that is, for $\Re(w)=\frac 12$ and $w\not=\frac
12$, a necessary and sufficient condition for $\lambda_w$ to be an
eigenvalue for the self-adjoint operator $\widetilde S_\theta$ is that
$\theta u^{}_{\theta,w}=0$.
An argument principle discussion (going back to \cite{Backlund}, see
\cite{T} pp. 212-213) shows that the asymptotic count of the zeros $w$
of $\theta u^{}_{\theta,w}$ with imaginary parts $0\le \Im(w)\le T$ is
$\frac{T}{\pi} \log T + O(T)$.

Thus, if we (erroneously) believe the confinement of zeros of $\theta
u^{}_{\theta,w}$ to $\Re(w)=\frac 12$ (and $[0,1]$), we (seem to) find that $\theta
u^{}_{\theta,w}$ asymptotically has $\frac{T}{\pi} \log T + O(T)$ zeros
from height $0$ to $T$ on the critical line. Since $\theta
u^{}_{\theta,w}=0$ on the critical line does truly entail $\theta E_w=0$,
we (seem to) conclude that the asymptotic count of zeros of $\theta
E_w$ on is $\frac{T}{\pi} \log T + O(T)$. For $\theta$ the Eisenstein-Heegner
distribution attached to a fundamental discriminant $-d<0$ and
$k=\mathbb Q(\sqrt{-d})$, this would (seem to) imply that $\zeta_k(s)$
has $100\%$ of its zeros on the critical line, the
asymptotic version of the Riemann Hypothesis for $\zeta_k(s)$.
%
%
%
The argument is fundamentally
flawed, at the point where one fallaciously imagines that
$u^{}_{\theta,1-w}=u^{}_{\theta,w}$: the functional equation of
$u^{}_{\theta,w}$ involves an extra term, which vanishes exactly for
$\theta E_w=0$, thus thwarting this naively optimistic approach.


\section{Spacing of spectral parameters}


\subsection{Exotic eigenfunction expansions of distributions}

Fix $a>1$. Let $\{f_j:j=1,2,\ldots\}$ be eigenfunctions for the
operator $\widetilde S_\Theta$ of section 4.2, with eigenvalues
$\lambda_{s_1}\le\lambda_{s_2}\le\ldots$, with $|f_j|_V=1$. As in
Corollary \ref{compare-norms}, $|f_j|_{\Eis^1}=\sqrt{\lambda_{s_j}}$. 
Let $i_\Theta:\Eis^1\cap \Eis^0_\Theta \to \Eis^1$ be the inclusion,
and $i_\Theta^*:\Eis^{-1}\to (\Eis^1\cap \Eis^0_\Theta )^*$ its
adjoint. By Theorem 
\ref{density-by-smooth-cutoffs} , $\Eis^1\cap \Eis^0_\Theta $ is the
$\Eis^1$-topology completion of $D\cap \Eis^0_\Theta $, where $D$ is the space
of pseudo-Eisenstein series with test-function data. Thus, $\{f_j\}$
is an orthonormal basis for $\Eis^0_\Theta $, and
$\{f_j/\sqrt{\lambda_{s_j}}\}$ is an orthonormal basis for $\Eis^1\cap
\Eis^0_\Theta $.

At first for finite sums, define Sobolev-like norms by
\[
\Big|\sum_j c_j\cdot f_j\Big|_{W_r}^2 = \sum_j |c_j|^2\cdot \lambda_{s_j}^r
\]
and let $W_r$ be the completion of the space of finite linear combinations
of the vectors $f_j$ with respect to this norm. Thus, $\{f_j\cdot
\lambda_{s_j}^{-r/2}\}$ is an orthonormal basis for $W_r$.

\begin{thm} The map $W_{-1}\approx (\Eis^1\cap \Eis^0_\Theta )^*$ given by (at
first for finite sums, and then extending by continuity)
\[
\Big(\sum_i a_i\cdot f_i\Big)\Big(\sum_j b_j\cdot f_j\Big)
= \sum_j a_i\cdot \overline{b_j}
\]
for $\sum_j b_j\cdot f_j\in \Eis^1\cap \Eis^0_\Theta $ is an isomorphism.
\end{thm}

\proof First, since $f_j\in \Eis^1$ by the Friedrichs construction, and
$f_j\in \Eis^0_\Theta $, it makes sense to apply such functionals in
$\Eis^{-1}$ to the eigenfunctions $f_j$.
By Plancherel and Corollary \ref{compare-norms}, $W_0=\Eis^0_\Theta $ and
$W_1=\Eis^1\cap \Eis^0_\Theta $. Generally, the natural pairing
\[
\langle \sum_j a_j\cdot f_j,\;\sum_j b_j\cdot f_j\rangle
= \sum_j a_j\,\overline{b_j}
\]
on $W_r\times W_{-r}$ puts these two spaces in duality. Thus, we have
a natural isomorphism $W_{-1}\to W_1^*=(\Eis^1\cap \Eis^0_\Theta )^*$
respecting the duality pairings.   
\qed

Thus, functionals $\mu\in (\Eis^1\cap \Eis^0_\Theta )^*$ admit expansions 
\[
\mu = \sum_j \mu(\,\overline{f}_j)\cdot f_j
\]
convergent in the $W_{-1}$ topology, and compatible with evaluation on
the corresponding expansions of elements in $W_1=\Eis^1\cap \Eis^0_\Theta $.   


\subsection{Exotic eigenfunction expressions for solutions of
equations}\label{exotic-expressions}

Let $\{f_j:j=1,2,\ldots\}$ be eigenfunctions for the operator
$\widetilde S_\Theta$ of section 4.2, with eigenvalues
$\lambda_{s_1}\le\lambda_{s_2}\le\ldots$, with $|f_j|_V=1$.

\begin{prop}\label{discrete-solve-by-division} The operator $S_\Theta^\#$ of Lemma
\ref{S-Theta-bd-variant} is a topological isomorphism
from $W_1=\Eis^1\cap \Eis^0_\Theta $ to $W_{-1}=(\Eis^1\cap \Eis^0_\Theta )^*$,
expressible as
\[
S_\Theta^\# \sum_j b_j\cdot f_j
=
\sum_j b_j\cdot S_\Theta^\#  f_j
=
\sum_j \lambda_{s_j}\cdot b_j\cdot f_j
\]
\end{prop}

\proof Since constants are excluded from $\Eis^0_\Theta $, $S_\Theta$ is
{\itshape strictly} positive, so we are in the situation of section
2.2, and its Friedrichs extension gives the asserted isomorphism, as observed
just prior to Proposition \ref{domains}. Because $S_\Theta^\#$ is
continuous $W_1\to W_{-1}$, it commutes with limits, giving the
formula. 
\qed

\begin{thm} For $\mu\in W_{-1}=(\Eis^1\cap \Eis^0_\Theta )^*$, for
$\lambda_w\not\in [0,+\infty)$, the equation $(\widetilde
S_\Theta-\lambda_w)v=\mu$ has a unique solution in $W_1=\Eis^1\cap
\Eis^0_\Theta $,  
given by
\[
v = \sum_j \frac{\mu(\overline{f}_j)\cdot f_j}{\lambda_{s_j}-\lambda_w}
\]
This expansion is holomorphic $W_1$-valued off $\Re(w)=\frac 12$, and
has a $W_1$-valued analytic continuation to $\mathbb C$ with the
discrete set of spectral parameters $s_j$ removed.
\end{thm}

\proof First, by the construction of the Friedrichs extension, all
eigenfunctions $f_j$ are in $\Eis^1$, so any $\theta\in \Eis^{-1}$ can be
sensibly applied to them via the $\Eis^1\times \Eis^{-1}$ pairing. The bound 
\[
\sum_j |\mu(\overline{f}_k)|^2\cdot \lambda_{s_j}^{-1} < \infty
\]
implies
\[
\sum_j
\Big|\frac{\mu(\overline{f}_j)}{\lambda_{s_j}-\lambda_w}\Big|^2\cdot
\lambda_{s_j} < \infty 
\]
so the given expansion for $v$ is indeed in $W_1$. Since
$S_\Theta^\#:W_1\to W_{-1}$ is continuous, it commutes with the
implied limits in the infinite sums. 

To show the holomorphy, from section 6.2 it suffices to prove that for
$w$ in a fixed compact $C$ not meeting $\mathbb R$, for every $\ve>0$
there is $i_o$ such that for all $i_1,i_2\ge i_o$ 
\[
\sup_{w\in C}
\Big|\sum_{i_1\le j\le i_2}
\frac{\mu(\overline{f}_j)}{\lambda_{s_j}-\lambda_w}
\Big|_{W_1}
< \ve
\]
This follows from $|\lambda_{s_j}-\lambda_w|\gg_C |\lambda_{s_j}|$.
By Theorem \ref{discretization}, $\widetilde S_\Theta$ has compact
resolvent, so the parameters $s_j$ have no accumulation point in
$\mathbb C$. Thus, the same inequality holds for $C$ compact not
meeting the set of spectral parameters $s_j$, giving the analytic
continuation. 
\qed


Let $\theta=\overline{\theta}=\theta_{\mathcal D} = \sum_d
\nu_d\,\theta_d$ be a finite 
real-linear combination of Heegner distributions $\theta_d$, as in
3.6, with fundamental discriminants $d<0$. As noted in Section
\ref{Heegner-distributions}, $\theta\in \Eis^{-1+\ve}$ for some
$\ve>0$. Let $m(\mathcal D)=\max_d |d|^{\frac 12}/2$ be the maximum 
over $d$ appearing with non-zero coefficient, and increase $a$ if
necessary so that $a>1$.

\begin{cor}\label{discretized-solution} The equation
$(S_\Theta^\#-\lambda_w)v=i^*_\Theta\theta$ has unique solution  
\[
v_{\theta,w} =
\sum_j \frac{\theta(\overline{f}_j)}{\lambda_{s_j}-\lambda_s}\cdot f_j
\]
in $W_1=\Eis^1\cap \Eis^0_\Theta $, a meromorphic $W_1$-valued
function of $w$ away from the spectral parameters $s_j$, and
\[
\theta v_{\theta,w} =
\sum_j \frac{|\theta f|^2}{\lambda_{s_j}-\lambda_s}\cdot f_j
\]
is uniformly absolutely convergent on compacts away from the spectral
parameters $s_j$, producing a holomorphic $\mathbb C$-valued function
there. 
\end{cor}

\proof The image $i^*_\Theta\theta$ is in $W_{-1}$. \qed

\begin{cor}\label{final-exotic-expression} A solution $u\in \Eis^1$ to
the equation $(-\Delta-\lambda_w)u=\theta$ is in fact in $\Eis^1\cap
\Eis^0_\Theta $, so is $v_{\theta,w}$ of the previous corollary, so has an
expansion    
\[
u = \sum_j \frac{\theta (\,\overline{f}_j)}{\lambda_{s_j}-\lambda_w}\cdot f_j
\]
convergent in $\Eis^1\cap \Eis^0_\Theta $, holomorphic in $w$ away from
spectral parameters $s_j$, and there
\[
\theta u = \sum_j \frac{|\theta(f_j)|^2}{\lambda_{s_j}-\lambda_w}
\]
\end{cor}

\proof First, we recall why such a solution $u$ must be of the form
$u^{}_{\theta,w}$. Certainly in $\Re(w)>\frac 12$ (and $w\not=1$) that
equation has the solution given by a spectral expansion converging in $\Eis^1$:
\[
u^{}_{\theta,w}
=
\frac{\theta(1)\cdot 1}{(\lambda_1-\lambda_w)\cdot
\langle 1,1\rangle}
+
\frac{1}{4\pi i}\int_{(\frac 12)} \frac{\theta E_{1-s}\cdot
E_s}{\lambda_s-\lambda_w} \;\d s 
\]
On $\Re(s)=\frac 12$, Theorem \ref{solvability-necessary-condition}
shows that the equation is solvable only if $\mathcal E\theta (w)=0$,
and Corollary \ref{remove-ambiguity-M} gives $\theta E_{1-w}=\mathcal
E\theta(w)=0$. Thus, by Theorem \ref{mero-contn-in-V}, at such $w$ the meromorphic
continuation of $u^{}_{\theta,w}$ is in $\Eis^1$. Since the homogeneous
form of the equation has no solution in $\Eis^1$, it must be that
$u=u^{}_{\theta,w}$. From Theorem \ref{eta-u}, the condition $\theta
E_w=0$ gives $u^{}_{\theta,w}\in \Eis^1\cap \Eis^0_\Theta $ for $a\ge m(\mathcal
D)$ (and $a>1$). The previous corollary applies to the image
$\mu=i_\Theta^*\theta$, noting that for $f\in \Eis^1\cap \Eis^0_\Theta
$,  
\[
(i_\Theta^*\theta)(f) = \theta(i_\Theta f) = \theta(f)
\]
suppressing explicit reference to the inclusion $i_\Theta$. Then
$u^{}_{\theta,w}$ must be the solution $v_{\theta,w}$ of the previous
corollary.
\qed


\subsection{Interleaving}\label{interleaving}

Continue to let $\theta=\overline{\theta}=\theta_{\mathcal D} = \sum_d
\nu_d\,\theta_d$ be a finite real-linear combination of Eisenstein-Heegner
distributions $\theta_d$, as in 3.6, so that all the previous results
apply. Let
\[
v_{\theta,w} =
\sum_j \frac{\theta(\overline{f}_j)}{\lambda_{s_j}-\lambda_s}\cdot f_j
\]
be as in Corollary \ref{discretized-solution}.

\begin{cor}\label{strict-increase} Let $s_j=\frac12+it_j$
and $s_{j+1}=\frac 12+it_{j+1}$ be two adjacent zeros of
$a^s+c_sa^{1-s}$ on the line $\Re(s)=\frac 12$, with
$0<t_j<t_{j+1}$. For $w=\frac12+i\tau$ with $t_j<\tau<t_{j+1}$, the
function $\tau\to\theta v_{\theta,w}$ is continuous, real-valued,
strictly increasing, goes from $-\infty$ as $\tau\to t_j$ to $+\infty$
as $\tau\to t_{j+1}$. 
\end{cor}

\proof For $\tau\in\mathbb R$, 
\[
\theta v_{\theta,\frac 12+i\tau}
=
\sum_j \frac{|\theta(f_j)|^2}{t_j^2-\tau^2}
\]
is real-valued (away from the $t_j$). Being the restriction of a
meromorphic function, $\tau\to \theta v_{\theta,\frac 12+i\tau}$ is
continuous away from the $t_j$. On one hand, for $\tau\in
(t_j,t_{j+1})$, as $\tau\to t_j^+$, the summand
$|\theta(f_j)|^2/(t_j^2-\tau^2)$ goes to $-\infty$ and the rest of the
sum remains finite. On the other hand, in that interval, as $\tau\to
t_{j+1}^-$, the summand $|\theta(f_{j+1})|^2/(t_{j+1}^2-\tau^2)$ goes
to $+\infty$ and the rest of the the sum remains finite.
Away from the poles, the derivative is 
\[
\frac{\partial}{\partial\tau}
\theta v_{\theta,\frac 12+i\tau}
= \frac{\partial}{\partial\tau}
\sum_j \frac{|\theta(f_j)|^2}{t_j^2-\tau^2}
=
2\tau
\cdot \sum_j \frac{|\theta(f_j)|^2}{(t_j^2-\tau^2)^2}
> 0
\]
so the function is strictly increasing.
\qed


The first corollary directly addresses on-the-line zeros of $\theta
E_w$:

\begin{cor}\label{nascent-spacing} Let $s_j=\frac12+it_j$
and $s_{j+1}=\frac 12+it_{j+1}$ be adjacent zeros of
$a^s+c_sa^{1-s}$ on the line $\Re(s)=\frac 12$, with
$0<t_j<t_{j+1}$. Let $w_1=\frac12+i\tau_1$ and $w_2=\frac12+i\tau_2$  
with $t_j<\tau_1<\tau_2<t_{j+1}$ with $\theta E_{w_j}=0$. Then $\theta
u^{}_{\theta,w_1}<\theta u^{}_{\theta,w_2}$.
\end{cor}

\proof As above, the condition $\theta E_w=0$ implies that the
meromorphic continuation of $u^{}_{\theta,w}$ is in $\Eis^1$, and also that
$u^{}_{\theta,w}\in \Eis^0_\Theta $. By Corollary \ref{discretized-solution},
$u^{}_{\theta,w}$ must be $v_{\theta,w}$, and $\theta v_{\theta,w}$ is
monotone in intervals $(t_j,t_{j+1})$.
\qed


This second corollary is a variant that more directly addresses
solutions of the equation $(-\Delta-\lambda_w)u=\theta$:

\begin{cor} Let $s_j$ and $s_{j+1}$ be adjacent zeros of
$a^s+c_sa^{1-s}$ on the line $\Re(s)=\frac 12$.  On the line
$\Re(w)=\frac 12$, between $s_j$ and $s_{j+1}$, there is at most one
$w$ such that a solution $u\in \Eis^1$ of the equation
$(-\Delta-\lambda_w)u=\theta$ satisfies $\theta u=0$.
\end{cor}

\proof Again, if a solution $u$ to $(-\Delta-\lambda_w)u=\theta$ is in
$\Eis^1$, then $\theta E_w=0$, by Theorem
\ref{solvability-necessary-condition}. And, again, the
meromorphically-continued $u^{}_{\theta,w}$ is in $\Eis^1$ at that point, and is the
unique solution to the equation, since the homogeneous equation
$(-\Delta-\lambda_w)u=0$ has no non-zero solution. Also, again, $\theta
E_w=0$ implies that $u^{}_{\theta,w}\in \Eis^0_\Theta $. Thus, by Corollary
\ref{final-exotic-expression}, $u^{}_{\theta,w}=v_{\theta,v}$, and then
\[
\theta u^{}_{\theta,w} =
\theta v_{\theta,w} = \sum_j \frac{|\theta(f_j)|^2}{\lambda_{s_j}-\lambda_w}
\]
The monotonicity assertion of Corollary \ref{strict-increase} shows
that there is exactly one $w=\frac 12+i\tau$ in the given interval
such that $\theta v_{\theta,w}=0$.
\qed

\begin{cor} There is at most one $w$ on $\Re(w)=\frac 12$ between
adjacent zeros $s_j$ and $s_{j+1}$ of $a^s+c_sa^{1-s}$ on
$\Re(s)=\frac 12$ such that $\lambda_w$ is an eigenvalue for
$\widetilde S_\theta$. 
\end{cor}

\proof Combine the previous corollary with Corollary
\ref{final-condition-for-eigenfunction}.
\qed


\subsection{Spacing of spectral parameters}

Continue to let $\theta=\overline{\theta}=\theta_{\mathcal D} = \sum_d
\nu_d\,\theta_d$ be a finite real-linear combination of Eisenstein-Heegner
distributions $\theta_d$, as in 3.6, so that the immediately previous
results apply.

From the previous section, the positions of the zeros of
$a^s+c_sa^{1-s}$ influence the positions of the parameters $w$
appearing in differential equations $(-\Delta-\lambda_w)u=\theta$ with
condition $\theta u=c\in\mathbb R$. The zeros $s$ of
$a^s+c_sa^{1-s}=0$ allow a degree of adjustment by choice of the
cut-off height $a>1$.

To anticipate the point of the present discussion, for $\log\log
\Im(s)$ large, the behavior of $\zeta(s)$ on the edge of the critical
strip is relatively regular, by \cite{T} (5.17.4) page 112 (in an
earlier edition, page 98). This gives an eventual regularity of
spacing of the zeros of $a^s+c_sa^{1-s}$, with implications for the
exotic eigenfunction expansions above.

As originally in \cite{Backlund}, or from \cite{T} pages
212--213, by the argument principle and Jensen's inequality,
\[
\#\Big(\hbox{zeros of $a^s+c_sa^{1-s}$ with $0\le \Im(s)\le T$}\Big)
\;\;=\;\; \frac{T}{\pi}\log T + O(T)
\]
Thus, near height $T$, the average gap is $\pi/\log T$.
As recalled in Corollary
\ref{constant-term-vanishing}, all the zeros are on $\frac 12+i\mathbb
R$ or $[0,1]$, and are in bijection with the discrete spectrum of
$\widetilde S_\Theta$ (modulo $w\leftrightarrow 1-s$) by $s\to
\lambda_s=s(1-s)$.

First, $|c_s|=1$ on $\Re(s)=\frac 12$, so $c_{\frac 12+it}=e^{-2i\psi(t)}$
with real-valued $\psi(t)=\arg \xi(1+2it)$. Since $\zeta(s)$
does not vanish on $\Re(s)=1$, $\psi(t)$ is differentiable. Letting
$s=\frac 12=it$, rearrange 
\[
\begin{split}
a^s+c_sa^{1-s}
=
\sqrt{a}\cdot e^{-i\psi( t)}\cdot \Big(e^{i t\log
a+i\psi( t)}+e^{-i t\log a-i\psi( t)}\Big)
\\
=
2\sqrt{a}\cdot e^{-i\psi( t)}\cdot
\cos\Big( t\log a+\psi( t)\Big)
\hskip80pt
\end{split}
\]
Thus, the on-the-line vanishing condition is $\cos( t\log a+\psi( t))=0$.

\begin{prop}\label{asymptotics-psi} $\psi( t)=\arg\xi(1+2i t)$ satisfies
\[
\psi( t) =  t\log t+O(\frac{ t\log  t}{\log\log  t})
\hskip30pt\hbox{and}\hskip30pt
\psi'( t) = \log  t + O(\frac{\log  t}{\log\log t})
\]
\end{prop}

\proof Of course,
\[
\psi( t)
=
\arg \xi(1+2i t)
=
- t\log \pi + \arg \Gamma(\frac 12+i t) + \arg\zeta(1+2i t)
\]
From $\log \Gamma(s)=(s-\frac 12)\log s - s +\frac 12\log 2\pi + O(s^{-1})$
\[
\begin{split}
\log \Gamma(\frac 12+i t)
=
i t\log(1+i t) - (\frac 12+i t) +\frac 12\log 2\pi + O\Big({1}{\frac 12+i t}\Big)
\\
=
i t\Big(\log t+O(\frac{1}{t^2})+i(\frac{\pi}{2}+O(\frac 1t))\Big)
 - (\frac 12+i t) +\frac 12\log 2\pi + O\Big(\frac{1}{\frac 12+i t}\Big)
\end{split}
\]
Thus,
\[
\Im \log \Gamma(\frac 12+i t)
\;=\;  t\log  t -  t + O(\frac 1t)
\]
and
\[
\psi( t)
\;=\;
(- t \log \pi) + ( t\log  t -  t) + \arg\zeta(1+2i t) +
O( t^{-1})
\]
Similarly, the asymptotic $\Gamma'(s)/\Gamma(s)=\log s+O(s^{-1})$ gives
\[
\Im \frac{d}{ds}\log \Gamma(\frac 12+i t) \;=\; \log  t + O( t^{-1})
\]
and
\[
\psi'( t)
\;=\;
(- \log \pi) + (\log  t) + \frac{d}{d t}\arg\zeta(1+2i t) + O( t^{-1})
\]
From \cite{T} (5.17.4) page 112 (in an earlier edition, page  
98), for $u\ge t$,
\[
\log \zeta(1+iu) - \log \zeta(1+it) \;=\; O\Big(\frac{\log t}{\log\log
t}\Big)\cdot (u-t)
\]
This gives the asymptotics for $\psi$ and $\psi'$. \qed


\begin{cor}\label{gap-asymptotic} Fix $a_o>1$. Given $\ve>0$ (with
$\ve<1$ for definiteness), for sufficiently large $T_o>0$,
for all $1<a\le a_o$, for consecutive real zeros
$t_j<t_{j+1}$ of $a^{\frac12+it}+c_{\frac12+it}a^{\frac12-it}$
with $T_o\le t_j<t_{j+1}$, 
\[
(1-\ve)\cdot \frac{\pi}{\log t_j}
\le
t_{j+1}-t_j
\le
(1+\ve)\frac{\pi}{\log t_j}
\]
\end{cor}

\proof This is elementary from the previous. \qed


\begin{thm}\label{derivative-of-spectral-parameters}
The real zeros $t=t(a)$ of $a^{\frac 12+it}+c_{\frac 12+it}a^{\frac
12-it}$, which are also the real zeros of $\cos(t\log a+\psi(t))$, are
differentiable functions of cut-off height $a$, with
\[
\frac{\partial t}{\partial a}
=
\frac{-t/a}{\log a + \psi'(t)}
\]
with non-vanishing denominator.
\end{thm}

\proof Implicit differentiation of $\cos(t\log a+\psi(t))=0$ gives the
formula.
Non-vanishing of the denominator at points where $\cos(t\log
a+\psi(t)=0$ follows from the Maa\ss-Selberg relation, as
follows. On one hand, for $\Im(s)>0$ the higher Fourier terms of $E_s$
do not vanish identically, so $\langle
\wedge^aE_s,\wedge^aE_s\rangle>0$.  On the other hand, the
Maa\ss-Selberg relation is 
\[
\begin{split}
\langle \wedge^a E_s,\wedge^a E_r\rangle
=
\frac{a^{s+{\overline r}-1}}{s+{\overline r}-1}
+ c_s\frac{a^{(1-s)+{\overline r}-1}}{(1-s)+{\overline r}-1}
\hskip50pt
\\
+ c_{\overline r}\frac{a^{s+(1-{\overline r})-1}}{s+(1-{\overline
r})-1}
+ c_s c_{\overline r} \frac{a^{(1-s)+(1-{\overline r})-1}}{
(1-s)+(1-{\overline r})-1}
\end{split}
\]
For $s=\frac 12+it+\ve$ and $r=\frac 12+it$, this is
\[
\langle \wedge^a E_s,\wedge^a E_r\rangle \;=\;
\frac{a^\ve}{\ve}
+ c_s\frac{a^{-2it-\ve}}{-2it-\ve}
+ c_{\overline r}\frac{a^{2it+\ve}}{2it+\ve}
+ c_s c_{\overline r} \frac{a^{-\ve}}{-\ve}
\]
As $\ve\to 0$, the two middle terms go to
\[
c_{\frac 12+it} \frac{a^{-2it}}{-2it} + c_{\frac 12-it} \frac{a^{2it}}{2it}
=
e^{-2i\psi(t)}\frac{a^{-2it}}{-2it} + e^{2\psi(t)}
\frac{a^{2it}}{2it}
=
\frac 1t\cdot \sin(2\log a+2\psi(t))
\]
The vanishing condition is $\cos(t\log a+\psi(t))=0$, so $t\log
a+\psi(t)\in \pi\mathbb Z$, so $2t\log a+2\psi(t)\in 2\pi\mathbb Z$,
and the sine function vanishes there. That is, the sum of these two
terms is $0$.

Modulo $O(\ve^2)$, the sum of the first and last terms is
\[
\begin{split}
\frac{1}{\ve}
(a^\ve -  c_{\frac 12+it+\ve} c_{\frac 12-it} a^{-\ve})
\;=\;
\frac{1}{\ve}\Big(
(1+\ve\log a) - (c_{\frac 12+it}+\ve c'_{\frac 12+it}) c_{\frac 12-it}(1-\ve\log a)
\Big)
\\
=
\frac{1}{\ve}\Big(
1+\ve\log a - (1+\ve c'_{\frac 12+it}c_{\frac 12-it})(1-\ve\log a)
\Big)
\;=\;
2\log a - \frac{c'_{\frac 12+it}}{c_{\frac 12+it}}
\end{split}
\]
since $c_sc_{1-s}=1$. Since $c_{\frac 12+it}=e^{-2i\psi(t)}$ and $c_s$ is
holomorphic, 
\[
c'_{\frac 12+it} \;=\; \frac{d}{d(it)}c_{\frac 12+it}
\;=\; -i \frac{d}{dt}c_{\frac 12+it}
\;=\; -i \frac{d}{dt}\,e^{-2i\psi(t)}
\;=\; -2\psi'(t)\cdot e^{-2i\psi(t)}
\]
Thus, $c'_{\frac 12+it}/c_{\frac 12+it}=-2\psi'(t)$. That is, $\ve\to
0$, these terms go to $2\log a+2\psi'(t)$. That is, when
$a^s+c_sa_s=0$, 
\[
0<\langle \wedge^a E_s,\wedge^aE_s\rangle=2\log a+2\psi'(t)
\]
giving the non-vanishing of the denominator. \qed

Thus, to change a given zero $t$ for cut-off height $a$ by the
expected gap amount $\pi/\log t$ at that height, change $a$ by
roughly   
\[
\frac{\hbox{gap}}{\hbox{derivative}}
\;\sim\;
\frac{
\pi/(\log a + \log t)
}{
(-t/a)/(\log a + \log t)
}
\;\sim\;
\frac{-\pi a}{t}
\]


\begin{cor} Fix $a_o>1$. Given 
$\ve>0$ (and $\ve<1$), there is $T_o$ sufficiently large such that,
for every $a$ in the range $1<a\le a_o$, for every real zero $t=t(a)$
of $a^{\frac 12+it}+c_{\frac 12+it}a^{\frac 12-it}$,
\[
(1-\ve)\cdot \frac{t/a}{\log t}
\le
-\,\frac{\partial t}{\partial a}
\le
(1+\ve) \frac{t/a}{\log t}
\]
\end{cor}

\proof Elementary from the above. \qed


Combining the previous with Corollary \ref{nascent-spacing}:

\begin{thm} Given $\ve>0$, there is $T_o>0$ sufficiently
large such that, for two zeros $w_1=\frac 12+i\tau_1$ and $w_2=\frac
12+i\tau_2$ of $\theta E_w$ with $T_o<\tau_1<\tau_2$, if 
$\theta u^{}_{\theta,w_1}\ge\theta u^{}_{\theta,w_2}$, then we have the
lower bound $|\tau_2-\tau_1|\ge (1-\ve)\pi/\log \tau_1$.
\end{thm}

\proof Elementary from the above. \qed

As a very special case, relevant to the discrete spectrum (if any) of
$\widetilde S_\theta$:

\begin{cor} Given $\ve>0$, there is $T>0$ sufficiently large such
that, for two zeros $w_1=\frac 12+i\tau_1$ and $w_2=\frac 12+i\tau_2$
of $\theta E_w$ with $T<\tau_1<\tau_2$, if
$\theta u^{}_{\theta,w_1}=0=\theta u^{}_{\theta,w_2}$, then we have the
lower bound $|\tau_2-\tau_1|\ge (1-\ve)\pi/\log T$. That is, the
discrete spectrum, if any, of $\widetilde S_\theta$, is $\lambda_{w_j}$
with a lower bound $|w_{j+1}-w_j|\ge (1-\ve)\pi/\log T$ for
adjacent $w_j$ and $w_{j+1}$ at height $T$. 
\end{cor}

\proof This is a special case of the previous theorem and Corollary
\ref{final-condition-for-eigenfunction}.  \qed


\subsection{Juxtaposition with pair-correlation}

The asymptotic lower bound $\pi/\log T$ on spacing of consecutive
spectral parameters $w_j$ at height $T$ for eigenvalues
$\lambda_{w_j}$ of $\widetilde S_\theta$ is in conflict with
Montgomery's pair-correlation \cite{Montgomery-pair}. We continue to
take $\theta$ to be a real-linear combination of Eisenstein-Heegner
distributions, so that all previous results apply.

For example, assuming the Riemann Hypothesis and {\it pair
correlation}:

\begin{cor}\label{ninety-four} At most 94\% of the zeros of $\zeta(s)$
  give eigenvalues 
$\lambda_s=s(1-s)$ for $\widetilde S_\theta$.
\end{cor}

\begin{rem}
This shows that the optimistic simple version of the conjecture at the
end of \cite{CdV-pseudo} cannot hold. Namely, with $\theta$ the
Eisenstain-Dirac $\delta$ at $e^{2\pi i/3}$, it
{\itshape cannot} be the case that all zeros $\rho_j$ of $\zeta(s)$
give eigenvalues for $\widetilde S_\theta$.
\end{rem}

\begin{rem}
Unless we believe that a subset of on-the-line zeros of
$\zeta(s)$ is naturally distinguished, this might suggest that the
discrete spectrum is empty.
\end{rem}

\proof Let the imaginary parts of zeros be
$\ldots \le \gamma_{-1}<0<\gamma_1\le\gamma_2\le\ldots$.
Montgomery's pair correlation conjecture is that, for $0\le
\alpha<\beta$, 
\[
\begin{split}
\#\Big\{m<n: 0\le \gamma_m,\gamma_n\le T\;\hbox{with}\;
\frac{2\pi \alpha}{\log T}\le\gamma_n-\gamma_m\le \frac{2\pi\beta}{\log T}\Big\}
\\
\;\sim\;
  \int_\alpha^\beta \Big(1 - \big(\frac{\sin \pi u}{\pi
  u}\big)^2\Big)\;\d u
  \hskip80pt
    \end{split}
\]
For example, the asymptotic fraction of pairs of zeros within half the
average spacing $\frac{2\pi}{\log T}$ up to height $T$ is
\[
  \int_0^{\frac 12} \Big(1 - \big(\frac{\sin \pi u}{\pi
  u}\big)^2\Big)\;\d u
  \;\;\approx\;\; 0.11315 \;\;>\;\; 0
\]
From the lower bound in the previous section, for at least one of every
such pair $(m,n)$ the corresponding zero cannot appear among discrete
spectrum parameters $w$ for $\widetilde S_\theta$.
\qed


\section{Spacing of zeros of $\zeta_k(s)$}\label{unconditional-spacing}

{\itshape Without} assumptions on zeros $\zeta_k(s)$ as spectral
parameters for any of the pseudo-Laplacians above, we have non-trivial
corollaries on spacing of those zeros.


\subsection{Exotic eigenfunction expansions and interleaving}

Let $\theta$ continue to be a finite real-linear combination of
Eisenstein-Heegner distributions. For this section, assume that the
cut-off height $a$ is above the highest Heegner point involved in
$\theta$, so that the determinant vanishing condition
$\eta(v_{w,a})\pt\theta(u^{}_{\theta,w}) - \eta_a(u^{}_{\theta,w})^2=0$
of \eqref{resolvent}, amplified as in the computations leading up to
Corollary \ref{key-cor}, simplifies to 
\begin{align}
a^{1-w}(a^{w}+c^{}_w a^{1-w})&\times \Big(
\frac{|\theta(1)|^2}{\langle 1,1\rangle\cdot \lambda_1-\lambda_w)}
+ \frac{1}{4\pi i}\int_{(\frac 12)}
\left|\theta E_s\right|^2  
 \frac{\d s}{\lambda_s-\lambda_w}\Big)\notag\\
&= \frac{1}{2w-1} \left( a^{1-w} \theta E_w) \right)^2
\end{align}
where
\begin{equation}
\theta E_s=\sum_d \nu_d\pt
\bigg(\frac{\sqrt{|d|}}{2}\bigg)^s \pt \frac{\zeta(s)}{\zeta(2s)}\pt L(s,\chi_d)
\end{equation}
To apply the meromorphic continuation results of subsections
\ref{relatively-elementary-mero-cont} or
\ref{mero-contn-spectral-synthesis-integrals}, let  
\[
J_{\theta,w} \;=\;
\frac{\theta(1)\cdot 1}{(\lambda_1-\lambda_w)\cdot \langle 1,1\rangle}
+
\frac{1}{4\pi i}\int\limits_{(\frac 12)}
\frac{\theta E_{1-s}\cdot E_s-\theta E_{1-w} \cdot 
E_w}{\lambda_s-\lambda_w}\;\d s
\]
This allows elementary algebraic rearrangement of the
determinant-vanishing condition to a symmetrized form, as in
subsection \ref{relatively-elementary-mero-cont}:
\[
(a^{w-\frac12}+c_wa^{\frac12-w})\cdot J_{\theta,w}
\;=\;
(a^{w-\frac12} - c_w a^{\frac12-w}) \cdot
\frac{\theta E_{1-w}\cdot \theta E_w}{2(1-2w)}
\]
On the critical line $w=\frac12+i\tau$, letting $\psi(\tau)=\arg
\xi(1+2i\tau)$ as above, this is
\[
\cos\Big(\tau\log a+\psi(\tau)\Big) \cdot J_{\theta,w}
\;=\;
\sin\Big(\tau\log a+\psi(\tau)\Big)
\cdot \frac{\theta E_{1-w}\cdot \theta E_w}{-4\tau}
\]
Let $S_{\Theta,\theta}$ be the restriction of $-\Delta$ to test
functions in the Lax-Phillips space $L^2_a(\GH)$ of pseudo-cuspforms
which are also annihilated by $\theta$. Let
$\widetilde{S}_{\Theta,\theta}$ be its Friedrichs extension. We have a
variant of earlier results: 


\begin{thm}
The parameters $w=\frac12+i\tau$ for eigenvalues
$\lambda_w=w(1-w)<-1/4$ of $\widetilde{S}_{\Theta,\theta}$ are the
solutions of 
\[
\cos\Big(\tau\log a+\psi(\tau)\Big) \cdot J_{\theta,w}
\;=\;
\sin\Big(\tau\log a+\psi(\tau)\Big)
\cdot \frac{\theta E_{1-w}\cdot \theta E_w}{-4\tau}
\]
For such $w$, the eigenfunction is a linear combination
$Au_{\theta,w}+Bv_{a,w}$ of the meromorphic continuations of
$u_{\theta,w}$ and $v_{a,w}$, where $A,B$ are not both $0$, such that
\[
\theta(Au_{\theta,w}+Bv_{a,w}) \;=\; 0 \;=\; \eta_a(Au_{\theta,w}+Bv_{a,w})
\]
\end{thm}

\proof The equation $(\Delta-\lambda)u=0$ has no solution $u\in
\Eis^2\cap \Eis^0_\Theta$ except $u=0$, so the only possible
$\lambda_w$-eigenvectors are solutions $u\in \Eis^1$ to equations
\[
(\Delta-\lambda_w)u \;=\; A\cdot \theta + B\cdot \eta_a
\]
(with not both constants $0$), with the additional condition $\theta
u=0=\eta_au$. In terms of spectral expansions of elements of $\Eis^1$,
as in subsection \ref{necessary-condition}, if $\lambda_w<-1/4$ is an
eigenvalue, then $(A\theta+B\eta_a)E_w=0$. From the vector-valued
meromorphic continuation results of
subsection \ref{mero-contn-spectral-synthesis-integrals},
the meromorphic continuation of
\[
Au_{\theta,w}+Bv_{a,w}
\;=\;
\frac{A\theta+B\eta_a)(1)\cdot 1}{(\lambda_1-\lambda_w)\cdot \langle
1,1\rangle}
+\frac{1}{4\pi i}\int_{(\frac12)} \frac{(A\theta+B\eta_a)E_{1-s}\cdot
E_s}{\lambda_s-\lambda_w}\;ds
\]
at $w$ with $\Re(w)\le \frac12$ is in $\Eis^1$ when
$(A\theta+B\eta_a)E_w=0$. Further, $\eta_b(Au_{\theta,w}+Bv_{a,w})=0$
for all $b\ge a\gg_\theta 1$. Thus, when $\lambda_w<-1/4$ is an
eigenvalue, some non-zero $Au_{\theta,w}+Bv_{a,w}$ is the
eigenfunction. There exist such $A,B$, not both $0$, exactly when the
relevant two-by-two determinant vanishes, as above. That is,
$\lambda_w<-1/4$ is an eigenvalue exactly when
\[
\cos\Big(\tau\log a+\psi(\tau)\Big) \cdot J_{\theta,w}
\;=\;
\sin\Big(\tau\log a+\psi(\tau)\Big)
\cdot \frac{\theta E_{1-w}\cdot \theta E_w}{-4\tau}
\]
and the eigenfunction is a linear combination $Au_{\theta,w}+Bv_{a,w}$
of the meromorphic continuations. \qed


Via the Baire category theorem, given $a\ge 1$ and $\ve>0$, we can
increase the cut-off height $a$ by an amount less than $\ve$ so that
$\theta E_{s_j}\not=0$ for all parameters $s_j$ with $\eta_a E_{s_j}=0$,
that is, for all eigenfunctions of the form $\wedge^a E_{s_j}$ for
$\Stil_\Theta$. In the sequel, assume that such an adjustment has been
made.

We have the {\itshape interleaving result:}

\begin{thm}
Between two consecutive zeros $\tau^{(a)}_j$ and
$\tau^{(a)}_{j+1}$ of $\cos(\tau\log a+\psi(\tau))$ there is a unique
$\tau$ such that $\lambda_{\frac12+i\tau}$ is an eigenvalue of
$\Stil_{\Theta,\theta}$. 
\end{thm}

\proof  Since $\eta_a(Au_{\theta,w}+Bv_{a,w})=0$ implies
that $f_w=Au_{\theta,w}+Bv_{a,w}$ is in $\Eis^1\cap \Eis^0_\Theta$,
this $f_w$ is expressible in terms of the 
orthonormal basis of eigenfunctions $\varphi_j$ for $\Stil_\Theta$, in
an expansion converging in $\Eis^1$:
\[
f_w \;=\; \sum_i \langle f_w, \varphi_i\rangle \cdot \varphi_i
\]
Since we have arranged that $\theta E_{s_j}\not=0$ for all $s_j$,
necessarily $A\not=0$, so $A=1$ without loss of generality.
As earlier, let $\varphi_i$ have eigenvalue $s_i(1-s_i)$, with $\Im(s_i)\le
\Im(s_{i+1})$, and all $\varphi_i$ real-valued. 
Let $K$ be the kernel of $\eta_a$ on $\Eis^1\cap \Eis^0_\Theta$, and $j:K
\to \Eis^1$ the inclusion, with adjoint $j^*:\Eis^{-1}\to K^*$.
Let $\Stil^\#_\Theta=:K\longrightarrow K^*$ be as in subsection
\ref{extns-of-restrns}. Then $(\Stil_\Theta^\#-\lambda_w)f_w=j^*(\theta+B\eta_a)$.
Conveniently,
\[
j^*(\theta+B\eta_a) \;=\; j^*\theta + B\cdot j^*\eta_a
\;=\; j^*\theta + B\cdot 0
\;=\; j^*\theta
\]
The image $j^*(\theta+B\eta_a)=j^*\theta$ in $K^*$ admits an
eigenfunction expansion
\[
j^*(\theta+B\eta_a)=j^*\theta
\;=\;
\sum_i
(j^*\theta)( \varphi_i)\cdot \varphi_i
\;=\;\sum_i  
\theta(j \varphi_i)\cdot \varphi_i
\]
and then
\[
\sum_i j^*(\theta+B\eta_a) (\varphi_i)\cdot \varphi_i
\;=\;
\sum_i (j^*\theta)(\varphi_i)\cdot \varphi_i
\;=\;
j^*\theta
\;=\;
(\Stil_\Theta^\#-\lambda_w) f_w
\]
\[
\;=\;
(\Stil_\Theta^\#-\lambda_w)
\sum_i \langle f_w, \varphi_i\rangle \cdot \varphi_i
\;=\;
\sum_i \langle f_w, \varphi_i\rangle \cdot (\Stil_\Theta^\# - \lambda_w)\varphi_i
\]
\[
\;=\; 
\sum_i \langle f_w, \varphi_i\rangle \cdot (\lambda_{s_i}-\lambda_w)\varphi_i
\]
Thus, $\langle f_w,\varphi_i\rangle=(j^*\theta)
\varphi_i/(\lambda_{s_i}-\lambda_w)$, and we have an expansion
convergent in $\Eis^1$:
\[
f_w \;=\; \sum_i \frac{(j^*\theta) (\varphi_i)\cdot \varphi_i}{\lambda_{s_i}-\lambda_w}
\]
The further condition $\theta(j f_w)=0$ for $f_w$ to be an
eigenfunction for $\Stil_{\Theta,\theta}$ is
\[
0 \;=\; \theta(jf_w)\;=\; (j^*\theta)f_w
\;=\; \sum_i \frac{|(j^*\theta) \varphi_i|^2}{\lambda_{s_i}-\lambda_w}
\]
By the intermediate value theorem, there
is at least one $w$ between any two $s_i$ and $s_{j+1}$. The
derivative of $\theta f_w$ never vanishes, and no
$(j^*\theta)(\varphi_i)$ is $0$, so between $s_i$ and $s_{j+1}$ there
is {\it exactly one} $w$ on $\frac12+i\mathbb R$ satisfying the
equation. \qed 


\subsection{Dependence of eigenvalues on cut-off height}

From the previous section, the necessary and sufficient condition on
$w=\frac12+i\tau$ for $\lambda_w$ to be an eigenvalue of
$\widetilde{S}_{a,\theta}$, equivalently, of $\widetilde{S}_{\ge
a,\theta}$ is the vanishing condition
\[
\cos\Big(\tau\log a +\psi(\tau)\Big) \cdot J(w)
\;=\;
\sin \Big(\tau\log a+\psi(\tau)\Big) \cdot
\frac{\theta E_{1-w}\cdot \theta E_w}{2\tau}
\]
A simultaneous zero of $J_{\theta,w}$ and $\theta E_w$ makes
$\lambda_w$ an eigenvalue for $\widetilde{S}_\theta$ already, hence of
$\widetilde{S}_{\Theta,\theta}$ and of $\widetilde{S}_{a,\theta}$, but we expect that
there are few or no eigenvalues for $\widetilde{S}_\theta$. 
For $\lambda_w$  {\itshape not} a simultaneous zero of $J_{\theta,w}$ and
$\theta E_w$, with $w=\frac12+i\tau$, this vanishing condition is
equivalent to 
\[
\tan\Big(\tau\log a +\psi(\tau)\Big)
\;=\; \frac{2\tau\cdot J(w)}{\theta E_w\cdot \theta E_{1-w}}
\]
At least for $0<\tau\in\mathbb R$, let
\[
R(\tau\log\tau) \;=\; \frac{2\tau\cdot J(\frac12+i\tau)}{\theta
E_{\frac12+i\tau}\cdot \theta E_{\frac12-i\tau}}
\]
denote the right-hand side.

\begin{prop} For an eigenvalue $\lambda_{\frac12+i\tau}$ of $\widetilde{S}_{\ge a,\theta}$
(equivalently, of $\widetilde{S}_{a,\theta}$),
\[
\frac{\partial \tau}{\partial a} \;=\;
\frac{\displaystyle \frac{\tau}{a}(R(\tau\log \tau)^2+1)}
{(\log \tau+1)\cdot R'(\tau\log \tau) - (\log a + \psi'(\tau)) (R(\tau\log \tau)^2+1)
}
\]
\end{prop}

\begin{cor} For {\itshape all} large $\tau$,
\[
R'(\tau\log \tau)
\;\le\;
\frac{\log a + \psi'(\tau)}{\log\tau+1} \cdot  (R(\tau\log \tau)^2+1)
\]
That is, for large $\tau$, the graph of $t\to R(t\log t)$ has slope
essentially bounded above by the slope of $t\to \tan(t\log t)$ at the
same height, since $\psi'(\tau)\sim \log \tau$.
\end{cor}

\proof {\itshape (Of claim)} It makes sense to rescale the right-hand side of
the eigenvalue condition, expressing it as a function of something
asymptotic to $\tau\log\tau$, both because of the analogous scaling of
the left-hand side, and because of the vertical asymptotics of zeros
of $\theta E_w$. One simple choice is $R(\tau\log \tau)$, although one
might use $R(\psi(\tau))$ instead. Using the simple rescaling,
differentiating with respect to cut-off height $a$,    
\[
\Big(
\frac{\tau}{a}+\frac{\partial \tau}{\partial a}(\log a + \psi'(\tau))
\Big)\cdot \sec^2\Big(\tau\log a +\psi(\tau)\Big)
\;=\;
\frac{\partial \tau}{\partial a}\cdot (\log\tau+1)\cdot R'(\tau\log \tau)
\]
so
\[
\frac{\partial \tau}{\partial a}
\Big((\log\tau+1)\cdot R'(\tau\log \tau) - (\log a + \psi'(\tau)) \sec^2(\tau\log a+\psi(\tau)) \Big)
\;=\; \frac{\tau}{a}\sec^2(\tau\log a+\psi(\tau))
\]
and
\[
\frac{\partial \tau}{\partial a} \;=\;
\frac{\displaystyle \frac{\tau}{a}\sec^2(\tau\log a+\psi(\tau))}
{(\log\tau+1)\cdot R'(\tau\log \tau) - (\log a + \psi'(\tau)) \sec^2(\tau\log a+\psi(\tau))
}
\]
Since $\sec^2\alpha=\tan^2\alpha+1$, by the relation $\tan(\tau\log
a+\psi(\tau))=R(\tau\log \tau)$, at such points $\sec^2(\tau\log
a+\psi(\tau))=R(\tau\log \tau)^2+1$. Thus,
\[
\frac{\partial \tau}{\partial a} \;=\;
\frac{\displaystyle \frac{\tau}{a}(R(\tau\log \tau)^2+1)}
{(\log \tau+1)\cdot R'(\tau\log \tau) - (\log a + \psi'(\tau)) (R(\tau\log \tau)^2+1)
}
\]
which is the assertion of the claim. \qed

\proof {\itshape (Of corollary)} From the min-max principle, when
$\partial \tau/\partial a$ it exists, it satisfies
$\partial \tau/\partial a\le 0$.
We also have the {\it interleaving} of parameters $\tau$ for
eigenvalues $\lambda_{\frac12+i\tau}$ of $\widetilde{S}_{a,\theta}$ between those of
$\widetilde{S}_\Theta$.  The parameter values $\tau$ for $\widetilde{S}_\Theta$ can
be adjusted by changing cut-off height $a$: changing $a$ by about
$\pi a/\tau$ moves the zeros by about the average gap amount
$\pi/\log \tau$.
Thus, given large $\tau$, by the intermediate value theorem, we need
increase $a$ only slightly to make $\lambda_{\frac12+i\tau}$ an eigenvalue for
$\widetilde{S}_{a,\theta}$. Thus, 
$$
R'(\tau\log \tau)
\;\le\;
\frac{\log a + \psi'(\tau)}{\log\tau+1} \cdot  (R(\tau\log \tau)^2+1)
$$
holds for {\it all} large $\tau$. \qed


\subsection{Sample unconditional results on spacing of zeros}

Without any assumption that zeros of $\theta E_w$ do or do not give
eigenvalues $\lambda_w$ for a self-adjoint operator, we have
illustrative corollaries about spacing of zeros of $\theta E_w$.

\begin{cor}\label{spacing-corollary} Let $t<t'$ be large such that
$\frac12+it$ and $\frac12+it'$ are adjacent on-line zeros of $\theta
E_w$, and neither a zero of $J_{\theta,w}$.  Suppose there is a {\it
unique} zero $\frac12+i\tau_o$ of $J_{\theta,\frac12+i\tau}$ between
$\frac12+it$ and $\frac12+it'$, and $\frac{\partial}{\partial
\tau}J_{\theta,\frac12+i\tau}>0$. Then $|t'-t|\ge \frac{\pi}{\log
t}\cdot (1+O({1\over \log\log t}))$. That is, in this configuration,
the distance between consecutive zeros must be {\it at least} the
average.  \end{cor}

\begin{cor} Let $t<t'$ be large such that $\frac12+it$ and $\frac12+it'$ are
adjacent on-line zeros of $\theta E_w$, and neither a zero of $J_{\theta,w}$.  Suppose
there is no zero of $J_{\theta,w}$ on the critical line between $\frac12+it$ and
$\frac12+it'$, but there is a {\it pair} of off-line zeros
$\frac12\pm\ve+i\tau_o$ of $J_{\theta,w}$ with $t<\tau_o<t'$ and $\ve$ very
small.  Then $|t'-t|\ge \frac{\pi}{2\log t}\cdot (1+O({1\over \log\log
t}))$. That is, in this configuration, the distance between
the consecutive zeros must be {\it least} essentially half the average.
\end{cor}


\proof {\it (Of both)} In the scenario of the first corollary, this bound implies that
$R(t\log t)$ cannot get from $-\infty$ (at $\frac12+i\tau$) to $+\infty$
(at $\frac12+i\tau'$) much faster than $\tan(t\log t)$, which is by a
change of $\pi/\log t$.

In the scenario of the second corollary, $R(t\log t)$ goes to
$\pm\infty$ (with the same sign) approaching $t$ from the right and
$t'$ from the left, and is very near $0$ at $\tau_o$. Depending on the
sign, the bound by comparison to the {\it tangent} function implies
that either from $t$ to $\tau_o$, or else from $\tau_o$ to $t'$, the
function $R(t)$ cannot change faster than $\tan(t\log t)$. \qed


\begin{rem} Similar, somewhat more complicated, and perhaps less
interesting, corollaries hold, as well, for off-the-line zeros.
\end{rem}


\end{document}